\documentclass[11 pt]{amsart}
\usepackage{amsmath}
\usepackage{latexsym,amsfonts,amssymb,mathrsfs}
\usepackage{mathdots}
\usepackage{color}
\usepackage[all,2cell,color]{xy}
\usepackage{enumitem} 
\usepackage{bbm}
\usepackage{thmtools}
\usepackage{mathtools}
\usepackage{verbatim}
\usepackage{cancel} 
\usepackage{tikz}
\usetikzlibrary{cd}
\usetikzlibrary{calc}
\usetikzlibrary{math}
\usetikzlibrary{arrows}
\usetikzlibrary{fit}
\usetikzlibrary{positioning}
\usetikzlibrary{decorations.pathmorphing, arrows.meta}
\usetikzlibrary{shapes} 
\usetikzlibrary{snakes} 
\usetikzlibrary{matrix} 
\tikzset{Rightarrow/.style={double equal sign distance,>={Implies},->},
triple/.style={-,preaction={draw,Rightarrow}},
quadruple/.style={preaction={draw,Rightarrow,shorten >=0pt},shorten >=1pt,-,double,double
distance=0.2pt}}

\usepackage[T1]{fontenc}
 \usepackage[utf8]{inputenc}
 \usepackage{array}
 \usepackage{stmaryrd} 
 \usepackage{blkarray} 
 \usepackage{mathtools} 
 \usepackage{longtable}
 \usepackage{multicol} 
 \usepackage{xr-hyper}
 \usepackage{hyperref}
\usepackage[capitalise]{cleveref}

\usetikzlibrary{positioning} 
\usetikzlibrary{decorations.pathreplacing}
\usetikzlibrary{arrows, decorations.markings} 

\UseAllTwocells
\usepackage[
textwidth=3cm,
textsize=small,
colorinlistoftodos]
{todonotes}

\newcolumntype{D}{>{\hfil$}p{1.2cm}<{$\hfil}}

\tikzset{%
    symbol/.style={%
        draw=none,
        every to/.append style={%
            edge node={node [sloped, allow upside down, auto=false]{$#1$}}}
    }
}

\tikzset{%
scalearrow/.style n args={3}{
  decoration={
    markings,
    mark=at position (1-#1)/2*\pgfdecoratedpathlength
      with {\coordinate (#2);},
    mark=at position (1+#1)/2*\pgfdecoratedpathlength
      with {\coordinate (#3);},
    },
  postaction=decorate,
  }
}

\theoremstyle{plain}   
\newtheorem{thm}{Theorem}[section] 
\makeatletter\let\c@thm\c@thm\makeatother
\newtheorem{cor}{Corollary}[section]
\makeatletter\let\c@cor\c@thm\makeatother
\newtheorem{lem}{Lemma}[section]
\makeatletter\let\c@lem\c@thm\makeatother
\newtheorem{prop}{Proposition}[section]
\makeatletter\let\c@prop\c@thm\makeatother

\makeatletter\let\c@claim\c@thm\makeatother

\makeatletter\let\c@conjecture\c@thm\makeatother

\makeatletter\let\c@wconjecture\c@thm\makeatother

\makeatletter\let\c@thmRestate\c@thm\makeatother

\newtheorem*{unnumberedtheorem}{Theorem}
\newtheorem*{unnumberedcorollary}{Corollary}

\theoremstyle{definition}

\newtheorem{defn}{Definition}[section]
\makeatletter\let\c@defn\c@thm\makeatother

\makeatletter\let\c@const\c@thm\makeatother
\newtheorem{notn}{Notation}[section]
\makeatletter\let\c@notn\c@thm\makeatother

\makeatletter\let\c@convention\c@thm\makeatother

\makeatletter\let\c@convention\c@thm\makeatother

\theoremstyle{remark}

\newtheorem{rmk}{Remark}[section]
\makeatletter\let\c@rmk\c@thm\makeatother
\newtheorem{ex}{Example}[section]
\makeatletter\let\c@ex\c@thm\makeatother

\makeatletter\let\c@observation\c@thm\makeatother

\makeatletter\let\c@warning\c@thm\makeatother
\newtheorem{digression}{Digression}[section]
\makeatletter\let\c@digression\c@thm\makeatother

\makeatletter\let\c@answ\c@thm\makeatother

\makeatletter\let\c@answ\c@thm\makeatother

\makeatletter\let\c@aside\c@thm\makeatother

\makeatletter
\let\c@equation\c@thm
\numberwithin{equation}{section}
\makeatother

\crefname{lem}{Lemma}{Lemmas}
\crefname{thm}{Theorem}{Theorems}
\crefname{defn}{Definition}{Definitions}
\crefname{notn}{Notation}{Notations}
\crefname{const}{Construction}{Constructions}
\crefname{prop}{Proposition}{Propositions}
\crefname{rmk}{Remark}{Remarks}
\crefname{cor}{Corollary}{Corollaries}
\crefname{equation}{Display}{Displays}
\crefname{ex}{Example}{Examples}
\crefname{thmalph}{Theorem}{Theorems}
\crefname{answ}{Answer}{Answers}
\crefname{question}{Question}{Questions}

\newcommand{\cC}{\mathcal{C}}
\newcommand{\cD}{\mathcal{D}}

\newcommand{\cO}{\mathcal{O}}

\newcommand{\cS}{\mathcal{S}}

\newcommand{\cat}{\cC\!\mathit{at}}
\newcommand{\set}{\cS\!\mathit{et}}
\newcommand{\sset}{\mathit{s}\set}

\newcommand{\adch}{\mathit{ad}\mathcal C\mathit{h}}

\newcommand{\msset}{m\sset}

 \newcommand{\omegacat}{\omega\cat}

 \newcommand{\comp}{\ast}

\newcommand{\inttrunc}[1]{\tau_{\leq n}^{\text{i}}}

\DeclareMathOperator{\id}{id}
\DeclareMathOperator{\im}{im}
\newcommand{\aamalg}[1]{\underset{#1}{{\amalg}}} 
\DeclareMathOperator{\op}{op}

\DeclareMathOperator{\pr}{pr}

\newcommand{\NRS}{N^{\operatorname{RS}}}

\newcommand{\atom}[1]{\langle#1\rangle}

\newcommand{\tabld}[2]{\begin{pmatrix}#1^-_0 &\dots &#1^-_{#2-1}
  &#1^-_{#2}\cr\noalign{\vskip 3pt} #1^+_0 &\dots &#1^+_{#2-1}
  &#1^+_{#2}\end{pmatrix}}

\newcommand{\supp}{\mathrm{supp}}

\newcommand{\dwedgevar}[2]{\prescript{#1}{#2}{\wedge}}

\input{square.tex}

\author{Andrea Gagna}
\address{Institute of Mathematics, Czech Academy of Sciences\\ \v{Z}itn\'a 25 \\115 67   Praha 1\\ Czech Republic}
\email{gagna@math.cas.cz}

\author{Viktoriya Ozornova}
\address{Max Planck Institute for Mathematics, Bonn, Germany}
\email{viktoriya.ozornova@mpim-bonn.mpg.de}

\author{Martina Rovelli}
\address{Department of Mathematics and Statistics, 
University of Massachusetts, MA 01003-9305
Amherst, USA
}
\email{rovelli@math.umass.edu} 

\title[Nerves and joins of ``free loop-free'' $\omega$-categories]{Nerves and cones of\\ ``free loop-free'' $\omega$-categories}

\begin{document}

\maketitle

\begin{abstract}
We show that the complicial nerve construction is homotopically compatible with two flavors of cone constructions when starting with an $\omega$-category that is suitably free and loop-free. An instance of the result recovers the fact that the standard $m$-simplex is equivalent to the complicial nerve of the $m$-oriental.
\end{abstract}

 \tableofcontents

\newpage 
\section*{Introduction}

Several flavors of higher categories have proven to be the mathematical structure governing interesting phenomena in geometry and theoretical physics. One of the most general setups is to study higher categories which consist of objects and morphisms in each dimension with a composition along lower dimensional morphisms that satisfies some version of the usual axioms.

Most examples occurring in nature are \emph{weak} $\omega$-categories, in which the axioms -- such as the associativity axioms -- are only satisfied up to a higher invertible cell. However, in order to develop a sound theory it is crucial to be able to recover also the examples of \emph{strict} higher categories, referred to as ``strict $\omega$-categories'' or just $\omega$-categories, in which the axioms are satisfied on the nose. On the one hand, strict $\omega$-categories classify most relevant features of weak $\omega$-categories; on the other hand, it is extremely complicated to formalize the combinatorics of weak $\omega$-categories, and the strict framework provides a first highly non-trivial playground.

The realization of a strict higher category $\cD$ as a weak higher category is often implemented by some kind of nerve construction $N\cD$. Efforts towards defining homotopical nerve constructions include work by different groups of authors and for different models and flavors of weak higher categories, e.g. \cite{VerityComplicialAMS,Nerves2Cat,CampbellHoCoherent,GHL,MoserNerve,LoubatonNerfs}.

Despite the potentially different frameworks, those nerve constructions are usually right adjoint functors by nature, and are therefore well-behaved with respect to constructions involving limits and right adjoints. Instead, they in general have an ill behaviour with respect to constructions involving colimits or left adjoint functors.

It is notably a hard problem to understand which colimit-related constructions for strict $\omega$-categories commute with the nerve construction at least up to homotopy, even when restricting one's attention to strict $n$-categories for low values of $n$. A few results in this direction, which treat primarily the compatibility of nerves of $1$- and $2$-categories with certain pushouts, include \cite{ThomasonModelCat,joyal,ORfundamentalpushouts,HORR}. But, otherwise, very little is currently known for $n$-categories and more generally $\omega$-categories.

In this paper we address the compatibility of the complicial nerve of $\omega$-categories with two types of cone constructions, which can be understood as specific types of colimit constructions, and play a role in detecting oplax limits and lax colimits
of $\cD$-valued diagrams, respectively.

Given a weak or strict $\omega$-category $\cD$, the cone constructions are obtained by adding an extra point that maps to or from every object of $\cD$, and by filling in higher morphisms that point towards the composites, as is illustrated in the picture below for the small case of $\cD=[1]$ being the walking morphism category.\footnote{There are two other types of cone constructions, in which the $2$-cells are reversed (see \cite[Remark~6.37]{AraMaltsiniotisJoin} and \cite[\textsection 4.2]{GHL3} for more details on how to construct them in $\omega$-categories and in scaled simplicial sets), that would detect lax limits and oplax colimits, but they are not as easy to construct and the techniques from this paper cannot be employed on the nose to cover those cases. 
}
\[
       \cD\star[0]=
        \fbox{
        \begin{tikzcd}[row sep=tiny, column sep=small,ampersand replacement=\&]
            x \ar[rd, ""{swap, name=s}] \ar[dd] \&       \\
                                                \& \top     \\
            y \ar[ur]
                \ar[Rightarrow, from=s, to=3-1, shorten >= 4pt]
        \end{tikzcd}
        }
          \quad\quad
        \cD=
        \fbox{
        \begin{tikzcd}[row sep=normal,ampersand replacement=\&]
            x \ar[d] \\ y
        \end{tikzcd}
        }
        \quad\quad
         [0]\star\cD=
        \fbox{
        \begin{tikzcd}[row sep=tiny, column sep=small,ampersand replacement=\&]
                                        \& x \ar[dd] \\
        \bot \ar[ur] \ar[dr, ""{name=s}]   \&           \\
                                        \& y         
                \ar[Rightarrow, from=s, to=1-2, shorten >= 4pt]
        \end{tikzcd}
        }
\]

More precisely, the two cone constructions, $\cD\star[0]$ and $[0]\star\cD$, arise an instance of a join construction $\star$, which can be understood as a functor cocontinuous  in each variable on connected diagrams. This was studied by Ara--Maltsiniotis \cite{AraMaltsiniotisJoin} in the context of strict $\omega$-categories, and by Verity \cite{VerityComplicialI} in the context of weak $\omega$-categories modelled by complicial sets.

As our main result, we prove that these two cone constructions are compatible with the nerve construction for particularly nice $\omega$-categories.
The precise niceness condition
on the $\omega$-category $\cD$ is that of being a ``strong Steiner $\omega$-category'', which is designed to encode the fact that $\cD$ is freely generated and does not contain any loop.

\begin{unnumberedtheorem}[{\cref{maintheorem,maintheoremdual}}]
If $\cD$ is a suitably free and loop-free $\omega$-category, there are equivalences of weak $\omega$-categories\footnote{The second statement cannot be formally obtained from the first one, the obstruction being that the ``op-dual'' $\cD^{\op}$ of a strong Steiner $\omega$-category is not generally strong Steiner. See~\cref{dualnotformal} for more details on this.}
\[N(\cD\star[0])\simeq N\cD\star N[0]\quad\text{and}\quad N([0]\star\cD)\simeq N[0]\star N\cD.\]
\end{unnumberedtheorem}

 Examples of interest include $\cD$ being the $n$-oriental $\cO[n]$, the $n$-disk $\cD[n]$, and more generally all objects of Joyal's \cite{JoyalDisks} cell category $\Theta_n$ for $n\ge0$. Furthermore, the class of strong Steiner $\omega$-categories is closed under join and lax Gray tensor product of $\omega$-categories \cite{AraMaltsiniotisJoin}. Originally, strong Steiner $\omega$-categories are defined using algebraic models \cite{SteinerEmbedding,AraMaltsiniotisJoin}, and an explicit characterization of strong Steiner $\omega$-categories in terms of their categorical structure is the subject of \cite{GOR2}.

The focus on this particular class of $\omega$-categories is due to the fact that we make use of Steiner's machinery from \cite{SteinerUniversal} and Ara--Maltsiniotis' join construction from \cite{AraMaltsiniotisJoin}. These tools allow one to recognize $m$-simplices of the nerve of $\cD\star\overline\cD$ as suitable chain maps $O[m]\to D\star\overline D$, and use the rich structure of the category of (augmented directed) chain complexes. Here, $O[m]$ is the normalized chain complex of the standard $m$-simplex, $\star$ is the join of (augmented directed) chain complexes and $D$ is Steiner's algebraic model of the $\omega$-category $\cD$.

The theorem is proven in the model of weak $\omega$-categories given by complicial sets,\footnote{The theorem is proven in the model structure for (non-saturated) complicial sets, and holds \emph{a fortiori} in any localization of it, in particular for the model structure of saturated complicial sets.} which was intensively developed by Verity \cite{VerityComplicialI}, and the nerve is the one introduced by Roberts and further studied by Street \cite{StreetFillers} and Verity \cite{VerityComplicialAMS}. In this model, a weak $\omega$-category consists of a simplicial set in which the simplices that should be thought of as identity cells or invertible cells are marked. Unlike in other models of weak higher categories,
it is easy to implement explicitly the two desired cone constructions.

The theorem is valuable on at least two fronts. On the one hand, it is one of the few theorems asserting compatibility of a construction between strict $\omega$-categories and weak $\omega$-categories, in the form of complicial sets. On the other hand, even when $n=1$, the theorem is a crucial ingredient to understand whether lax limits of diagrams of $1$-categories can be equivalently computed in the $2$-category or in the $(\infty,2)$-categories of $1$-categories. Even more generally, it helps treat the same question for colimits in an $n$-category and colimits computed in the $(\infty,n)$-category given by its nerve.

Furthermore, the theorem provides an explicit fibrant replacement of a cone $\Delta[0]\star N\cD$.
As a special case, when $\cD=\cO[n-1]$ is an oriental $\omega$-category, the theorem also specializes to the following corollary.

\begin{unnumberedcorollary}[{\cref{orientalsimplex}}]
For $n\ge0$, there is an equivalence of weak $\omega$-categories
\[\Delta[n]\simeq N\cO[n].\]
\end{unnumberedcorollary}
This case was treated independently -- with different techniques that do not apply to the general case -- by Maehara and appeared in \cite{MaeharaOrientals} during the preparation of this article.

 \addtocontents{toc}{\protect\setcounter{tocdepth}{1}}
\subsection*{Acknowledgements}
We are thankful to L\'eonard Guetta for clarifying the relation between the different notions of freeness for $\omega$-categories, and to Georges Maltsiniotis for pointing out that strong Steiner categories are not closed under op-duality. The first-named author gratefully acknowledges the support of Praemium Academiae of M.~Markl and RVO:67985840.

 \addtocontents{toc}{\protect\setcounter{tocdepth}{2}}

\section{Strict and weak \texorpdfstring{$\omega$}{w}-categories}

In this section we briefly recap for the reader the notion of a strict $\omega$-category, with a particular focus on strong Steiner categories, and present a model of weak $\omega$-categories in terms of complicial sets.

\subsection{Strict \texorpdfstring{$\omega$}{w}-categories}

While we refer the reader to e.g.~\cite{StreetOrientedSimplexes} for a traditional approach to the definition of an $\omega$-category, we briefly recall the main features here.

The data of an \emph{$\omega$-category} $\cC$ consists of a collection of sets $\cC_q$, ${q \geq 0}$,
where $\cC_0$ is called the set of \emph{objects} of $\cC$ and $\cC_q$, $q>0$,
is the set of \emph{$q$-cells} or \emph{$q$-arrows} or cells of \emph{dimension} $q$ of $\cC$, together with:
\begin{itemize}[leftmargin=*]
    \item \emph{source} and \emph{target} operators $s_q, t_q \colon \cC_p \to \cC_q$
    for all $p > q \geq 0$;
    \item \emph{identity} operators $\id_q \colon \cC_p \to \cC_{q}$ for all $q> p\ge0$; 
    \item \emph{composition} operators $\comp_p \colon \cC_q \times_{\cC_p} \cC_q \to \cC_q$
    defined for all $q > p \geq 0$ and all pairs of $q$-cells $(g, f)$ such that $s_p(g) = t_p(f)$.
\end{itemize}
Notice that this is equivalent to endow every pair $(\cC_p, \cC_q)$, $q > p \geq 0$,
with the structure (but not the axioms, yet) of a category.
For all $r > q > p \geq 0$, we ask that the triple $(\cC_p, \cC_q, \cC_r)$
together with all the relevant source, target, identity and composition operators
is a $2$-category.

An \emph{$\omega$-functor} $F \colon \cC \to \cC'$ between $\omega$-categories $\cC$ and $\cC'$ is a collection of maps $F_q \colon \cC_q \to \cC'_q$ for
$q \geq 0$ that preserves source, target, identity, and composition.
We denote by $\omegacat$ the category of (small) $\omega$-categories and $\omega$-functors.

A cell in an $\omega$-category $\cC$ which is the identity of a lower dimensional cell is said to be \emph{trivial}. An $\omega$-category in which all $q$-cells are trivial for $q>n$ is an $n$-category, and an $\omega$-functor between $n$-categories reduces to an $n$-functor. We denote by $n\cat$ the (full) subcategory of $\omega\cat$ given by $n$-categories and $n$-functors.

\begin{rmk}
\label{truncation}
The canonical inclusion $n\cat\hookrightarrow \omega\cat$ for $n\ge0$ admits a right adjoint $\varpi_{n}\colon \omega\cat\rightarrow n\cat$,
called the $n$-truncation functor and treated e.g.~in~\cite[\textsection 1.2]{AraMaltsiniotisJoin}, which produces an $n$-category $\varpi_{n}\cC$ by forgetting all non-trivial $q$-cells of $\cC$ for $q> n$
and does not change the underlying $n$-category of $\cC$.
\end{rmk}

    Ara and Maltsiniotis \cite{AraMaltsiniotisJoin} introduced a join monoidal structure
    for $\omega$-categories, whose monoidal unit is the empty $\omega$-category.
    More details on how to construct the join $\cD\star\overline{\cD}$ of two $\omega$-categories $\cD$ and $\overline{\cD}$ are postponed until~\cref{sec:nu}. For now, we mention the following instances of $\omega$-categories obtained as an iterated join construction: Street's orientals from \cite{StreetOrientedSimplexes}.

\begin{ex}
\label{ex:oriental}
For $n\ge0$, the $n$-th oriental $\cO[n]$ is an $n$-category. The first few orientals can be depicted as follows.
      \[
			\cO[0] =
			\fbox
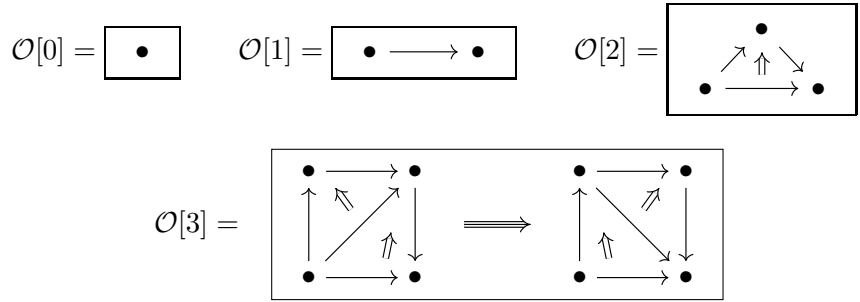
{
			\begin{tikzcd}[ampersand replacement=\&]
				\bullet
			\end{tikzcd}
			}
\quad\quad
			\cO[1] =
		\fbox{	\begin{tikzcd}[ampersand replacement=\&]
				\bullet \arrow[r] \& \bullet
			\end{tikzcd}
			}
		\quad\quad
		   	\cO[2] = 
				\fbox{
			 \begin{tikzcd}[row sep=small, column sep=tiny,ampersand replacement=\&]
                                            \& \bullet \ar[rd]   \&         \\
        \bullet \ar[ru] \ar[rr, ""{name=s}] \&                   \& \bullet
                \ar[Rightarrow, from=s, to=1-2, shorten >= 2pt]
        \end{tikzcd}
}
		\]
\begin{center}
    \centering
    \begin{tikzpicture}
    \draw (-3.0,-1) rectangle (3.0,1);
			\square{
				/square/label/.cd,
	     			0=$\bullet$, 1=$\bullet$, 2=$\bullet$, 3=$\bullet$,
	     			01={}, 12={}, 23={}, 03={},
	     			012={}, 023={}, 123={}, 013={},
	     			0123={}
     			}
	\end{tikzpicture}
\end{center}
For a precise account, we refer the reader to Street's original construction of $\cO[n]$ from \cite{StreetOrientedSimplexes}. Alternatively, following {\cite[\textsection7]{AraMaltsiniotisJoin}} or
\cite{BuckleyGarnerOrientals}, the $n$-oriental can be understood as the $\omega$-category obtained as the
iterated join of $n+1$ copies of the terminal ($\omega$-)category $[0]$ 
with itself, that is
\[\cO[n]\cong\underbrace{[0]\star\dots\star[0]}_{n+1},\]
which happens to be an $n$-category.
\end{ex}

The $n$-th oriental $\cO[n]$ should be thought as ``the free $n$-category over an $n$-simplex'' and a consequence of the main theorem,~\cref{orientalsimplex}, will formalize the fact that it is also ``the free $(\infty, n)$-category over an $n$-simplex''.

For $n\ge0$, another important family of examples of $n$-categories is the collection of objects of Joyal's cell category $\Theta_n\subseteq n\cat$ from \cite{JoyalDisks}.

\subsection{Complicial sets and weak \texorpdfstring{$\omega$}{w}-categories}

In this paper, we  consider a model of weak higher categories due to Verity based on the following mathematical object, considered by Verity in \cite[\textsection 5.1]{VerityComplicialAMS} under the name of ``stratified simplicial set''.

\begin{defn}
A \emph{simplicial set with marking}
is a simplicial set endowed with a subset of simplices of strictly positive dimensions -- called \emph{thin} or \emph{marked} -- that contains all degenerate simplices.
\end{defn}

We denote by $m\sset$ the category of simplicial sets with marking and marking-preserving simplicial maps.

The category of simplicial sets with marking hosts a model structure (in fact more a priori than one) for the homotopy theory of $(\infty,n)$-categories.

\begin{thm}
\label{modelstructurewithsaturation}
Let $n\in\mathbb N\cup\{\infty\}$.
\label{modelstructureondiscretepresheaves}
The category $m\sset$ supports a
model structure where the fibrant objects are precisely the \emph{$(\infty,n)$-categories} and the cofibrations are precisely the monomorphisms (of underlying simplicial sets).
\end{thm}

The definition of an $(\infty,n)$-category will not be needed in this paper, but we nevertheless provide a little context for the interested reader.

For $n\in\mathbb N$, one can define an \emph{$(\infty,n)$-category} to be a simplicial set with marking
\begin{enumerate}[leftmargin=*]
    \item that is a weak complicial set, as in \cite[Definition~15]{VerityComplicialI},
        \item in which all simplices are marked in dimension higher than $n$, and
        \item that satisfies a saturation condition up to dimension $n$, as in \cite[\textsection 3.2]{EmilyNotes}.
\end{enumerate}
In technical terms, this is referred to as an ``$n$-trivial saturated weak complicial set'', or an \emph{$n$-complicial set} for short.
 The interpretation is that, in a saturated $n$-complicial set, $k$-simplices represent $k$-morphisms with a specified boundary decomposition, and according to this correspondence the marked $k$-simplices
 are precisely the $k$-equivalences. We refer the reader e.g.~to \cite{EmilyNotes} for further elaboration on this viewpoint. The model structure for $n\in\mathbb N$ is the one from \cite[Theorem~1.28]{or}.

The case $n=\infty$ is subtle, as there are a few (a priori non-equivalent!) ways to make sense of an $(\infty,\infty)$-category, each of them coming with a supporting (different!) model structure. Comparing with the definition of an $(\infty,n)$-category for finite $n$, an \emph{$(\infty,\infty)$-category} consists of
\begin{enumerate}[leftmargin=*]
    \item a weak complicial set
    \item without any analog
    of the second condition, and
    \item potentially with an analog of the third condition.
\end{enumerate}
Three popular choices for how to replace the third condition are: to not include any, to request the saturation condition for every dimension $n$, or to request a stronger saturation condition that strictly generalizes all finite dimensional ones at once. The three choices lead respectively to Verity's model of $(\infty,\infty)$-categories from \cite{VerityComplicialI}, to an \emph{inductive} model of $(\infty,\infty)$-categories, and to a \emph{coinductive} model of $(\infty,\infty)$-categories, using the terminology from \cite{jf}. In the first case, the supporting model structure is the one from \cite[Example~104]{VerityComplicialI}, and in the second and third case the model structures are further (left Bousfield) localizations of this one.

We recall for the reader the terminology that will be needed throughout the paper.

\begin{defn}[{\cite[Notation~101]{VerityComplicialAMS}}]
A sub-simplicial set with marking $X$ of a simplicial set with marking $Y$ is \emph{regular} if a simplex of $X$ is marked in $X$ precisely when it is marked in $Y$.
\end{defn}

\begin{notn}[{\cite[Notation~105, Definition~120]{VerityComplicialAMS}}]
\label{preliminarynotation}
We denote
\begin{itemize}[leftmargin=*]
    \item by $\Delta^k[m]$, for $0\leq k \leq m$, the standard $m$-simplex in which a non-degenerate simplex is marked if and only if it contains the vertices $\{k-1,k,k+1\}\cap [m]$;
    \item by $\Delta^k[m]'$, for $0\leq k \leq m$, the standard $m$-simplex with marking obtained from $\Delta^k[m]$ by additionally marking the $(k-1)$-st and $(k+1)$-st face of $\Delta[m]$;
    \item by $\Delta^k[m]''$, for $0\leq k \leq m$, the standard $m$-simplex with marking obtained from $\Delta^k[m]'$ by additionally marking the $k$-th face of $\Delta[m]$;
    \item by $\Lambda^k[m]$, for $0\leq k \leq m$, the regular sub-simplicial set of $\Delta^k[m]$ with marking
    whose simplicial set is the $k$-horn $\Lambda^k[m]$.
\end{itemize}
\end{notn}

We fix the following terminology (cf.~\cite[Definition~15]{VerityComplicialI}).

\begin{defn}
A map of simplicial sets with marking $X\to Y$ is a \emph{complicial inner anodyne extension} if it can be written as a retract of a transfinite composition of pushouts of maps of the following form:
\begin{enumerate}[leftmargin=*]
 \item  complicial inner horn extensions
 $$\Lambda^k[m]\to \Delta^k[m]\text{ for $m> 1$ and $0< k< m$},$$
 \item complicial thinness extensions
$$\Delta^k[m]' \to \Delta^k[m]''\text{ for $m\geq 2$ and $0< k < m$}.$$
\end{enumerate}
\end{defn}

We will produce several 
acyclic cofibrations using the following one.

\begin{lem}[{\cite[Lemma~1.12]{ORfundamentalpushouts}}]
  \label{CompMarkAtOnce} For $m\geq 2$ and $0<k<m$, let $\Lambda^k[m]'$ denote the regular subset of $\Delta^k[m]'$ whose underlying simplicial set is given by the $k$-horn $\Lambda^k[m]$. The inclusion
    \[\Lambda^k[m]'\to \Delta^k[m]''\text{ for $m\geq 2$, $0<k<m$}\]
is a complicial inner anodyne extension.
\end{lem}

\begin{rmk}
 \label{underlyingcomplicialinner}
One can prove with standard model categorical techniques the following formal properties of complicial inner anodyne extensions.
\begin{enumerate}[leftmargin=*]
\item Any complicial inner anodyne extension is an acyclic cofibration in the model structure(s) from~\cref{modelstructurewithsaturation}.
\item The underlying simplicial map of a complicial inner anodyne extension is an inner anodyne extension of simplicial sets.
     \item The class of complicial inner anodyne extensions is closed under transfinite composition and pushouts.
\end{enumerate}
\end{rmk}

\subsection{Join of simplicial sets with marking}

Given simplicial sets $X$ and $\overline X$, the \emph{join} is the 
simplicial set $X\star \overline X$ with set of $q$-simplices
given by
\[(X\star \overline X)_q=\coprod_{\substack{ k+ 1 + \overline k = q \\ k \geq -1, \overline k \geq -1}}X_k\times \overline X_{\overline k},\]
assuming that in degree $-1$ we have the singletons $X_{-1}=\overline X_{-1}=\{\emptyset\}$,
and simplicial structure for which faces and degeneracies of a simplex $(\sigma, \tau)\in X_{k}\times \overline X_{\overline k}\subseteq (X\star \overline X)_q$ are given by
\[
d_i(\sigma, \overline\sigma)=\left\{\begin{array}{ll}
(d_i\sigma, \overline\sigma) &\mbox{ if }0\leq i\leq k,\\
(\sigma, d_{i-k-1}\overline\sigma) &\mbox{ if }k+1\leq i \leq q=k+1+\overline k,
\end{array}
\right.
\]
and 
\[
s_i(\sigma,\overline\sigma)=
\left\{\begin{array}{ll}
(s_i\sigma,\overline\sigma) & \mbox{ if }0\leq i\leq k,\\
(\sigma, s_{i-k-1}\overline\sigma)& \mbox{ if }k+1\leq i \leq q=k+1+\overline k.
\end{array}
\right.
\]

We will care about the situation when $\overline X=\Delta[\top]$, resp.~$X=\Delta[\bot]$, is a simplicial set no non-degenerate simplices and a unique $0$-simplex denoted $\top$, resp.~$\bot$. In this case,
the non-degenerate simplices of $X\star\Delta[\top]$, resp.~$\Delta[\bot]\star X$, are
either of the form $(\sigma,\emptyset)$ or $(\sigma,\top)$, resp.~of the form $(\emptyset,\overline\sigma)$ or $(\bot,\overline\sigma)$.

\begin{defn}[{\cite[\textsection 3.1]{VerityComplicialI}}]
Given simplicial sets with marking $X$ and $\overline X$, the \emph{join} is the simplicial set $X\star \overline X$ endowed with the marking in which a simplex $(\sigma,\overline\sigma)$ is marked in $X\star \overline X$ if and only if either $\sigma$ is marked in $X$ or $\overline\sigma$ is marked in $\overline X$ (or both).
\end{defn}

\subsection{Complicial nerve of \texorpdfstring{$\omega$}{w}-categories}

The geometry of orientals is such that the construction $m\mapsto \cO[m]$ defines a cosimplicial object $\cO[\bullet]$ in $\omega\cat$, and in particular it makes sense to define the following nerve construction $N\colon\omega\cat\to\sset$ originally due to Street {\cite{StreetOrientedSimplexes}}.
The \emph{Street nerve} $N\cD$ of an $\omega$-category $\cD$ is the simplicial set in which the set of $m$-simplices is the set of $\omega$-functors $\mathcal O[m]\to \cD$, namely
\[N_m\cD=\omega\cat(\cO[m],\cD),\]
and the simplicial structure is induced by the geometry of orientals.

The Street nerve can be endowed with the following marking, originally considered by Roberts in unpublished work and Street in \cite{StreetOrientedSimplexes}, further studied by Verity in \cite{VerityComplicialAMS}, and later discussed by Riehl in \cite{EmilyNotes}, obtaining a functor $\NRS\colon\omega\cat\to\msset$.

\begin{defn}
\label{RSnerve}
Let $\cD$ be an $\omega$-category.
The \emph{Roberts--Street nerve} is the simplicial set with marking $\NRS\cD$, in which the underlying simplicial set is the Street nerve $N\cD$, and an $m$-simplex of $N\cD$ is marked in $\NRS\cD$ if and only if the corresponding $\omega$-functor $\cO[m]\to\cD$ sends the unique non-trivial $m$-cell of $\cO[m]$ to a trivial $m$-cell of $\cD$.
\end{defn}

 In particular, in the Street nerve of an $n$-category $\cD$ all simplices in dimension at least $n+1$ are marked.
 
 \begin{rmk}
 Let $n\in\mathbb N\cup\{\infty\}$.
As a consequence of work by Loubaton \cite{LoubatonNerfs}, for any $n$-category $\cD$ that does not contain any equivalences (in the sense of~\cite[\S 4.2]{LMW}) the Roberts--Street nerve $\NRS\cD$ is fibrant in the model structure(s) from~\cref{modelstructurewithsaturation}, and is therefore an $(\infty,n)$-category.\footnote{When starting with an arbitrary $\omega$-category $\cD$, the Roberts--Street nerve is not generally an $(\infty,n)$-category, but can be made into one by modifying suitably the marking, which is also treated in Loubaton's work.}
 \end{rmk}

\begin{rmk}
\label{comparisonmap}
Given $\omega$-categories $\cD$ and $\overline{\cD}$, there is a canonical inclusion of simplicial sets with marking
\[
\NRS\cD\star\NRS\overline{\cD}\hookrightarrow \NRS(\cD\star\overline{\cD})
\]
given by
regarding a pair of $\omega$-functors $\alpha\colon\cO[k]\to\cD$ and $\overline\alpha\colon\cO[\overline k]\to\overline{\cD}$ for $k+1+\overline k=m$ as a functor
\[
\begin{tikzcd}
\cO[m]\cong\cO[k]\star\cO[\overline k]\xrightarrow{\alpha\star\overline\alpha}\cD\star\overline{\cD}.
\end{tikzcd}
\]
This inclusion is natural in $\cD$ and $\overline{\cD}$.
When $\overline\cD=[0]$, resp.~$\cD=[0]$, the construction specializes to inclusion of simplicial sets with marking
\begin{equation}
\label{eqcomparisonmap1}
\NRS\cD\star\Delta[0]\hookrightarrow \NRS(\cD\star[0]),\ \text{resp.}\ \Delta[0]\star\NRS\overline{\cD}\hookrightarrow \NRS([0]\star\overline{\cD}).
\end{equation}
Iterating either of the constructions above, one also obtains an inclusion of simplicial sets with marking
\begin{equation}
\label{eqcomparisonmap2}\Delta[m]\cong\underbrace{\Delta[0]\star\dots\star\Delta[0]}_{m+1}\hookrightarrow \NRS(\underbrace{[0]\star\dots\star[0]}_{m+1})\cong\NRS\cO[m].\end{equation}
\end{rmk}

The main goal of this paper is to show that inclusions from (\ref{eqcomparisonmap1}) are weak equivalences when $\cD$, resp.~$\overline\cD$, is a strong Steiner
$\omega$-category, and that the inclusion from (\ref{eqcomparisonmap2}) is a weak equivalence for every $m\ge0$.

To this end, we will take a detour to obtain a handy algebraic description of the simplicial sets $N(\cD\star\overline\cD)$ and $N\cO[m]$, generalizing and recovering ideas from \cite{SteinerEmbedding,SteinerUniversal}.

\section{Steiner's algebraic models}
\label{sec:Steiner}

We briefly recall the main notation and results from Steiner's theory that will be needed in the rest of the paper.

\subsection{Augmented directed chain complexes}
By \emph{chain complex} $C$ we will always mean an
$\mathbb{N}$-graded chain complex of abelian
groups with homological indexing, that is
a family $(C_q)_{q\geq 0}$ of abelian groups $C_q$ for $q\ge0$, together with homomorphisms
$\partial_{q} \colon C_{q+1} \to C_{q}$ satisfying $\partial_q \partial_{q+1}=0$.

Given chain complexes $C$ and $\overline C$, a \emph{chain map} or \emph{morphism of chain complexes} $\phi\colon C\to \overline C$ consists of a family of homomorphisms $(\phi_q\colon C_q\to \overline C_q)_{q\geq 0}$ such that for every $q\ge0$ the homomorphism $\phi_q$ commutes with the differentials, namely $\overline\partial_{q} \phi_{q+1}=\phi_q\partial_{q}$ for every $q \geq 0$.

An \emph{augmented chain complex} is a pair $(C,\epsilon)$
of a chain complex $C$ and an augmentation, namely a homomorphism $\epsilon\colon C_0\to\mathbb Z$ such that $\epsilon\partial_0=0$.

An augmented chain map $\phi \colon (C,\epsilon) \to (\overline C, \overline\epsilon)$ between augmented chain complexes $(C,\epsilon)$ and $(\overline C,\overline\epsilon)$ consists of a chain map $\phi\colon C\to\overline C$ that is moreover compatible with the augmentations, namely such that $\overline\epsilon \phi_0 = \epsilon$.

\begin{defn}[{\cite[Definition~2.2]{SteinerEmbedding}}]
 An \emph{augmented directed complex}
 is a triple\footnote{We will mostly denote an augmented directed chain complex
by its underlying chain complex, if the augmentation
and the positivity submonoids are clear from the context.}
 $(C, C^+, \epsilon)$ where $(C, \epsilon)$ is
 an augmented chain complex and $C^+ =  (C^+_q)_{q\ge0}$ is a collection of commutative monoids,
 where $C^+_q$ is a submonoid of $C_q$ called the
 \emph{positivity submonoid} of $C_q$.
 
 A \emph{morphism of augmented directed chain complexes}, or an \emph{augmented directed chain map} $\phi \colon (C, C^+, \epsilon) \to (\overline C, \overline C^+, \overline\epsilon)$ between augmented directed chain complexes $(C, C^+, \epsilon)$ and $(\overline C, \overline C^+, \overline\epsilon)$
 is an augmented chain map~$\phi \colon (C, \epsilon) \to (\overline C, \epsilon')$ that moreover respects the positivity submonoids\footnote{The differentials of an augmented directed chain complex
    \emph{need not} respect the positivity submonoid.}, namely such that $\phi_q(C^+_q) \subseteq \overline C^+_q$
 for all $q\ge0$.
\end{defn}

We denote by $\adch$ the category of augmented directed chain complexes and augmented directed chain maps.

We now introduce an important example of augmented directed chain complexes, based on standard simplices.

\begin{notn} Let $n\ge-1$ and $q\ge0$. We denote
\begin{itemize}[leftmargin=*]
    \item by $\Delta[n]_q=\Delta([q],[n])$ the set of $q$-simplices of the standard simplex\footnote{We follow the convention that $[-1]$ is the empty category, and $\Delta[-1]$ is the initial simplicial set, which is levelwise empty.} $\Delta[n]$. A generic $q$-simplex in $\Delta[n]$
    is of the form
\[[\mathbf a]=[a_0,a_1,\dots,a_q]\ ,\quad0\leq a_0 \le a_1 \le \ldots \le a_q\leq n.\]
\item by $B[n]_q\subseteq\Delta[n]_q $ the set of non-degenerate $q$-simplices of $\Delta[n]$. A generic $q$-simplex in $\Delta[n]$ is of the form
\[[\mathbf a]=[a_0,a_1,\dots,a_q]\ ,\quad0\leq a_0 < a_1 < \ldots <a_q\leq n.\]
\item by $O[n]_q=\mathbb Z[B[n]_q]$ the abelian group freely generated by non-degenerate $q$-simplices of $\Delta[n]$. The generic element of $O[n]_q$ is a formal sum
\[
        c = \sum_{[\mathbf a] \in B[n]_q} c_{[\mathbf a]} \cdot [\mathbf a]\ ,
        \quad c_{[\mathbf a]} \in \mathbb{Z}
    \]
    where finitely many coefficients $c_{[\mathbf a]}$
    are non-zero.
\item by $O[n]_q^+=\mathbb N[B[n]_q]$ the commutative monoid freely generated by non-degenerate $q$-simplices of $\Delta[n]$. The generic element of $O[n]_q^+$ is a formal sum\[
        c = \sum_{[\mathbf a] \in B[n]_q} c_{[\mathbf a]} \cdot [\mathbf a]\ ,
        \quad c_{[\mathbf a]} \in \mathbb{N}
    \]
     where only finitely many coefficients $c_{[\mathbf a]}$
    are non-zero.
\end{itemize}
There are canonical inclusions $\Delta[n]_q\supseteq B[n]_q\subseteq O[n]_q^+\subseteq O[n]_q$.
\end{notn}

The augmented directed chain complex $O[n]$ is the algebraic model of the $n$-oriental $\cO[n]$, in a sense that will be made precise in~\cref{orientalandnu}.

\begin{ex}[{\cite[Example~3.8]{SteinerEmbedding}}]
\label{AlgebraicOrientals}
For $n\ge-1$, we define an augmented directed chain complex $(O[n],O[n]^+,\epsilon)$. Here,  the chain complex $O[n]$ is given in dimension $q$ by $O[n]_q$ and with differential maps $\partial_{q-1}\colon O_{q}[n]\to O_{q-1}[n]$ defined on a basis element $[\mathbf a]\in B[n]_{q}$ by
\[\partial_{q-1}[\mathbf a]=\partial_{q-1}[a_0,\dots,a_q]=\sum_{i=0}^{q}(-1)^i\cdot[a_0,\dots,\widehat{a_i},\dots,a_q]\in O[n]_{q-1}.\]
The chain complex $O[n]$ comes with an augmentation map $\epsilon\colon O[n]_0\to\mathbb Z$ defined on a basis element $[a]\in B[n]_0$ by
\[\epsilon[a]=1\in\mathbb Z.\]
\end{ex}

\subsection{Strong Steiner complexes}
\label{sec:steiner_complexes}

We introduce a class of particularly nice augmented directed chain complexes, that will turn out to be algebraic models for certain $\omega$-categories.

\begin{defn}[{\cite[Definition~3.1]{SteinerEmbedding}}]
Let $(C, C^+, \epsilon)$ be an augmented directed chain complex, and $B = (B_q)_{q\ge0}$ a family of subsets $B_q\subseteq C_q$ for $q\ge0$. The $\mathbb{N}$-graded set $B$ is a \emph{basis} for $C$ if for every $q \geq 0$ the set $B_q \subseteq C_q$
\begin{itemize}[leftmargin=*]
    \item  is a basis
    of the abelian group $C_q$, namely $C_q\cong\mathbb Z[B_q]$, and
    \item is a basis of the commutative monoid $C^+_q$, namely $C_q^+\cong\mathbb N[B_q]$.
\end{itemize} 
\end{defn}

If $C$ has a basis $B$, there are inclusions $B_q\subseteq C^+_q\subseteq C_q$.

\begin{defn}[{\cite[\textsection2]{SteinerEmbedding},\ \cite[\textsection2.7]{AraMaltsiniotisJoin}}]
\label{SupportComplex}
    Let $(C, C^+, \epsilon)$ be an augmented directed chain complex
    with basis $B$. For $q \geq 0$, an element $c$ in $C_q$
    has a canonical decomposition
    \[
        c = \sum_{b \in B_q} c_b \cdot b,
        \quad c_b \in \mathbb{Z}.
    \]
    The \emph{positive support}, resp.~\emph{negative support}, of $c$ is the finite subset of $B_q$ given by
    \[\supp^+c=\{c_b:c_b>0\},\quad\text{resp.~}\quad\supp^-c=\{c_b:c_b<0\}.\]
    The \emph{support} $\supp(c)$ is the union of the positive and negative support, that is the set of generators appearing in the linear expansion of $c$ under the basis $B_q$.
    The \emph{positive part}, resp.~\emph{negative part}, of $c$ is
    \[
        c^+ = \sum_{\substack{b \in B_q \\ b \in\supp^+c}} c_b \cdot b
       \in C_q^+, \quad\text{resp.~}\quad
        c^- = -\sum_{\substack{b \in B_q \\ b \in\supp^-c}} c_b \cdot b\in C_q^+.
    \]
    We denote
    \[\partial^+c:=(\partial c)^+\quad\text{and}\quad\partial^-c:=(\partial c)^-.\]
\end{defn}

The following symbol keeps track of the iterated positive and negative parts of the boundaries of a chain.

\begin{notn}[{\cite[\textsection2.8]{AraMaltsiniotisJoin}}]
\label{def:atom}
	Let $(C, C^+, \epsilon)$ be an augmented directed chain complex.
	For $\alpha=+,-$ and $c\in C_q$, the symbol $\atom{c}^\alpha_{p}$ is defined inductively on $p=q,\dots,0$ by
	\[\atom{c}^\alpha_p=\left\{\begin{array}{cr}
	    c & p=q\  \\
	    \partial^\alpha(\atom{c}^\alpha_{p+1}) & 0\le p<q.
	\end{array}\right.\]
	\end{notn}

	\begin{defn}[{\cite[Definition 3.4]{SteinerEmbedding},~\cite[\textsection2.8]{AraMaltsiniotisJoin}}]
	\label{defn:unital}
	Let $(C, C^+, \epsilon)$ be an augmented directed chain complex with basis $B = (B_q)_{q\ge 0}$
	The basis $B$ of $C$ is said to be \emph{unital}
	if \[\epsilon(\atom{c}^-_0) = 1 = \epsilon(\atom{c}^+_0)\]
	for all $c$ in $B_q$ and $q\ge0$.
\end{defn}

\begin{lem}[{\cite[\textsection3]{SteinerEmbedding}}]
    If an augmented directed chain complex $C$ admits a basis $B$,
    then the basis is uniquely determined.
\end{lem}
    
    In particular, it makes sense to say that $C$ admits a (unital) basis, without further specification.

\begin{notn}[{\cite[Definition 3.6]{SteinerEmbedding}}]
\label{loopfree_preorder_complexes}
 Let $C$ be an augmented directed chain complex with basis $B$. We denote by $\leq_{\mathbb N}$ the preorder relation on $\coprod_{q\ge0} B_q$ generated by the condition
    \[\text{$a\leq_{\mathbb N}b$ if $a\in B_p$, $b\in B_q$, and $a\in \supp(\partial^-_{q-1}(b))$ or $b\in \supp (\partial^+_{p-1}(a))$}.\]
\end{notn}

\begin{defn}[{\cite[\textsection3.6]{SteinerEmbedding}}]
\label{defbasisalgebraic}
    Let $C$ be an augmented directed chain complex with a unital basis $B$. If the preorder relation $\leq_{\mathbb{N}}$ is a partial order, we say that $C$ admits a \emph{strongly loop-free unital basis}.
\end{defn}

\begin{defn}[{\cite[\textsection2.15]{AraMaltsiniotisJoin}}]
\label{defbasisalgebraic}
    An augmented directed chain complex $C$ with a strongly loop-free unital basis $B$ is a \emph{strong Steiner complex}.
\end{defn}

\begin{ex}[{\cite[Theorems 3.1-3.2]{SteinerOrientals}}]
For $n\ge0$, the augmented directed chain complex $O[n]$ is a strong Steiner complex.
\end{ex}

The next subsection is devoted to explaining the connection between augmented directed chain complexes and $\omega$-categories.

\subsection{Steiner's functors}
\label{sec:nu}

Steiner constructs a pair of adjoint functors \[\lambda\colon \omega\cat\rightleftarrows\adch\colon\nu,\]
so that for any $\omega$-category $\cC$ and any augmented directed chain chain complex $C'$
there is a natural bijection
\begin{equation}
\label{eq:adjunction}
\adch(\lambda \cC,C')\cong\omega\cat(\cC,\nu C').
\end{equation}
To give a bit of context, we recall the basic data of these constructions, and refer the reader to the original sources for more details.

For $C$ an augmented directed  chain complex, 
the set of $q$-cells  of the
$\omega$-category $\cC:=\nu C$ for $q \geq 0$ is given by
\[(\nu C)_q=\{x:x\text{ Steiner table in }C\},\]
where a \emph{Steiner table} in $C$ is a matrix
\[
	x=\tabld{x}{q}
\]
such that, for $\alpha=+,-$ 
and $0\le p\le q$, the following hold:
\begin{enumerate}
	\item $x^\alpha_p$ belongs to $C^+_p$;
	\item $\partial(x^\alpha_p) = x^+_{p-1} - x^-_{p-1}$ for $0<p\leq q$;
	\item $\epsilon(x^\alpha_0) = 1$;
	\item $x_q^- = x_q^+$.
\end{enumerate}
We refer to~\cite[Definition~2.8]{SteinerEmbedding} or to~\cite[\textsection 2.4]{AraMaltsiniotisJoin}
for a full description of the $\omega$-categorical structure of $\nu C$.

For $\cC$ an $\omega$-category, the abelian group $(\lambda\cC)_q$ of $q$-chains of $\lambda\cC$ is the quotient of $\mathbb Z[\cC_q]$ given by
\begin{equation}
\label{lambda}
(\lambda\cC)_q=\frac{\mathbb Z[\cC_q]}{\left<[x\ast_p y]_q-[x]_q-[y]_q:x,y\in\cC_q;p<q\right>}.\ 
\end{equation}
The positivity submonoid $(\lambda\cC)_q^+$ is the submonoid of $(\lambda\cC)_q$ generated by the collection of elements $[f]_q$ for $f$ a $q$-cell of $\cC$.
The differential map $\partial_{q-1}\colon(\lambda\cC)_q\to(\lambda\cC)_{q-1}$ is determined by the condition on generators $f\in\cC_q$ given by
\[\partial_{q-1}([f]_{q})=[t_q(f)]_{q-1}-[s_q(f)]_{q-1},\]
and the augmentation map $\epsilon\colon(\lambda\cC)_0\to\mathbb Z$ by the condition on generators $x\in\cC_0$ given by
\[\epsilon([x]_0)=1.\]

Let's illustrate a couple of general phenomena about the construction $\lambda\cD$ by exploring an explicit example.

\begin{ex}
Consider the $\omega$-category $\cO[3]$.
\begin{center}
    \centering
    \begin{tikzpicture}
    \draw (-3.0,-1) rectangle (3.0,1);
			\square{
				/square/label/.cd,
	     			0 = {$x$}, 1={$y$}, 2={$z$}, 3={$w$},
                    01={$f$}, 12={$g$}, 23={$h$},
                    02={$i$}, 13={$j$}, 03={$k$},
                    023={$\alpha$}, 012={$\beta$},
                    013={$\gamma$}, 123={$\delta$},
                    0123={$\Gamma$}
     			}
	\end{tikzpicture}
\end{center}
\begin{enumerate}[leftmargin=*]
    \item Although $[f]_1\neq0\in(\lambda\cO[3])_1$, we have
    $[\id_2(f)]_2=0\in(\lambda\cO[3])_2$. This illustrates a general principle: if a $q$-cell $x$ in $\cC$ is trivial the the class $[x]_q$ vanishes in $(\lambda\cC)_q$.
    \item We have $[\beta]_2=[\id_2(h) \ast_0 \beta]_2$. More generally, if two $q$-cells of $\cC$ differ by a $(q-1)$-dimensional whisker, they represent the same class in $(\lambda\cC)_q$.
    \item We have $[(\id_2(h)\ast_0 \beta) \ast_1 \alpha]_2=[\beta]_2+[\alpha]_2$. 
    Roughly speaking, in $\lambda\cC$ composition becomes addition.
\end{enumerate}
\end{ex}

The core result of Steiner's theory that we need is the following.

\begin{thm}[{\cite[Theorem~5.11]{SteinerEmbedding}}]
\label{nuFullyFaithful}
The functor $\nu$ is an equivalence 
of categories when restricted to the full subcategory of strong Steiner complexes and corestricted to its essential image. In particular, it induces a natural bijection
\[\adch(C,C')\cong \omega\cat(\nu C,\nu C')\]
for any strong Steiner complexes $C$ and $C'$.
\end{thm}


\subsection{Strong Steiner $\omega$-categories}

\cref{nuFullyFaithful} motivates the study of the following class of $\omega$-categories, which are modelled by strong Steiner complexes.

 \begin{defn}[{\cite[\textsection2.15]{AraMaltsiniotisJoin}}]
\label{defnstrongSteiner}
An $\omega$-category $\cC$ is a \emph{strong Steiner $\omega$-category} if there exists a strong Steiner complex $C$ such that
\[\cC\cong\nu C.\]
\end{defn}

This means
that the adjunction $\lambda\colon \omega\cat\rightleftarrows\adch\colon\nu$
restricts to an equivalence of categories
\begin{equation}
\label{SteinerEquivalence}  \lambda\colon\{\cC\text{ strong Steiner $\omega$-category}\}\simeq\{C\text{ strong Steiner complex}\}\colon\nu
\end{equation}
between the full subcategory of strong Steiner complexes and the full subcategory 
of strong Steiner $\omega$-categories.

So we may regard $C$ as an algebraic model of the categorical object $\cC$. As an important instance, according to this identification, $O[m]$ corresponds to $\cO[m]$.

\begin{ex}[{\cite[Theorem~3.2]{SteinerOrientals}}]
\label{orientalandnu}
For $n\ge0$, there are isomorphisms of $n$-categories\footnote{One could also take~\cref{orientalandnu} as a definition of the $n$-th oriental $\cO[n]$.} and of augmented directed chain complexes:
\[\nu O[n]\cong\cO[n]\quad\text{ and }\quad O[n]\cong\lambda\cO[n].\]
\end{ex}

In particular, we can use Steiner's correspondence (\ref{SteinerEquivalence}) to describe nerves of strong Steiner $\omega$-categories, a task that is in general highly nontrivial, given the complicated geometry of $\cO[n]$, as opposed to $O[n]$.

\begin{rmk}
\label{descriptionnerve}
For any strong Steiner $\omega$-category -- or more generally any $\omega$-category of the form $\cC=\nu C$ -- there are natural bijections
\[\begin{array}{rcllll}
    N_m\cC & = & \omega\cat(\cO[m],\cC)&\text{\cref{RSnerve}}\\
    &\cong&\omega\cat(\cO[m],\nu C)&\\
    &\cong&\adch(\lambda\cO[m], C)&\text{\eqref{eq:adjunction}}\\
    &\cong & \adch(O[m],C).&\text{\cref{orientalandnu}}
\end{array}\]
One could verify (cf.~\cite[\textsection 3]{SteinerComplicialNerves}) that, under this identification, a simplex $x\colon\cO[m]\to\cC$, regarded as an augmented directed chain map $x\colon O[m]\to C$,
is marked if and only if \[x[0,1,\dots,m]=0\in C_m.\]
\end{rmk}

The class of strong Steiner $\omega$-categories contains many other examples of interest.

\begin{rmk}
Examples of strong Steiner $\omega$-categories include:
\begin{itemize}[leftmargin=*]
    \item the $n$-disk $\cD[n]$ and the $n$-th sphere $\cS[n]$ for $n\ge0$;
    \item the $n$-th oriental $\cO[n]$ for $n\ge0$, as shown in \cite[Example~3.8]{SteinerEmbedding} and recalled in \cref{ex:oriental};
    \item the object
    $[m|k_1,\dots,k_m]$ of Joyal's \cite{JoyalDisks} cell category $\Theta_n$ for $n\ge0$, $m\ge0$, $i=1,\dots,m$ and $k_i\ge0$, as shown in \cite{SteinerSimpleOmega};
    \item the join $\cD\star\overline\cD$ of two strong Steiner $\omega$-categories $\cD$ and $\overline\cD$, as a direct consequence of how the join is constructed in \cite[Th\'eor\`eme~6.29]{AraMaltsiniotisJoin} and recalled in \cref{joincompatible};
    \item the lax Gray
    tensor product $\cD\otimes\overline\cD$ of two strong Steiner $\omega$-categories $\cD$ and $\overline\cD$, as a direct consequence of how the join is constructed in \cite[Th\'eor\`eme~A.15]{AraMaltsiniotisJoin}.
\end{itemize}
\end{rmk}

However, there are various interesting phenomena which cannot be captured by strong Steiner $\omega$-categories.

\begin{rmk}
Examples of $\omega$-categories that are not strong Steiner $\omega$-categories include:
\begin{itemize}[leftmargin=*]
    \item any $\omega$-category $\cC$ that contains non-trivial endomorphisms;
    \item any $\omega$-category $\cC$ that is not freely generated by a polygraph in the sense of \cite{BurroniWordProblems} (or a computad in the sense of \cite{Street2Computads,BataninComputadsMonads}). 
    \item the odd dual $\cC^{\op}$ (in the sense of \cite[\textsection1.8]{AraMaltsiniotisJoin}) of the $\omega$-category $\cC$ from \cite[Exemple~2.16]{AraMaltsiniotisJoin}, as was pointed out to us by Maltsiniotis; in particular, the op-dual of a strong Steiner $\omega$-category is generally not a strong Steiner $\omega$-category.
\end{itemize}
\end{rmk}

\begin{rmk}
In virtue of the equivalence of categories \eqref{SteinerEquivalence}, if $\cC$ is a strong Steiner $\omega$-category then necessarily $\cC\cong\nu C$ where $C=\lambda\cC$. In particular, determining whether a given $\omega$-category is a strong Steiner category involves a priori understanding\footnote{More precisely, Steiner \cite{SteinerEmbedding} shows that $\cC$ is a strong Steiner $\omega$-category if and only if $\cC$ admits a \emph{strongly loop-free atomic basis}, where all the involved notions are defined in terms of $\lambda\cC$.}
the associated augmented directed chain complex $\lambda\cC$. A detailed analysis of what it means for an $\omega$-category $\cC$ to be strong Steiner purely in terms of its $\omega$-categorical structure -- and without explicit mention of $\lambda\cC$ -- is the subject of our paper \cite{GOR2}.
\end{rmk}


\subsection{Join of chain complexes and of $\omega$-categories}

 Given chain complexes $(D,\epsilon)$ and $(\overline D,\overline\epsilon)$, we recall the construction, from \cite[\textsection 6.5]{AraMaltsiniotisJoin}, of the \emph{join} augmented chain complex $(D\star\overline D,\widetilde\epsilon)$.
 
The abelian group $(D\star \overline D)_q$ of $q$-chains
is given by
\begin{equation}
 \label{FormulaJoin}   
(D\star \overline D)_q=\bigoplus_{\substack{k+ 1 + \overline k = q \\ k \geq -1, \overline k \geq -1}} D_{k}\otimes \overline D_{\overline k},\end{equation}
following the convention that $D_{-1}=\overline D_{-1}=\mathbb{Z}[\emptyset]$,
which is the free group generated by one element denoted $\emptyset$.

The differential map $\widetilde\partial_q\colon (D\star\overline D)_{q+1}\to(D\star\overline D)_q$ is defined on elements of the form $d\star\overline d\in D_k\star\overline D_{\overline k}\subseteq (D\star\overline D)_q$ by
\[
 \widetilde\partial(d \star \overline d)=\left\{
\begin{array}{lll}
    \partial d \star \overline d + (-1)^{k+1} \cdot d \star \overline\partial\overline d &  k,\overline k>0\\
    \epsilon d \star \overline d - d \star \partial\overline d &k=0,\overline k>0 \\
    \partial d \star \overline d + (-1)^{k+1} \cdot d \star \overline\epsilon\overline d &k>0,\overline k=0 \\
     \emptyset \star \overline\partial\overline d&k=-1,\overline k=0\\
     \partial d \star \emptyset& k=0,\overline k=-1.
\end{array}
\right.
\]
Here, we use the fact that the join, being an instance of a tensor product, is linear in each variable, namely we have 
$(a+b)\star c=a\star c+b\star c$.

The augmentation map $\widetilde\epsilon\colon (D\star\overline D)_0\to\mathbb Z$ is given on elements of the form $d\star\emptyset\in D_0\star\mathbb Z[\emptyset]$ and $\emptyset\star\overline d\in\mathbb Z[\emptyset]\star\overline D_0$ by
\[
    \widetilde\epsilon(d \star \emptyset) = \epsilon d
    \quad \text{and} \quad
    \widetilde\epsilon(\emptyset \star \overline d) = \overline\epsilon\overline d.
\]

The join construction extends to the context of augmented directed chain complexes.

\begin{defn}[{\cite[\textsection 6.7]{AraMaltsiniotisJoin}}]
Given augmented directed chain complexes $D=(D,D^+,\epsilon)$ and $\overline D=(\overline D,\overline D^+,\overline\epsilon)$, the \emph{join} augmented directed chain complex $((D\star \overline D),(D\star \overline D)^+,\widetilde\epsilon)$ is the augmented chain complex $(D\star \overline D,\widetilde\epsilon)$ endowed with the positivity submonoid given by
\[
    (D \star \overline D)^+_q = \{d \star \overline d : d \in D^+_{k}, \overline d\in \overline D^+_{\overline k},\ k, \overline k \geq -1,\ k+1+\overline k = q\},
\]
following the convention that $D^+_{-1} = \overline D^+_{-1} = \mathbb{N}[\emptyset]$, which is the free monoid generated by one element denoted $\emptyset$.
\end{defn}
The join construction defines a functor $\star\colon\adch\times\adch\to\adch$, which gives a monoidal structure on $\adch$ whose monoidal unit is $O[-1]$, the trivial chain complex with trivial augmentation from \cref{AlgebraicOrientals}.

For the main result of this paper, we will be interested in the case of one of the two factors of the join being a copy of $O[0]$. More precisely, we will consider the situation where $D=O[\top]$, resp.~$\overline D=O[\bot]$, is the (unique) strong Steiner complex with no non-zero chains of positive dimension and generated in dimension $0$ by a unique positive chain denoted $\top$, resp.~$\bot$.

When $\overline D=O[\top]$, resp.~$D=O[\bot]$, the formula from (\ref{FormulaJoin}) specializes to
\[(D\star O[\top])_q \cong \bigoplus_{k=0}^{q}(D_{k}\otimes\mathbb Z[\emptyset]) \oplus ( D_{k-1}\otimes \mathbb Z[\top]),\]
\[\text{resp.}\ (O[\bot]\star \overline D)_q \cong \bigoplus_{\overline k=0}^{q}(\mathbb Z[\bot]\otimes \overline D_{\overline k-1}) \oplus (\mathbb Z[\emptyset]\otimes\overline D_{\overline k}).\]

\begin{ex}[{\cite[Lem.~7.6]{AraMaltsiniotisJoin}}]
\label{orientaljoin}
For any $n\ge0$, there is an isomorphism
\[O[n]\cong\underbrace{O[0]\star\dots\star O[0]}_{n+1}.\]
It follows that, for any $k\ge-1$ and $\overline k\ge-1$ for which $k+1+\overline k=m$, there is a canonical isomorphism of augmented directed chain complexes
\[O[k]\star O[\overline k]\cong O[m].\]
\end{ex}

\begin{thm}[{\cite[Corollary~6.22]{AraMaltsiniotisJoin}}]
\label{thmjoin}
The join $D\star\overline D$ of two strong Steiner complexes $D$ and $\overline D$ is a strong Steiner complex.
\end{thm}

\label{rmk:join_basis}
 If $D$ and $\overline D$ are two strong
Steiner complexes with bases $B$ and $\overline B$, then the basis of $D \star \overline D$
is given by the set
\[
 B \star \overline B = \{d \star \emptyset : d \in B\} \cup \{d \star \overline d : d \in B,\ \overline d \in \overline B\}
                    \cup \{\emptyset \star \overline d : \overline d \in \overline B\},
\]
or, under the convention $B_{-1} = \overline B_{-1} = \emptyset$, by the set
\[B \star \overline B=\{d \star \overline d : d \in B,\ \overline d \in \overline B\}.\]

The following fact is a crucial property of the join of strong Steiner complexes. Recall that there is a copy of $D$ that lives inside the join $D\star\overline D$, namely $D\cong D\star O[-1]\subseteq D\star\overline D$, and similarly for $\overline D$. Roughly speaking, the proposition formalizes the fact that in a join of strong Steiner complexes $D\star\overline D$ if a chain map $x$ ``ends'' in the copy of $D$ then it must fully take values there.

\begin{prop}
\label{prop:last_vertex_limit}
    Let $D$ and $\overline D$ be two strong Steiner complexes, $m\geq 0$ and $x \colon O[m] \to D \star \overline D$ an augmented directed chain map.
    \begin{itemize}[leftmargin=*]
        \item If $x[m]$ -- namely the value of $x$ on the last $0$-chain of $O[m]$ -- belongs to $D\star O[-1]\subset D\star\overline D$, then so does $x[\mathbf a]$ for every $[\mathbf a]\in B[m]_q$.
        \item If $x[0]$ -- namely the value of $x$ on the first $0$-chain of $O[m]$ -- belongs to $O[-1]\star \overline D\subseteq D\star\overline D$, then so does  $x[\mathbf a]$ for every $[\mathbf a]\in B[m]_q$.
    \end{itemize}
\end{prop}

To prove the proposition, we need a preliminary lemma.

\begin{lem}
\label{NoZeroBoundary}
Let $C$ be a strong Steiner complex, $q\ge0$ and $c\in C^+_{q+1}$ a positive $(q+1)$-chain in $C$. Then we have
\[\partial_q^+c\neq0\in C_q\quad\text{ and }\quad\partial^-_qc\neq0\in C_q.\]
\end{lem}

\begin{proof}[Proof of \cref{NoZeroBoundary}]
We prove that $\partial^+c\neq0$ arguing by contradiction. The argument for $\partial^-c\neq0$ is analogous.

If $c\in B_{q+1}$ is a basis element of $C$, we see looking at~\cref{def:atom} that $\partial^+c=0$ would imply that $\atom{c}^+_0=0$, contradicting the unitality of the basis $B$ of $C$. So the statement holds for $c$ being a basis element.

Now we turn to the general case. Let $b_0$ be a maximal element of the finite set $
\supp^+c$ with respect to the partial order $\leq_{\mathbb N}$.
By the first part of the proof, there exists at least one $b_0'\in \supp^+(\partial b_0)$, so in particular $b_0\leq_{\mathbb{N}} b_0'$, while $b_0\neq b_0'$ for dimension reasons.

If the expression of $c$ is
\[c=c^+=\sum_{b\in \supp^+c} c_b \cdot b\ ,\quad c_b\geq 0,\]
we obtain
\[
\partial c=\sum_{b\in \supp^+c} c_b \cdot\partial b=\sum_{b\in \supp^+c} c_b \cdot(\partial^+ b-\partial^- b)=\sum_{b\in\supp^+c} c_b\cdot  \partial^+b - \sum_{b\in \supp^+c} c_b\cdot  \partial^-b.
\]
Since we assumed $\partial^+c=0$, all the summands of the first sum need to be cancelled by summands in the second sum.
In particular, there has to be a $b_1 \in \supp^+c$ with $b_0'\in\supp( \partial^-b_1)$ in order to cancel the summand corresponding to $b_0'$ of the first sum. This implies in particular $b_0'\leq_{\mathbb{N}} b_1$ and
 so by transitivity $b_0 \leq_{\mathbb{N}}b_0'
 \leq_{\mathbb{N}} b_1$ while $b_1\neq b_0$, contradicting the maximality of $b_0$. This shows that $\partial^+c\neq0$, as desired.
\end{proof}

With the notation of the proof above,
it is not true in general that \[\supp(\partial^+ b_0)\subseteq\supp(\partial^+ c),\] since some cancelling might occur if $b_0$ is not maximal.

\begin{rmk}
Consider the  augmented directed chain complex $C=\lambda\cC$, where $\cC$ is the $1$-category
    \[\cC=\fbox{$x \xrightarrow{f} y \xrightarrow{g} z$}.\]
Let $c=[g]_1 + [f]_1=[g\circ f]_1$ and $b_0=[f]_1$, which is a non-maximal element of $\supp^+(\partial c)$.
In this situation, we have that
\[\supp^+(\partial[f]_1)=\{[y]_0\}\quad\supp^+(\partial([f]_1+[g]_1))=\{[z]_0\}\]
and in particular
\[[y]_0\in\supp^+(\partial[f]_1)\setminus\supp^+(\partial([f]_1+[g]_1)).\]
\end{rmk}

\begin{proof}[Proof of~\cref{prop:last_vertex_limit}]
The statement is trivial for $m=0$ and for $m=1$ it follows immediately
by the fact that there is no basis element $b$
of $D \star \overline{D}$ of dimension $1$ such that
$\partial^+_0 b $ is in $D_0$ and $\partial^-_0 f$ is in $\overline{D}_0$
by definition.

Assume $q>1$ is the smallest dimension in which we have a basis element $[\mathbf{a}]$ of the left-hand side mapped to a chain in $ D\star \overline D\setminus D\star \emptyset$, say,
\[
x[\mathbf{a}]=\sum_{i\in I} \lambda_i\cdot  c_i\star \overline c_i, 
\]
with at least one $\lambda_i> 0$ for some corresponding $\overline c_i$ of dimension $\geq 0$. Let $0\leq p\leq q$ be the maximal dimension of such $\overline c_i$, and let $I'\subseteq I$ be the subset corresponding to all summands for which the dimension of their $\overline c_i$ is exactly~$p$. 

Then on the one hand, the chain $\partial x[\mathbf{a}]= x\partial[\mathbf{a}]$ is in $D\star \emptyset$ by minimality of $q$.
Thus, the sum
\[
\partial \left (\sum_{i\in I} \lambda_i \cdot c_i\star \overline c_i\right)=\sum_{i\in I} \lambda_i\cdot  (\partial c_i)\star \overline c_i+\sum_{i\in I} (-1)^{\lvert c_i\rvert+1}\lambda_i\cdot  c_i\star (\partial\overline c_i),
\]
turns out to be in $D\star \emptyset$.
\begin{equation}
 \label{eq:join_no_way_back}
   \sum_{i\in I'} \lambda_i\cdot (\partial c_i)\star \overline c_i=0
\end{equation}
has to vanish. Indeed, by the maximality of $p$ such terms cannot be canceled out by any other term of the sum above. Now, fix a $\overline{c}_j$ for
some $j$ in $I'$ and denote by $I_j$ the subset of $I'$ consisting of
the indices $i$ such that $\overline{c}_i = \overline{c}_j$;
said otherwise, if $\overline{c}_\bullet \colon I' \to \overline{D}_p$
denotes the indexing function, let $I_j$ be $(\overline{c}_\bullet)^{-1}(\overline{c}_j)$.
As a consequence of the vanishing of the sum~\eqref{eq:join_no_way_back},
we must have that
\[
 \sum_{i\in I_j} \lambda_i\cdot (\partial c_i)\star \overline c_j = 0,
\]
which is equivalent to saying that \[\sum_{i\in I_j} \lambda_i\cdot (\partial c_i)=0.\]
But this last sum is precisely the boundary of the positive chain $\sum_{i\in I'} \lambda_i \cdot c_i$ of $D$ and the fact that it must vanish contradicts~\cref{NoZeroBoundary}. This completes the proof.
\end{proof}

Given Steiner's equivalence from (\ref{SteinerEquivalence}), the join construction
for strong Steiner complexes,
\[\star\colon\{\text{strong Steiner cpx}\}\times\{\text{strong Steiner cpx}\}\to\{\text{strong Steiner cpx}\},\]
can be used to define a join construction for strong Steiner $\omega$-categories,
\[\star\colon\{\text{strong Steiner $\omega$-cat}\}\times\{\text{strong Steiner $\omega$-cat}\}\to\{\text{strong Steiner $\omega$-cat}\}.\]
Ara and Maltsiniotis showed that we can extend this join construction to the whole category $\omegacat$
in an essentially unique fashion.

\begin{thm}[{\cite[Theorem~6.29]{AraMaltsiniotisJoin}}]
\label{joincompatible}
There exists a unique -- up to unique monoidal isomorphism -- monoidal structure $\star\colon\omega\cat\times\omega\cat\to\omega\cat$ on $\omegacat$, called the join of $\omega$-categories, such that
\begin{itemize}[leftmargin=*]
    \item for any strong Steiner complexes $D$ and $\overline D$
there is a natural isomorphism of $\omega$-categories
\[\nu D\star\nu \overline D\cong\nu(D\star \overline D);\]
\item the functor $-\star-$ commutes with colimits in each variable.
\end{itemize}
\end{thm}

The functor $\nu$ provides us with the following algebraic descriptions for nerves of joins of strong Steiner $\omega$-categories.

\begin{rmk}
By~\cref{joincompatible}, for any strong Steiner $\omega$-categories $\cD=\nu D$ and $\overline \cD=\nu \overline D$, there are natural bijections
\[\begin{array}{lll}
    N_m(\cD\star\overline\cD) & = & \omega\cat(\cO[m],\cD\star\overline\cD)\\
    &\cong&\omega\cat(\cO[m],\nu D\star\nu \overline D)\\
    & \cong & \omega\cat(\nu O[m],\nu(D\star \overline D))\\
    &\cong & \adch(O[m],D\star \overline D),
\end{array}\]
with a similar description of the marking as the one from~\cref{descriptionnerve}.

When $\overline\cD=\cO[\top]$, resp.~$\cD=\cO[\bot]$, is the terminal category with a single object $\top$, resp.~$\bot$, we obtain the bijections
\[N_m(\cD\star\cO[\top])\cong\adch(O[m],D\star O[\top]),\]
\[\text{resp.}\  N_m([\bot]\star\overline\cD)\cong\adch(O[m], O[\bot]\star \overline D).\]

Finally, when $\cD=\cO[n-1]$, using~\cref{orientaljoin} we obtain the bijection
\[\begin{array}{llll}
N_m\cO[n]&\cong&N_m(\cO[n-1]\star\cO[\top])\\
&\cong&\adch(O[m], O[n-1]\star O[\top])\\
&\cong&\adch(O[m], O[n]).
\end{array}\]
\end{rmk}

\section{Features of simplices of nerves of cones}

Throughout the document, let $\cD=\nu D$ be a strong Steiner $\omega$-category. We introduce important features and constructions based on simplices of $N_m(\cD\star\cO[\top])$, which can be identified using~\cref{descriptionnerve,thmjoin} with augmented directed chain maps
\begin{equation}
\label{simplex}
 x\colon O[m]\to D\star O[\top].
\end{equation}

Most constructions and results are inspired by and generalize Steiner's constructions from \cite{SteinerUniversal}, that focuses on the case of $\cD=\cO[n]$
being an oriental.

An analogous analysis applies to simplices of $N_m(\cO[\bot]\star\cD)$,
which have been identified with chain maps
\begin{equation}
\label{dualsimplex}
x\colon O[m]\to O[\bot]\star D.\end{equation}

For most constructions and results in the paper we leave it to the reader to adapt the statements for this case, but we will point out when there are subtleties in adapting constructions, statements or proofs.

\begin{digression}
\label{dualnotformal}
If $D$ is such that its \emph{odd dual} $D^{\op}$, in the sense of \cite[\textsection2.18]{AraMaltsiniotisJoin}, is a strong Steiner complex, then all constructions and results for a simplex of the form (\ref{dualsimplex}) can be obtained by applying the odd dual $(-)^{\op}$ to the corresponding statements for (\ref{simplex}), and using the compatibility of the op-dual construction with the join construction with swapped factors.  However, in general the odd dual of a strong Steiner $\omega$-category
is \emph{not} a strong Steiner $\omega$-category.\footnote{An example of this phenomenon, which was pointed out by Maltsiniotis, is the odd dual $\cC^{\op}$ of the $\omega$-category $\cC$ from \cite[Exemple~2.16]{AraMaltsiniotisJoin}.} Thus, all properties of simplices of the form (\ref{dualsimplex}) and the proof of~\cref{maintheoremdual} need to be treated on their own, and cannot be deduced formally from the properties of simplices of the form (\ref{simplex}) and from the proof of~\cref{maintheorem}.
\end{digression}

\subsection{Faces and degeneracies of a simplex}
Let $J$ be an augmented directed chain complex, and we consider chain maps of the form $x\colon O[m]\to J$, having in mind the relevant situation where $J$ is of the form $J=D\star O[\top]$ or $J=O[\bot]\star D$, for $D$ any augmented directed chain complex.

As we let $m$ vary, one can define the following structure maps.
Let $s^j\colon[m]\to[m-1]$ and $d^j\colon[m]\to[m+1]$ with usual ranges denote coface and codegeneracy maps in $\Delta$.

\begin{defn}
\label{DefFaceDeg}
Let $x\colon O[m]\to J$ be a chain map.
\begin{enumerate}[leftmargin=*]
    \item For $m\ge1$ and $0\le i\le m$, the \emph{$i$-th face} of $x$ is the chain map
    \[d_ix\colon O[m-1]\to J\] 
    given on basis elements by
    \[(d_ix)[\mathbf a]=x[d^i(\mathbf a)].\]
    Alternatively,
    \[( d_ix)[\mathbf a_{<i},\mathbf a_{\ge i}]=x[\mathbf a_{<i},\mathbf a_{\ge i}+1],\]
    where $\mathbf{a}_{<i}\subseteq\{0,\dots,i-1\}$ and  $\mathbf{a}_{\ge i}\subseteq\{i,\dots,m-1\}$.
\item For $m\ge0$ and $0\le i\le m$, the \emph{$i$-th degeneracy} of $x$ is the chain map
\[s_ix\colon O[m+1]\to J\] given on basis elements by
\[(s_ix)[\mathbf a]=x[s^i(\mathbf a)].\]
Alternatively,
\[
\left\{
\begin{array}{lcl}
(s_ix)[\mathbf{a}_{<i},\mathbf{a}_{>{i+1}}]=x[\mathbf{a}_{<i},\mathbf{a}_{>{i+1}} -1]
\\
(s_ix)[\mathbf{a}_{<i},i, \mathbf{a}_{>{i+1}}]=(s_ix)[\mathbf{a}_{<i},i+1,\mathbf{a}_{>{i+1}}]=x[\mathbf{a}_{<i},i,\mathbf{a}_{>{i+1}}-1]
\\
(s_ix)[\mathbf{a}_{<i},i,i+1,\mathbf{a}_{>{i+1}}]=0,
\end{array}
\right.
\]
where $\mathbf{a}_{<i}\subseteq\{0,\dots,i-1\}$ and  $\mathbf{a}_{>{i+1}}\subseteq\{i+2,\dots,m+1\}$. 
\end{enumerate}
\end{defn}

\begin{rmk}
 \label{SimpDegLinear}
Let $x,y\colon O[m]\to J$ be chain maps.
\begin{enumerate}[leftmargin=*]
 \item For $m\ge1$ and $0\le i\le m$, we have
    $d_i(x+y)= d_ix+ d_iy$.
        \item For $m\ge0$ and $0\le i\le m$, we have
     $s_i(x+y)= s_ix+s_iy$.
\end{enumerate}

\end{rmk}

\begin{rmk}
\label{NerveOsSet}
Let $x\colon O[m]\to J$ be an augmented directed  chain map.
\begin{enumerate}[leftmargin=*]
 \item For $m\ge1$ and $0\le i\le m$, the map
    $d_ix$ is an augmented directed chain map.
        \item For $m\ge0$ and $0\le i\le m$,
        the map
    $s_ix$ is an augmented directed chain map.
\end{enumerate}
\end{rmk}

Since the assignment  $m\mapsto O[m]$ defines a cosimplicial augmented directed chain complex $O[\bullet]$, we know that $\adch(O[\bullet],J)$ is a simplicial set, whose face and degeneracy maps are precisely those described explicitly in~\cref{DefFaceDeg}. In particular, these face and degeneracy maps satisfy the usual simplicial identities, which we recall for the reader's convenience.

\begin{prop}
\label{SimpIdentities}
Let $x\colon O[m]\to J$ be a chain map, and $0\le i,j\le m$. The following identities hold.
\begin{itemize}
    \item
$s_is_jx=s_{j+1}s_ix\quad\text{ if }\quad0\leq i\leq j\leq m$
\item 
$d_is_jx=\left\{
\begin{array}{lll}
s_{j-1}d_ix&\text{ if }\quad 0\leq i<j \leq m\\
x&\text{ if }\quad i \in\{j,j+1\} \\
s_jd_{i-1}&\text{ if }\quad 1\leq j+1 < i \leq m
\end{array}\right.$
\item
$d_id_jx=d_{j-1}d_ix\quad\text{ if }\quad0\leq i<j\leq m.$
\end{itemize}
\end{prop}

We will also make use of the following ``iterated'' simplicial identities, which can be easily deduced.

\begin{rmk}
\label{IteratedSimpIdentities}
Let $x\colon O[m]\to J$ be a chain map.
The following identities hold.
\begin{itemize}
    \item $s_j^{\ell}x=s_{j+\ell-1}s_{j+\ell-2}\ldots s_{j+1}s_jx$, for $0 \leq j \leq m$, $\ell \geq 1$,
    \item $d_j^{\ell}x=d_jd_{j+1}\ldots d_{j+\ell-1}x$, for $0 \leq j < m$, $1 \leq \ell \leq m-j+1$,
    \item $s_{j+1}^{\ell}s_jx=s_j^{\ell+1}x$ for $0\leq j\leq m$, $\ell \geq 0$. 
    \item $d_js_k^{\ell}x=s_k^{\ell-1}x$ for $0\leq k\leq m$, $k \leq j\leq \ell + k$. 
    \item $d_k^{\ell}d_jx=d_k^{\ell+1}x$ for $0\leq k+\ell\leq m$, $k\leq j\leq \ell + k$. 
    \item $s_{j}s_k^{\ell}x=s_k^{\ell+1}x$ for $0\leq k \leq m$, $k<j\leq m+\ell$, $\ell \geq 0$
\end{itemize}
\end{rmk}

\subsection{The last factor of a simplex}

Throughout this subsection, let $D$ be a strong Steiner complex. We introduce in detail constructions $\gamma x$ and $\beta x$ based on simplices $x$ of the form (\ref{simplex}), which mostly generalize Steiner's from \cite[\textsection4]{SteinerEmbedding}, and indicate briefly in~\cref{dualgamma} how the analogous constructions work for the case of simplices $x$ of the form (\ref{dualsimplex}).

Recall there are canonical inclusions \[D\cong D\star O[-1]\hookrightarrow D\star O[\top]\hookleftarrow O[-1]\star O[\top]\cong O[\top]\]
of $D$ and $O[\top]$ into the join $D\star O[\top]$.
The following definition generalizes \cite[Definition~4.2]{SteinerEmbedding} in the relevant cases.\footnote{More precisely, our definition of rank agrees with Steiner's when $x$ does not factor through $D$.}

\begin{defn}
Let $x\colon O[m]\to D\star O[\top]$ be an augmented directed chain map.
\begin{itemize}[leftmargin=*]
    \item The \emph{rank} of $x$ is the number $r$ of elementary $0$-chains $[a]$ in $O[m]$ such that $x[a]\in D\subseteq D\star O[\top]$, namely
    \[r:=\#\{ 0\le a\le m:x[a]\in D\subseteq D\star O[\top]\}.\]
\item The \emph{corank} $s$ is one less than the number of elementary $0$-chains $[a]$ in $O[m]$ such that $x[a]\in O[\top]\subseteq D\star O[\top]$, namely
\[s:=\#\{ 0\le a\le m:x[a]\in O[\top]\subseteq D\star O[\top]\}-1.\]
\end{itemize}
\end{defn}

By definition, one always has $r+s=m$, with $0\le r\le m+1$ and $-1\le s\le m$.

\begin{rmk}
\label{remarkrank}
Given an augmented directed chain map $x\colon O[m]\to D\star O[\top]$, the rank of $x$ can also be understood as
\[r=\left\{\begin{array}{cc}
     \min\{0\le a\le m\colon x[a]=\top\}&\text{if }\{0\le a\le m\colon x[a]=\top\}\neq\emptyset  \\
    m+1 &  \text{otherwise.}
\end{array}\right.\]
\end{rmk}
\begin{rmk}
\label{totallydegenerate}
The map $x$ factors through $O[\top]$ if and only if
it is \emph{totally degenerate}, by which we mean $x$ is an iterated degeneracy of the $0$-simplex.
\end{rmk}
\begin{rmk}
The map $x$ factors through $O[\top]$ if and only if
it has rank $r=0$ and $s=m$.
Instead, $x$ factors through $D$ if and only if it has rank $r=m+1$ and corank $s=-1$.
\end{rmk}

\begin{rmk}
\label{rmk:face_rank}
    If $x\colon O[m]\to D\star O[\top]$ is an augmented directed chain map of rank $r$
    and corank $s\geq 0$,
    then for any $r < i \leq m$ we have that $d_i x \colon O[m-1]\to D\star O[\top]$
    has still rank $r$ but corank $s-1$.
\end{rmk}

\begin{defn}\label{DefGammaBeta}
Let $x\colon O[m]\to D\star O[\top]$ be an augmented directed chain map with rank $r\le m$
and corank $s\ge0$ which is given on a basis element $[\mathbf a,r]$ for $[\mathbf a]\in B[m]_q$ by
\[
    x[\mathbf a,r]=\sum_{ b \in B_{q+1}} x'_{\mathbf a, b}\cdot b\star\emptyset+ \sum_{ b \in B_q} x''_{\mathbf a, b}\cdot b\star \top.
\] 
\begin{itemize}[leftmargin=*]
    \item The \emph{last factor} of $x$ is the augmented directed chain map \[\gamma x\colon O[r]\to D\star O[\top]\] given on basis elements by
\[(\gamma x)[r]= \emptyset\star\top
\ , \quad 
(\gamma x)[\mathbf a,r]=\sum_{ b \in B_q}x''_{\mathbf a, b}\cdot b\star \top
\ , \quad(\gamma x)[\mathbf a]=\sum_{ b \in B_{q}}x''_{\mathbf a, b} \cdot b\star\emptyset,
\]
for $[\mathbf a]$
in $B[r-1]_q \subseteq B[m]_q$, and $q\ge 0$.
\item The \emph{normalized last factor} of $x$ is the augmented directed chain map
\[\beta x\colon O[m]\to D\star O[\top]\] given by
\[\beta x=s_{r}^{s}\gamma x.\]
\end{itemize}
\end{defn}

\begin{lem}
\label{gammaproperties}
Let $D$ be a strong Steiner complex. For $ x\colon O[m]\to D\star O[\top]$ with rank $r\le m$, the formula above defines an augmented directed chain map $\gamma x\colon O[r]\to D\star O[\top]$ with rank $r$ and corank $0$.
\end{lem}

In order to prove the lemma, it is handy to rewrite the formulas for $\gamma x$ in a slightly more abstract way.

Let $J:=D\star O[\top]$
denote the join and $q\ge0$.
Denote by $J(\top)_q$ the subgroup of $J_q$ generated by all basis elements of the form $b\star \top$ with $b\in B_{q-1}\subseteq D_{q-1}$,
and by $J(\neg\top)_q$ the subgroup of $J_q$ generated by all basis elements of the form $b=b\star \emptyset$ with $b\in B_{q}\subseteq D_q$.
Denote by
\[\pr_{\top}\colon J_q\to J(\top)_{q}
\text{ and }\pr_{\neg\top}\colon J_q\to J(\neg\top)_{q}\]
the canonical projections. In this language, we observe that the formulas for $\gamma x$ can be rewritten to
\[\gamma x[r]= \emptyset\star\top
\, ,\;
\gamma x[\mathbf a,r]=\pr_{\top}x[\mathbf a,r]
\, ,
\; \gamma x[\mathbf a]=(-1)^{q+1}\pr_{\neg\top}\partial\pr_{\top} x[\mathbf{a},r].
\]
for $[\mathbf a]$
in $B[r-1]_q \subseteq B[m]_q$, and $q\geq 0$.

Before proving the lemma, we need the following auxiliary fact.

\begin{lem}
\label{ProjAlmostChainMap}
Let $q\geq 0$. With the notation above, the following holds:
\[\partial_{q} \pr_{\neg\top} \partial_{q+1} \pr_{\top}=-\pr_{\neg\top} \partial_{q}  \pr_{\top}  \partial_{q+1}\colon J_{q+2}\to J(\neg\top)_{q}.\]
\end{lem}
\begin{proof}
Under the canonical isomorphism $J_{q+1}\cong J(\top)_{q+1}\oplus J(\neg\top)_{q+1}$, the differential 
$\partial_{q+1}\colon J_{q+2}\to J_{q+1}$ 
can be represented as a matrix operator
\[
\begin{pmatrix}
u'_{\top,\top} & u'_{\top,\neg\top}\\
u'_{\neg\top,\top} & u'_{\neg\top,\neg\top}\\
\end{pmatrix}
\colon J(\top)_{q+1}\oplus J(\neg\top)_{q+1}\to J(\top)_{q}\oplus J(\neg\top)_{q},\]
and similarly we write $\begin{pmatrix}
u_{\top,\top} & u_{\top,\neg\top}\\
u_{\neg\top,\top} & u_{\neg\top,\neg\top}\\
\end{pmatrix}$ for the matrix operator representing $\partial_{q}$. 
By definition, one may verify that 
\[u'_{\top,\neg\top}=0\quad \text{and}\quad u_{\top, \neg\top} = 0,\]
so the expression above simplifies to
\[
\begin{pmatrix}
u'_{\top,\top} & 0\\
u'_{\neg\top,\top} & u'_{\neg\top,\neg\top}\\
\end{pmatrix}
\colon J(\top)_{q+1}\oplus J(\neg\top)_{q+1}\to J(\top)_{q}\oplus J(\neg\top)_{q}.\]
Next, the condition $0=\partial_q  \partial_{q+1}$ reads as
\[
\begin{pmatrix}
0 & 0\\
0 & 0\\
\end{pmatrix}=\begin{pmatrix}
u_{\top,\top} & 0\\
u_{\neg\top,\top} & u_{\neg\top,\neg\top}\\
\end{pmatrix}
\begin{pmatrix}
u'_{\top,\top} & 0\\
u'_{\neg\top,\top} & u'_{\neg\top,\neg\top}\\
\end{pmatrix}
=
\]
\[=
\begin{pmatrix}
u_{\top,\top}u'_{\top,\top} & 0\\
u_{\neg\top,\top}u'_{\top,\top}+u_{\neg\top,\neg\top}u'_{\neg\top,\top} & u_{\neg\top,\neg\top}u'_{\neg\top,\neg\top}\\ 
\end{pmatrix},
\]
so in particular
\begin{equation}
  \label{auxiliaryfact}  
0=u_{\neg\top,\top}u'_{\top,\top}+u_{\neg\top,\neg\top}u'_{\neg\top,\top}.
\end{equation}

The matrix corresponding to the left hand side $\partial_{q} \pr_{\neg\top} \partial _{q+1} \pr_{\top}$ is
\[\begin{pmatrix}
u_{\top,\top} & 0\\
u_{\neg\top,\top} & u_{\neg\top,\neg\top}\\
\end{pmatrix}
\begin{pmatrix}
0 & 0\\
0 & \id\\
\end{pmatrix}
\begin{pmatrix}
u'_{\top,\top} & 0\\
u'_{\neg\top,\top} & u'_{\neg\top,\neg\top}\\
\end{pmatrix}
\begin{pmatrix}
\id & 0\\
0 & 0\\
\end{pmatrix}
=
\begin{pmatrix}
0 & 0\\
u_{\neg\top,\neg\top}u'_{\neg\top,\top} & 0\\
\end{pmatrix},\]
and the matrix corresponding to the right hand side $-\pr_{\neg\top} \partial_{q}  \pr_{\top}  \partial_{q+1}$ is
\[
-\begin{pmatrix}
0 & 0\\
0 & \id\\
\end{pmatrix}
\begin{pmatrix}
u_{\top,\top} & 0\\
u_{\neg\top,\top} & u_{\neg\top,\neg\top}\\
\end{pmatrix}
\begin{pmatrix}
\id & 0\\
0 & 0\\
\end{pmatrix}
\begin{pmatrix}
u'_{\top,\top} & 0\\
u'_{\neg\top,\top} & u'_{\neg\top,\neg\top}\\
\end{pmatrix}
=
\begin{pmatrix}
0 & 0\\
-u_{\neg\top,\top}u'_{\top,\top} & 0\\
\end{pmatrix},\]
so the desired identity follows using (\ref{auxiliaryfact}).
\end{proof}

We are now ready to prove that $\gamma x$ is a chain map with the desired properties.

\begin{proof}[Proof of~\cref{gammaproperties}]
To prove that $\gamma x$ is a map of chain complexes, it suffices to check that $\gamma x  \partial=\partial \gamma x$ on each element of the basis in degree $q$, so we distinguish three cases.
 \begin{itemize}[leftmargin=*]
     \item Consider the chain $[\mathbf{a}]$ with $a_q<r$. On the one hand, using~\cref{ProjAlmostChainMap},
     we obtain that 
\[
\begin{array}{lll}
   \partial(\gamma x) [\mathbf{a}]   &=& \partial (-1)^{q+1}\pr_{\neg\top}\partial\pr_{\top} x[\mathbf{a},r]\\
 &=    & (-1)^{q} \pr_{\neg\top}\partial \pr_{\top} \partial x[\mathbf{a},r]. 
\end{array}
\]

On the other hand, if
\[\partial[\mathbf a] =
    \sum_{\mathbf{a'} \in B_{q-1}[r-1]} \lambda_{\mathbf{a'}}\cdot[\mathbf{a'}],\]
    where $\lambda_{\mathbf{a'}} \in \mathbb{Z}$, we can write
\[
    \partial [\mathbf{a},r]=(-1)^{q+1}\cdot[\mathbf{a}]+ \partial[\mathbf a] \star [r] =
    (-1)^{q+1}\cdot[\mathbf{a}]+
    \sum_{\mathbf{a'} \in B_{q-1}[r-1]} \lambda_{\mathbf{a'}}\cdot[\mathbf{a'},r].
\] 
Thus
\[
(\gamma x) \partial[\mathbf{a}]=(-1)^{q}\sum_{\mathbf{a'} \in B_{q-1}[r-1]} \lambda_{\mathbf{a'}}\cdot\pr_{\neg\top}\partial \pr_{\top}x[\mathbf{a'},r].
\]
 Moreover, since $x$ is of rank $r$, by~\cref{prop:last_vertex_limit} we know that $\pr_{\top} x[\mathbf{a''}]=0$ 
 for any $q$-dimensional $[\mathbf{a''}]$ with $a_q''<r$. In particular, we can insert \[0=(-1)^{q+1}\pr_{\top} x[\mathbf{a}]\] in the equation above, so that  \[
 \begin{array}{cl}
(\gamma x) \partial[\mathbf{a}] & =(-1)^{q}\pr_{\neg\top}\partial\left((-1)^{q+1}\pr_{\top} x[\mathbf{a}]+\pr_{\top}\sum\limits_{\mathbf{a'} \in B[r-1]_{q-1}} \lambda_{\mathbf{a'}}\cdot x[\mathbf{a'},r]\right)\\
&= (-1)^{q}\pr_{\neg\top}\partial\pr_{\top} x\left((-1)^{q+1}\cdot[\mathbf{a}]+\sum\limits_{\mathbf{a'} \in B[r-1]_{q-1}} \lambda_{\mathbf{a'}}\cdot[\mathbf{a'},r]\right)\\
&= (-1)^{q}\pr_{\neg\top}\partial\pr_{\top} x\partial[\mathbf{a},r].
\end{array}
\]
Since $x$ was a chain map, we have $\partial x[\mathbf{a},r]= x \partial[\mathbf{a},r]$, thus yielding the desired claim.

\item Consider the chain  $[\mathbf{a},r]$ with $a_{q}<r$    and $q > 0$. On the one hand we have that
\[
\partial (\gamma x) [\mathbf{a},r]= \partial\pr_{\top}(x [\mathbf{a},r]).
\]
On the other hand, using a similar reasoning to the previous case, we have that
\[
\begin{split}
 (\gamma x) \partial[\mathbf{a},r] 
 &= (\gamma x)\Bigl((-1)^{q+1} [\mathbf{a}] +
    \sum_{\mathbf{a'} \in B_{q-1}[r-1]} \lambda_{\mathbf{a'}} \cdot [\mathbf{a'}, r]\Bigr) \\
 &=(-1)^{q+1} (-1)^{q+1}\pr_{\neg\top}\partial(\pr_{\top} x [\mathbf{a}, r]) + \\
 & \phantom{=} +\sum_{\mathbf{a'} \in B_{q-1}[r-1]} \lambda_{\mathbf{a'}}\cdot\pr_{\top}x[\mathbf{a'},r].
 \end{split}
\]

Using that $x$ is a chain map once again, we have 
\[
\partial x[\mathbf{a},r] = (-1)^{q+1}\cdot x[\mathbf{a}]+\sum_{\mathbf{a'} \in B_{q-1}[r-1]} \lambda_{\mathbf{a'}} x[\mathbf{a'},r].
\]
Inserting this into the previous formula yields
\[
 (\gamma x) \partial[\mathbf{a},r]=\pr_{\neg\top}\partial(\pr_{\top} x [\mathbf{a}, r])+\pr_{\top}\bigl(\partial x[\mathbf{a},r]- (-1)^{q+1} x[\mathbf{a}]\bigr).
\]
Using once again that $x$ is of rank $r$ and~\cref{prop:last_vertex_limit},
we have
\[
 \pr_{\top}\left( (-1)^{q+1} x[\mathbf{a}]\right)=0.
\] 
Next, we can observe 
\[
\pr_{\top}\partial x[\mathbf{a},r]=\pr_{\top}\partial \bigl(\pr_{\top} x [\mathbf{a}, r] + \pr_{\neg\top} x [\mathbf{a}, r]\bigr).
\]
Since $J(\neg\top) \cong D$ is a sub-chain complex, this simplifies to 
\[
\pr_{\top}\partial x[\mathbf{a},r]=\pr_{\top}\partial \bigl(\pr_{\top} x [\mathbf{a}, r]\bigr).
\]
Inserting this into last expression for $(\gamma x) \partial[\mathbf{a},r]$, we obtain
\[
\begin{split}
     (\gamma x) \partial[\mathbf{a},r]
        &=\pr_{\neg\top}\partial\bigl(\pr_{\top} x [\mathbf{a}, r]\bigr)
            + \pr_{\top}\partial \bigl(\pr_{\top} x [\mathbf{a}, r]\bigr) \\
        & =\partial\pr_{\top} x [\mathbf{a}, r],
\end{split}
\]
as desired.

\item Consider the chain $[a,r]$ with $0 \leq a<r$ and $q=0$.
On the one hand, we have that
\[
(\gamma x) \partial [a, r] = (\gamma x) ([r]-[a]) = \emptyset\star\top - (\gamma x) [a]. 
\]
By definition we have
\[
(\gamma x) [a] = (-1)^{1}\pr_{\neg\top}(\partial(\pr_{\top} x[a,r])). 
\]
Since $x$ is of rank $r\le m$
\[\top=\pr_{\top} \partial(\pr_{\top} x[a,r]).\]
So we obtain
\[
(\gamma x) \partial [a,r] =  \pr_{\top} \partial(\pr_{\top} x[a,r]) +\pr_{\neg\top}(\partial(\pr_{\top} x[a,r])). 
\]
On the other hand, by definition we have 
\[
\partial (\gamma x) [a,r]= \partial \pr_{\top}(x [a,r]),
\]
as desired.

 \end{itemize}

The fact that the chain map $\gamma x$ is directed is by construction, and we now prove that $\gamma x$ respects the augmentation, i.e., that any elementary $0$-chain of
$O [m]$ is sent to an element of the basis of $B_0$ or to $\top$.
By definition, this holds for $[r]$, so it is enough to check it for $[i]$, with $0 \leq i < r$.
We claim that for any $1$-chain $[i, r]$ of $O[m]$, the element
$\gamma x[i, r]$ is a basis element of $D \star O[\top]$.
Indeed, write
\[
    x[i,r]=\sum_{ b \in B_{1}} x'_{\mathbf a, b}\cdot b\star \emptyset+ \sum_{ b \in B_0} x''_{\mathbf a, b}\cdot b\star \top.
\]
We wish to show that there is exactly one summand in the second sum.
This follows from the fact that the differential of every $1$-chain of the type $b \star \top$ 
with $b \in B_0$ has $\top$ as summand, while the differential of an element $b' \in B_1$
has of course no summand of the form $\top$. Hence the differential of $x[i, r]$
has a summand of the form $q\cdot  \top$, where $q$ 
is the sum of all $x''_{\mathbf a, b}$. As $x$ respects the positivity submonoids, all the $x''_{\mathbf a, b}$ are non-negative and so
$q$ is a non-negative number.
Since $x$ is a chain map, we have that the differential of $x[i, r]$ is equal to
$x\partial[i, r] = x[r] - x[i] = \top - x[i]$ and therefore we must have $q=1$. Thus, only one of the $x''_{\mathbf a, b}$ is $1$ and all the others vanish. If this summand is $b_0 \star \top$, 
then by definition $\gamma x[i] = b_0\star\emptyset$, which proves that $\gamma x$ respects the augmentation.
\end{proof}

We record how the constructions $\gamma$ and $\beta$ behave with certain degeneracies.

\begin{lem}
\label{gammavsdegeneracies}
Let $x\colon O[m]\to D\star O[\top]$ be an augmented directed chain map with rank $0<r\leq m$ and corank $m>s\geq 0$.
The following identities hold.
\begin{enumerate}[leftmargin=*]
    \item $\gamma s_r^s\gamma x=\gamma x$, 
    \item $\beta^2x=\beta x$,
    \item $\gamma s_jx=s_j\gamma x$ for $j<r$,
    \item $\gamma s_j\gamma x=s_j\gamma x$ for $j<r$,
    \item $\beta s_jx=s_j\beta x$ for $j<r$,
    \item $\gamma s_jx=\gamma x$ for $j>r$.
\end{enumerate}
\end{lem}

\begin{proof}[Proof of~\cref{gammavsdegeneracies}]
We begin by proving $(1)$ and by observing that $(1)$ implies $(2)$.
The map $s^s_r \gamma x$ has source $O[m]$, has
rank $r>0$ and corank $s$, so that the map $\gamma s^s_r \gamma x$
has same rank $r>0$ and has corank $0$,
but has source $O[r]$.
It follows then immediately from the definition that
\[\gamma s^s_r \gamma[r] = \top = \gamma [r].\]
Let $B$ denote the basis of $D$. For any $[\mathbf a]$ in $B[r-1]_{q}$, 
we can write
\[
    x[\mathbf a, r] = 
        \sum_{ b \in B_{q+1}} x'_{\mathbf a, b}\cdot b\star\emptyset+ \sum_{ b \in B_q} x''_{\mathbf a, b}\cdot b\star \top,
\]
so that
\[
    \gamma x[\mathbf a] = \sum_{ b \in B_q} x''_{\mathbf a, b}\cdot b\star\emptyset
    \quad\text{and}\quad
    \gamma x[\mathbf a, r] = \sum_{ b \in B_{q+1}} x''_{\mathbf a, b}\cdot b \star \top.
\]
Regarding $[\mathbf a]$ as an element of $B[m]_{q}$, we get by definition that
\[
    s^s_r\gamma x[\mathbf a] = \gamma x[\mathbf a]
    \quad\text{and}\quad
    s^s_r \gamma x[\mathbf a, r] = \gamma x[\mathbf a, r].
\]
Applying $\gamma$, we can conclude by its definition that
$\gamma s^s_r \gamma x= \gamma x$ 
as augmented directed chain maps
from $O[r]$ to $D \star O[\top]$.
We remark that the normalizing factor $s^s_r$ plays no role in the proof,
so that in fact we have also proved the equality $\gamma \gamma x = \gamma x$.

We now prove $(3)$, $(4)$ and $(5)$, so we fix $0 \leq j < r$. Notice that
the property $(4)$ is a consequence of $(1)$ and $(3)$ as
\[
    \gamma s_j \gamma x = s_j \gamma \gamma x = s_j \gamma x.
\]
As for equality $(5)$, using $(3)$ and~\cref{IteratedSimpIdentities} we have
\[
    \beta s_j x = s^s_{r+1} \gamma s_j x = s_j s^s_r \gamma x = s_j \beta x.
\]
We now prove (3). Observe that
\[
    s_j \gamma x[r+1] = \gamma x[r] = \top = \gamma (s_j x) [r+1],
\]
as $s_j x$ has rank $r+1\le m+2$.
Consider an element $[\mathbf a]$ of $B[r]_{q}$. We have to examine two cases.
Suppose first that there is $0 \leq k < r$ such that
$a_k \leq j < j+1 < a_{k+1}$ or $a_k < j < j+1 \leq a_{k+1}$ and set
\[
 [\mathbf c] = s^j [\mathbf a] = [a_0, \dots,a_k, a_{k+1} - 1, \dots, a_q - 1].
\]
Notice that $a_q-1 < r$ and write
\[
    x[\mathbf c,r] = 
        \sum_{ b \in B_{q+1}} x'_{\mathbf c, b}\cdot b\star\emptyset + \sum_{ b \in B_q} x''_{\mathbf c, b}\cdot b\star \top.
\]
Then
\[
    s_j(\gamma x)[\mathbf a] = \gamma x[\mathbf c] = \sum_{ b \in B_q} x''_{\mathbf c, b}\cdot b\star\emptyset
\]
and, since $s_j x[\mathbf a, r+1] = x[\mathbf c, r]$,
\[
    \gamma(s_j x) [\mathbf a] = \gamma x[\mathbf c] = s_j(\gamma x)[\mathbf a].
\]
A completely similar argument works for a basis element of $O[r+1]$ of the form
$[\mathbf a, r+1]$. Suppose now we are in the case where $a_k = j < j+1 = a_{k+1}$,
so that we have $s_j(\gamma x)[\mathbf a] = \gamma x(0) = 0$
and, since $s_j x [\mathbf a] = x(0) = 0$, $\gamma(s_j x)[\mathbf a] = \gamma(0) = 0$.
The same holds for basis elements of the form $[\mathbf a, r+1]$, so that in all the cases
$s_j \gamma x = \gamma s_j x$.

For $r < j \leq m$, relation $(6)$ follows from the definition of the last vertex operator $\gamma$,
as $\gamma x$ is defined just on the basis elements $[a_0, \dots, a_q]$ of $O[m]$ with $a_q \leq r$.
\end{proof}

\begin{rmk}
\label{conecharacterization}
For a non-degenerate chain map $x\colon O[m]\to D\star O[\top]$, the following are equivalent:
\begin{itemize}[leftmargin=*]
    \item $x$ is \emph{$\beta$-invariant}, by which we mean that $\beta x=x$
    \item $x$ can be expressed as
    \[x=\rho(z)\colon O[m]\cong O[m-1]\star O[\top]\xrightarrow{z\star O[\top]} D\star O[\top].\]
\end{itemize}
If the equivalent conditions are met, we say that $x$ is \emph{conical}, and 
we otherwise say it is  \emph{non-conical}.
\end{rmk}

\subsubsection{The statements for the dual case}
\label{dualgamma}
We devote the rest of this subsection to providing guidelines on how to define the construction $\gamma x$ for an augmented directed chain map
\[x\colon O[m]\to O[\bot]\star D.\]
The details can be filled in with similar techniques to those employed in the original case.

Recall there are canonical inclusions 
\[D\cong O[-1]\star D\hookrightarrow O[\bot]\star D \hookleftarrow O[\bot]\star O[-1]\cong O[\bot]\]
of $D$ and $O[\bot]$ into the join $O[\bot]\star D$.
The analogous definitions of rank and corank are as follows.
\begin{defn}
Let $x\colon O[m]\to O[\bot]\star D$ be an augmented directed chain map.
\begin{itemize}[leftmargin=*]
    \item The \emph{rank} of $x$ is the number  of elementary $0$-chains $[a]$ in $O[m]$ such that $x[a]\in D\subseteq O[\bot]\star D$, namely
    \[r:=\#\{ 0\le a\le m:x[a]\in D\subseteq O[\bot]\star D\}.\]
\item The \emph{corank} $s$ is one less than the number of elementary $0$-chains $[a]$ in $O[m]$ such that $x[a]\in O[\bot]$, namely
\[s:=\#\{ 0\le a\le m:x[a]\in O[\bot]\subseteq O[\bot]\star D\}-1.\]
\end{itemize}
By definition, one always has $r+s=m$, with $0\le r\le m+1$ and $-1\le s\le m$.
\end{defn}

The construction $\gamma x$ in the dual context is as follows.

\begin{defn}\label{DefGammaBeta2}
Let $x\colon O[m]\to O[\bot]\star D$ be an augmented directed chain map with rank $r>0$ and corank $s\geq 0$ which is given on a basis element $[s, \mathbf a+s]$ for $[\mathbf a]\in B[m]_q$ by
\[
    x[s, \mathbf a+s]=\sum_{ b \in B_{q+1}} x'_{\mathbf a+s, b}\cdot \emptyset\star b+ \sum_{ b \in B_q} x''_{\mathbf a+s, b}\cdot \bot \star b.
\] 
\begin{itemize}[leftmargin=*]
    \item The \emph{first factor} of $x$ is the augmented directed chain map
    \[\gamma x\colon O[r]\to O[\bot]\star D\] given on basis elements by
\[\gamma x[s]= \bot\star\emptyset
\ , \quad 
\gamma x[s, \mathbf a+s]=\sum_{ b \in B_{q+1}}x''_{\mathbf a+s, b}\cdot \bot \star b
\ , \quad\gamma x[\mathbf a]=\sum_{ b \in B_{q}}x''_{\mathbf a+s, b} \cdot \emptyset\star b,
\]
for $[\mathbf a]$
in $B[r-1]_q$, and $q\ge0$.
\item The \emph{normalized first factor} of $x$ is the augmented directed chain map
\[\beta x\colon O[m]\to O[\bot]\star D\] given by
\[\beta x=s_{0}^{s}\gamma x.\]
\end{itemize}
\end{defn}

We now prove that $\gamma x$ is indeed a map of augmented direct chain complexes, building on some preliminary facts.

\begin{lem}
\label{gammaproperties2}
Let $D$ be a strong Steiner complex.
For $ x\colon O[m]\to O[\bot]\star D$ with rank $r>0$, the formula above defines an augmented directed chain map $\gamma x\colon O[r]\to O[\bot]\star D$ with rank $r>0$ and corank $0$.
\end{lem}

Again, it is handy to rewrite the formulas for $\gamma x$ in a slightly more abstract way. 
Let $J:= O[\bot]\star D$
and $q\ge0$.
Denote by $J(\bot)_q$ the subgroup of $J_q$ generated by all basis elements of the form $\bot\star b$ with $b\in B_{q-1}\subseteq D_{q-1}$,
and by $J(\neg \bot)_q$ the subgroup of $J_q$ generated by all basis elements of the form $b=\emptyset\star b$ with $b\in B_{q}\subseteq D_q$.
Denote by
\[\pr_{\bot}\colon J_q\to J(\bot)_{q}
\text{ and }\pr_{\neg\bot}\colon J_q\to J(\neg\bot)_{q}\]
the canonical projections. In this language, we observe that the formulas for $\gamma x$ can be rewritten to 
\[(\gamma x)[s]= \bot\star\emptyset
\, ,\;
(\gamma x)[s,\mathbf a+s]=\pr_{\bot}x[s,\mathbf a+s]
\, ,
\; (\gamma x)[\mathbf a+s]=\pr_{\neg\bot}\partial\pr_{\bot} x[s,\mathbf{a}+s].
\]
for $[\mathbf a]$
in $B[r-1]_q$, and $q\geq 0$.

\begin{lem}
\label{ProjAlmostChainMap2}
Let $q\geq 0$. With the notation above, the following holds: 
\[\partial_{q} \pr_{\neg\bot} \partial _{q+1} \pr_{\bot}=-\pr_{\neg\bot} \partial_{q}  \pr_{\bot}  \partial_{q+1}:J_{q+2}\to J(\neg\bot)_{q}.\]
\end{lem}

We record how the constructions $\gamma$ and $\beta$ behave with certain degeneracies, analogously to what we did in~\cref{gammavsdegeneracies}.

\begin{lem}
\label{gammavsdegeneracies2}
Let $x\colon O[m]\to O[\bot]\star D$ be an augmented directed chain map with rank $0<r\leq m$ and corank $0\leq s <m$.
The following identities hold. 
\begin{enumerate}[leftmargin=*]
    \item $\gamma s_0^s\gamma x=\gamma x$, 
    \item $\beta^2x=\beta x$,
    \item $\gamma s_{m-j}x=s_{r-j}\gamma x$ for $j<r$,
    \item $\gamma s_{r-j}\gamma x=s_{r-j}\gamma x$ for $j<r$,
    \item $\beta s_{m-j} x=s_{m-j}\beta x$ for $j<r$,
    \item $\gamma s_{m-j}x=\gamma x$ for $m \geq j > r$.
\end{enumerate}
\end{lem}

The two lemmas can be proven with arguments similar to~\cref{ProjAlmostChainMap,gammavsdegeneracies}, respectively.~\cref{ProjAlmostChainMap2} can be used to prove~\cref{gammaproperties2}, by adapting the argument for~\cref{gammaproperties}.

\begin{rmk}
\label{conecharacterizationv2}
For a non-degenerate chain map $x\colon O[m]\to O[\bot]\star D$, the following are equivalent:
\begin{itemize}[leftmargin=*]
    \item $x$ is \emph{$\beta$-invariant}, by which we mean that $\beta x=x$
    \item $x$ can be expressed as 
    \[x=\rho(z)\colon O[m]\cong O[\bot]\star O[m-1]\xrightarrow{O[\bot]\star z} O[\bot]\star D.\]
\end{itemize}
If the equivalent conditions are met, we say that $x$ is \emph{conical}, and 
we otherwise say it is  \emph{non-conical}.
\end{rmk}

\subsection{Interpolations of a simplex}

Like before, let $D$ be a strong Steiner complex. We introduce in detail constructions $w_jx$, $v_jx$ and $\alpha_j x$ based on simplices $x$ of the form (\ref{simplex}), which mostly generalize Steiner's from \cite[\textsection4]{SteinerEmbedding}, and indicate briefly in~\cref{dualinterpolator} how the analogous constructions work for the case of simplices $x$ of the form (\ref{dualsimplex}).

\begin{defn}
\label{closedformulas}
Let $x\colon O[m]\to D\star O[\top]$ be an augmented directed chain map
of rank $r$.
\begin{enumerate}[leftmargin=*]
\item For $0< j\le r$,
the \emph{$j$-th witness} of $x$ is the chain map \[w_jx\colon O[m+1]\to D\star O[\top]\] given by
\[w_jx=s_j^{r-j+1} d_{j}^{r-j}(x-\beta x)+s_{j-1}\beta x.\]
\item For $0\le j\le r$, the \emph{$j$-th interpolator} of $x$ is the chain map \[v_jx\colon O[m]\to D\star O[\top]\] given by
\[v_jx=s_j^{r-j} d_{j}^{r-j}(x-\beta x)+\beta x.\]
\item For $0\le j<r$, the \emph{$j$-th approximator} is the chain map \[\alpha_jx\colon O[j+1+s]\to D\star O[\top]\] given by
\[\alpha_jx= d_{j+1}^{r-j-1}(x-\beta x)+s_j d_{j+1}^{r-j}\beta x.\]
\end{enumerate}
\end{defn}

These chain maps are actually augmented and directed chain maps.

\begin{prop}\label{InterpolatingOperations}
Consider $x$ and $j$ as in the situation of~\cref{closedformulas}.
\begin{enumerate}[leftmargin=*]
    \item The $j$-th witness $w_jx\colon O[m+1]\to  D\star O[\top]$ is an augmented directed chain map.
    \item The $j$-th interpolator $v_jx\colon O[m]\to D\star O[\top]$ is an augmented directed chain map.
\item The $j$-th approximator  $\alpha_jx\colon O[j+1+s]\to D\star O[\top]$ is an augmented directed chain map.
\end{enumerate}
\end{prop}

To prove the proposition, we make use of certain recursive formulas for $v_j$ and $w_j$ based on the wedge construction. These are the original formulas used by Steiner for $v_j$ and $w_j$ in the case $D=O[n-1]$, appearing in \cite[Proposition~ 4.11]{SteinerUniversal} and \cite[Proposition~5.3]{SteinerUniversal}.

\begin{defn}[{\cite[Definition~3.5]{SteinerUniversal}}] 
\label{defwedge}
Let $C$ be a strong Steiner complex, and $x,x'\colon O[m]\to C$ augmented directed chain maps with $d_jx=d_{j+1}x'$. The \emph{$j$-th wedge} is the chain map
\[x\wedge_jx'\colon O[m+1]\to C\] given by
\[x\wedge_jx'=s_{j+1}x-s_j^2d_{j+1}x'+s_jx'=s_{j+1}x-s_j^2d_{j}x+s_jx'.\]
\end{defn}

\begin{lem}
Consider $x$ and $x'$ as in~\cref{defwedge}. The chain map $x\wedge_jx'\colon O[m+1]\to C$ 
is an augmented directed chain map.
\end{lem}

The case where $C = O[n]$, for $n\geq 0$, is \cite[Proposition~6.8]{SteinerAlgebraNerves}.

\begin{proof} 
To show that $x\wedge_j x'$ preserves the augmentation and is therefore augmented, it suffices to compute its value on elementary $0$-chains, which is as follows:
\begin{equation}
\label{verticesedge}
(x\wedge_j x')[i]=\left
\{\begin{array}{lll}
x[i]=x'[i]&\text{if } j<i\\
x[j]&\text{if } j=i\\
x[j+1]=x'[j]&\text{if } j+1=i\\
x'[j+1]&\text{if } j+2=i\\
x[i-1]=x'[i-1]&\text{if } j+2<i.\\
\end{array}
\right.
\end{equation}

To show that $x\wedge_j x'$ is directed, we observe that for any elementary chain $[\mathbf a]$, the chain $x\wedge_j x'[\mathbf a]$ is a linear combination of elementary chains with only non-negative coefficients.
\begin{itemize}[leftmargin=*]
    \item If $[\mathbf a]$ contains at most one between $j$ and $j+2$, then
    $x\wedge_jx'[\mathbf a]\in C^+$ is a positive chain
    because $d_j(x\wedge_jx')=x'$ and $d_{j+2}(x\wedge_jx')=x$.
        \item If $[\mathbf a]$ contains $j$, $j+1$ and $j+2$, then
    $x\wedge_jx'[\mathbf a]=0$.
        \item If $[\mathbf a]$ contains $j$ and $j+2$ but not $j+1$, then $s_j^2d_jx[\mathbf a]=d_jx[(s^j)^2\mathbf{a}]=0$, implying
    \[x\wedge_jx'[\mathbf a]=s_{j+1}x[\mathbf{a}]+s_jx'[\mathbf{a}]=x[s^{j+1}\mathbf a]+x'[s^j\mathbf a],\]
    which only contains non-negative summands because $x$ and $x'$ are directed.\qedhere
\end{itemize}
\end{proof}

 We now provide a geometric intuition on the wedge operation.
 
\begin{rmk}
\label{rmk:wedge_geometric}
    The condition $d_j x = d_{j+1}x'$ can be expressed by the commutativity of the
    following diagram of augmented directed chain complexes: 
    \[
        \begin{tikzcd}
        {O[m-1]} \ar[r, "{{O[d^{j+1}]}}"] \ar[d, "{{O[d^{j}]}}"{swap}]    & {O[m]} \ar[d, "{x'}"] \\
        {O[m]} \ar[r, "x"']                                         & C.
        \end{tikzcd}
    \]
    Notice that
    \[
        \begin{tikzcd}
        {O[m-1]} \ar[r, "{{O[d^{j+1}]}}"] \ar[d, "{{O[d^{j}]}}"{swap}]    & {O[m]} \ar[d, "{{O[d^{j}]}}"] \\
        {O[m]} \ar[r, "{{O[d^{j+2}]}}"']                                         & {O[m+1]}
        \end{tikzcd}
    \]
    commutes, so that we get a map
    \[
        \varphi^j \colon O[m] \aamalg{O[m-1]} O[m] \to O[m+1].
    \]
    We also have a map
    \[
        (x, x') \colon O[m] \aamalg{O[m-1]} O[m] \to C.
    \]
    It is possible to define a map
    \[
        \psi^j \colon O[m+1] \to O[m] \aamalg{O[m-1]} O[m]
    \]
    which is a retraction of $\varphi^j$ and such that $x \wedge_j x' = (x, x')   \psi^j$.
    The map $\psi^j$ is precisely ${O[d^{j+2}]} \wedge_j O[d^j]$.
\end{rmk}

\begin{defn}\label{DefIteratedWedge} Let $u\colon O[m-\ell+1]\to C$ and $v\colon O[m]\to C$ be augmented directed chain maps with $d_ju=d_{j+1}^{\ell}v$ for some $\ell\ge1$
and $0 \leq j \leq m-\ell$.
The \emph{$\ell$-fold $j$-th wedge} is the chain map
\[u\wedge^\ell_jv\colon O[m+1]\to C\] given by
\[
\begin{array}{ccl}
    u\wedge_j^\ell v = s_{j+1}^{\ell}u-s_j^{\ell+1}d_{j+1}^\ell v+s_jv.
\end{array}\]
\end{defn}

 Just as for the simple wedge, we can provide a geometric intuition of the $\ell$-fold wedge.
 
\begin{rmk}
\label{rmk:iterated_wedge_geometric}
    Let $u\colon O[m-\ell+1]\to C$ and $v\colon O[m]\to C$ be augmented directed chain maps
    with $d_ju=d_{j+1}^{\ell}v$ for some $\ell\ge1$
    and $0 \leq j \leq m-\ell$.
    The condition $d_j u = d_{j+1}^\ell v$
    can be expressed by the commutativity of the
    following diagram of augmented directed chain complexes:
    \[
        \begin{tikzcd}[column sep=3em]
        {O[m-\ell]} \ar[r, "{O[(d^{j+1})^\ell]}"] \ar[d, "{{O[d^{j}]}}"{swap}]   & {O[m]} \ar[d, "v"] \\
        {O[m-\ell+1]} \ar[r, "u"']                                             & C.
        \end{tikzcd}
    \]
    Notice that
    \[
        \begin{tikzcd}[column sep=3em]
        {O[m-\ell]} \ar[r, "{O[(d^{j+1})^\ell]}"] \ar[d, "{{O[d^{j}]}}"{swap}]   & {O[m]} \ar[d, "{{O[d^{j}]}}"] \\
        {O[m-\ell+1]} \ar[r, "{O[(d^{j+2})^\ell]}"']                             & {O[m+1]}
        \end{tikzcd}
    \]
    commutes, so that we get a map
    \[
        \varphi^j_\ell \colon O[m-\ell+1] \aamalg{O[m-\ell]} O[m] \to O[m+1].
    \]
    We also have a map
    \[
        (u, v) \colon O[m-\ell+1] \aamalg{O[m-\ell]} O[m] \to C.
    \]
    It is possible to define a map
    \[
        \psi^j_\ell \colon O[m+1] \to O[m-\ell+1] \aamalg{O[m-\ell]} O[m]
    \]
    which is a retraction of $\varphi^j_\ell$ and such that $u \wedge_j^\ell v = (u, v)   \psi^j_\ell$.
    The map $\psi^j_\ell$ is precisely $O[(d^{j+2})^\ell] \wedge^\ell_j O[d^j]$. 
\end{rmk}

The following lemma formalizes the fact that an $\ell$-fold wedge can be seen as an iterated ordinary wedge. The intuition is that the description of the $\ell$-fold wedge as a retraction from~\cref{rmk:iterated_wedge_geometric} is built up as an iterated single step retractions as in~\cref{rmk:wedge_geometric}, that correspond precisely to instances of the ordinary wedge.

The following lemma, of which we omit the details, can be proven borrowing some ideas from \cite[Proposition\ 4.10]{SteinerUniversal}.

\begin{lem}
\label{lem:iterated_wedge}
Consider $u$ and $v$ as in~\cref{DefIteratedWedge}. The chain map $u\wedge^{\ell}_j v\colon O[m+1]\to C$ can be written as an iterated wedge:
\[u\wedge_j^\ell v = u\wedge_jd_{k+1}^{\ell-1}v\wedge_jd_{j+1}^{\ell-2}v\wedge_j\dots\wedge_j d_{j+1}v\wedge_j v,\]
and is in particular an augmented directed chain map.
\end{lem}

\begin{lem}\label{VvsW}
Consider $x$ and $j$ as in the situation of~\cref{closedformulas}. Then the following hold.
\begin{enumerate}[leftmargin=*] 
    \item $v_jx= d_j w_jx$.
    \item $d_j\alpha_jx=d_{j+1}^{r-j}v_jx$, so in particular $\alpha_jx\wedge^{r-j}_jv_jx$ exists.
    \item $w_{j+1}x=\alpha_jx\wedge^{r-j}_jv_jx$. 
\end{enumerate}
\end{lem}

\begin{proof}[Proof of~\cref{VvsW}]
(1) follows from simplicial identities.

For (2), on the one hand we have that
\[
\begin{array}{ccll}
     d_j\alpha_j x&=&d_j(d_{j+1}^{r-j-1}(x-\beta x)+s_j d_{j+1}^{r-j}\beta x)&\text{\cref{closedformulas}}  \\
     &=&d_{j}^{r-j}(x-\beta x)+ d_{j+1}^{r-j}\beta x,&\text{\cref{IteratedSimpIdentities}}
\end{array}
\]
and the other hand we have that
\[
\begin{array}{ccll}
     d_{j+1}^{r-j}v_jx&=&d_{j+1}^{r-j}(s_j^{r-j} d_{j}^{r-j}(x-\beta x)+\beta x)&\text{\cref{closedformulas}}  \\
     &=&d_{j}^{r-j}(x-\beta x)+d_{j+1}^{r-j}\beta x,&\text{\cref{SimpIdentities}}
\end{array}
\]
as desired.

For (3),
we have that
\[\begin{array}{rll}
   \alpha_{j}x\wedge^{r-j}_{j}v_{j}x  &  =  s_{j+1}^{r-j}\alpha_{j}x-s_{j}^{r-j+1}d_{j+1}^{r-j}v_{j}x+s_{j}v_{j}x &\text{\cref{DefIteratedWedge}}\\
   &=s_{j+1}^{r-j}d_{j+1}^{r-j-1}(x-\beta x)+s^{r-j}_{j+1}s_{j}d_{j+1}^{r-j}\beta x&\text{\cref{closedformulas}}\\
   &\hphantom{=}-s_{j}^{r-j+1}d_{j+1}^{r-j}(s_{j}^{r-j}d_{j}^{r-j}(x-\beta x)+\beta x)&\\
   &\hphantom{=}+s_{j}(s_{j}^{r-j}d_{j}^{r-j}(x-\beta x)+\beta x)&\\
   &=s_{j+1}^{r-j}d_{j+1}^{r-j-1}
   (x-\beta x) +s_{j}^{r-j+1}d_{j+1}^{r-j}\beta x &\text{\cref{IteratedSimpIdentities}}\\
&\hphantom{=}-s_{j}^{r-j+1}d_{j}^{r-j}(x-\beta x)-s_{j}^{r-j+1}d_{j+1}^{r-j}\beta x& \text{\cref{SimpIdentities}}\\
&\hphantom{=}+s_{j}^{r-j+1}d_{j}^{r-j}(x-\beta x)+s_{j}\beta x&\\
&=s_{j+1}^{r-j}d_{j+1}^{r-j-1}(x-\beta x)+s_{j}\beta x&\\
&= w_{j+1}x&\text{\cref{closedformulas}}
\end{array}\]  
as desired.
\end{proof}

We can now prove the proposition.

\begin{proof}[Proof of~\cref{InterpolatingOperations}]
We want to show that the chain maps $w_jx$, $v_jx$ and $\alpha_jx$ preserve 
the augmentation and the positivity submonoids, i.e., are directed.
Using Lemmas~\ref{VvsW} and~\ref{lem:iterated_wedge} we see that we proceed
inductively: knowing that $\alpha_jx$ and $v_jx$ are augmentation-preserving and directed,
so is $w_{j+1}x$, as this is a $(r-j)$-fold wedge of the former maps.
Moreover, $v_{j+1}x$ is a face of $w_{j+1}$ and hence augmentation-preserving and directed itself. Thus, it is enough to show that $v_0x$ and the $\alpha_jx$, for $0\leq j < r$,
are augmented directed chain maps.

Let $[\mathbf a] = [a_0, \dots, a_q]$ be an element in the basis of $O[m]$.
Applying $s^0$ for $r$ times we get $[a_0-r, \dots, a_q-r]$ where here
by convention $a_i - r = 0$ whenever $a_i < r$. If we
now apply $d^0$ for $r$ times we get a chain $[a_0', \dots, a_q']$
where each component is greater or equal than $r$ (with possibly repeated components).
Therefore $s^r_0d^r_0(x -\beta x)[\mathbf a]$ always maps to $0$ and we obtain
that $v_0x = \beta x$, which is an augmented directed chain map by~\cref{gammaproperties}.

Let us turn to $\alpha_j x$. By inspecting the defining identity we find that
\[\alpha_jx[i]=\left\{\begin{array}{lllll}
    x[i] & 0\leq i\leq j \\
        \beta x[j] & i=j+1\\
    \top & j+1< i \leq j+s+1\\
\end{array}\right.\]
This entails that $\alpha_jx$
is augmentation-preserving.
The argument that Steiner gives to prove the first part of \cite[Proposition~4.8]{SteinerUniversal} generalizes verbatim to prove that $\alpha_jx$ is a directed map, thereby finishing
the proof of the lemma.
\end{proof}

The following lemma records certain features of $\alpha_jx$.

\begin{lem}
\label{steiner4.8}
Let $x\colon O[m]\to D\star O[\top]$ be an augmented directed chain map
that
has rank $r\le m$
and corank $s$, and let $0\le j<r$. Then
$\alpha_jx$ has rank $j+2$ and corank $s-1$.
\end{lem}

\begin{proof}
The argument that Steiner gives to prove the second part of \cite[Proposition~4.8]{SteinerUniversal} generalizes verbatim to prove the current statement.
\end{proof}

The following lemma records the value of certain constructions on $d_ix$.

\begin{lem}
\label{steiner5.1}
\label{betaandfaces}
Let $x\colon O[m]\to D\star O[\top]$ be an augmented directed chain map with rank $r$. Let $0\le i< r$. The following hold.
\begin{itemize}[leftmargin=*]
    \item $d_ix$ has rank $r-1$.
    \item $\gamma d_ix=d_i\gamma x$.
    \item $\beta d_ix=d_i\beta x$.
    \item
    $ \alpha_jd_ix=\left\{\begin{array}{lll}
       \alpha_jx & \text{for $0\le j<i$} \\
       d_i\alpha_{j+1}x  & \text{for $i\le j<r-1$.} 
    \end{array}
    \right.$
\end{itemize}
\end{lem}

\begin{proof}
If $r=0$ or $r=1$, the first point is trivial. Otherwise by~\cref{remarkrank}
$x[r] = \top$ and $x[r-1] \in D_0$ imply
$d_ix[r-1] = \top$ and $d_ix[r-2] \in D_0$, so that $d_ix$ has rank $r-1$.
The argument that Steiner gives to prove \cite[Proposition~5.1]{SteinerUniversal} generalizes verbatim to prove the second and the fourth point.

Finally, we prove the third point. By the second point
we have that $\gamma d_i x= d_i\gamma x$. Taking the rank considerations from the first point
into account as well as the simplicial identities from~\cref{SimpIdentities}, we obtain that 
\[\beta d_ix=s_{r-1}^s \gamma d_ix= s_{r-1}^sd_i \gamma x=d_is_{r}^s\gamma x= d_i\beta x.\qedhere\]\end{proof}

The following lemma records certain features of $w_jx$ and $v_jx$.

\begin{lem}
\label{steiner5.3}
Let $x\colon O[m]\to D\star O[\top]$ be an augmented directed chain map that
has rank $r\le m$,
and $0<j\le r$. Then $w_jx$ and $v_jx$ have the same corank as $x$.
\end{lem}

\begin{proof}
The argument that Steiner gives to prove \cite[Proposition~5.3]{SteinerUniversal} generalizes verbatim to prove the current statement.
\end{proof}

The following justifies the terminology of ``interpolator'' for $v_jx$.

\begin{lem}\label{vstabilize}
Consider $x$ as in the situation of~\cref{closedformulas}. The following hold.
\begin{enumerate}[leftmargin=*]
\item $v_0x=\beta x$. 
\item $v_rx=x$.
\item if $v_jx=x$ then $v_{j+1}x=x$.
\end{enumerate}
\end{lem}

This means that the sequence $(v_jx)_{j\ge0}$ starts at $\beta x$, will at some point become $x$, and never change anymore. So it makes sense to regard $v_jx$ as an interpolation between $\beta x$ and $x$ in $\adch(O[m], D\star O[\top])$.

\begin{proof}
Statement (1) is a consequence of the simplicial identities and of~\cref{totallydegenerate} and was already treated while proving~\cref{InterpolatingOperations},
and Statement (2) can be obtained by direct computation.

For (3), we assume $v_jx=x$ with $j<r$ and will show then $v_{j+1}x=x$. Since $x=v_jx$, by definition of $v_jx$ we obtain that
\begin{equation}
\label{ConsequenceIndexj}
x-\beta x=v_jx-\beta x=s_j^{r-j}d_j^{r-j}(x-\beta x).
\end{equation}
This implies that
\begin{equation}
\label{ConsequenceAlpha}
\begin{array}{ccll}
   \alpha_jx  & = &  d_{j+1}^{r-j-1}(x-\beta x)+s_j d_{j+1}^{r-j}\beta x&\text{\cref{closedformulas}}\\
   & = &  d_{j+1}^{r-j-1}s_j^{r-j}d_j^{r-j}(x-\beta x)+s_j d_{j+1}^{r-j}\beta x&\text{Formula \eqref{ConsequenceIndexj}}\\
     & =& s_jd_j^{r-j}(x-\beta x)+s_j d_{j+1}^{r-j}\beta x&\text{\cref{SimpIdentities}}\\
     &=& s_jd_j (d_{j+1}^{r-j-1}(x-\beta x) + s_jd_{j+1}^{r-j}\beta x)&\text{\cref{SimpDegLinear}}\\
     &=& s_jd_j\alpha_jx.&\text{\cref{closedformulas}}
\end{array}
\end{equation}
It follows that
\[\begin{array}{ccll}
v_{j+1}x&=&d_{j+1}w_{j+1}x&\text{\cref{VvsW}}\\
&=&d_{j+1}(\alpha_jx\wedge^{r-j}_jx)&\text{\cref{VvsW}}\\
&=&d_{j+1}(s_{j+1}^{r-j}\alpha_jx-s_j^{r-j+1}d_{j+1}^{r-j}x+s_{j}x)&\text{\cref{DefIteratedWedge}}\\
&=&s_{j+1}^{r-j-1}\alpha_jx-s_j^{r-j}d_{j+1}^{r-j}x+x&\text{\cref{SimpIdentities}}\\
&=&s_{j+1}^{r-j-1}(s_jd_j\alpha_jx)-s_j^{r-j}d_{j+1}^{r-j}x+x&\text{Formula \eqref{ConsequenceAlpha}}\\
&=&s_{j+1}^{r-j-1}(s_jd_j\alpha_jx)-s_j^{r-j}d_{j+1}^{r-j}v_jx+x&\text{Assumption}\\
&=&s_{j+1}^{r-j-1}(s_jd_j\alpha_jx)-s_j^{r-j}d_j\alpha_jx+x&\text{\cref{VvsW}}\\
&=&x&\text{\cref{SimpIdentities}},
\end{array}\]
as desired.
\end{proof}

\subsubsection{The statements for the dual case}
\label{dualinterpolator}
We devote the rest of this section to providing guidelines on how to define the constructions $w_jx$, $v_jx$ and $\alpha_jx$ for an augmented directed chain map
\[x\colon O[m]\to O[\bot]\star D.\]
The details can be filled in with similar techniques to those employed in the original case.

\begin{defn}
\label{closedformulas2}
Let $x\colon O[m]\to O[\bot]\star D$ be an augmented directed chain map
of rank $r$.
\begin{enumerate}[leftmargin=*]
\item For $0< j\le r$,
the \emph{$j$-th witness} of $x$ is the chain map \[w_jx\colon O[m+1]\to O[\bot]\star D\] given by
\[w_jx=s_{m-r}^{r-j+1} d_{m-r+1}^{r-j}(x-\beta x)+s_{m-j+1}\beta x.\]
\item For $0\le j\le r$, the \emph{$j$-th interpolator} of $x$ is the chain map \[v_jx\colon O[m]\to  O[\bot]\star D\] given by 
\[v_jx=s_{m-r}^{r-j} d_{m-r+1}^{r-j}(x-\beta x)+\beta x.\]
\item For $0\le j<r$, the \emph{$j$-th approximator} is the chain map \[\alpha_jx\colon O[s+j+1]\to O[\bot]\star D\] given by
\[\alpha_jx= d_{m-r+1}^{r-j-1}(x-\beta x)+s_{m-r} d_{m-r}^{r-j}\beta x.\]
\end{enumerate}
\end{defn}

These chain maps are actually augmented and directed chain maps.

\begin{prop}\label{InterpolatingOperations2}
Consider $x$ and $j$ as in the situation of~\cref{closedformulas2}.
\begin{enumerate}[leftmargin=*]
    \item The $j$-th witness $w_jx\colon O[m+1]\to O[\bot]\star D$ is an augmented directed chain map.
    \item The $j$-th interpolator $v_jx\colon O[m]\to O[\bot]\star D$ is an augmented directed chain map.
\item The $j$-th approximator  $\alpha_jx\colon O[s+j+1]\to O[\bot]\star D$ is an augmented directed chain map.
\end{enumerate}
\end{prop}

Like in the original case treated earlier, to prove the proposition we make use of certain recursive formulas for $v_j$ and $w_j$ based on a variant of the wedge construction.

\begin{defn}\label{DefDualIteratedWedge}
Let $u\colon O[m-\ell+1]\to C$ and $v\colon O[m]\to C$ be augmented directed chain maps
with $d_{m-\ell-j+1} u= d_{m-\ell-j}^{\ell}v$ for some $1 \leq\ell \leq m$
and $0 \leq j \leq m-\ell$.
The \emph{dual $\ell$-fold $j$-th wedge} is the chain map
\[
    u \dwedgevar{\ell}{j} v\colon O[m+1]\to D
\]
given by
\[
\begin{array}{ccl}
    u \dwedgevar{\ell}{j} v =
    s_{m-\ell-j}^{\ell} u - s_{m-\ell-j}^{\ell+1}d_{m-\ell-j}^\ell v+s_{m-j} v.
\end{array}\]
\end{defn}

As for the iterated $\ell$-fold wedge, also the dual iterated
$\ell$-fold wedge can be seen as an iterated ordinary wedge,
with a similar underlying geometric intuition. 
Notice that if $u$ and $v$ are as in the previous definition,
then we have the identities $d_{m-j+1}(u \dwedgevar{\ell}{j} v)=v$ and $d^{\ell}_{m-j-\ell}(u \dwedgevar{\ell}{j} v)=u$.

\begin{lem}
\label{lem:dual_iterated_wedge}
Consider $u$ and $v$ as in~\cref{DefDualIteratedWedge}. The chain map $u\dwedgevar{\ell}{j} v\colon O[m+1]\to C$ can be written as an iterated wedge: 
\begin{multline*}
    u \dwedgevar{\ell}{j} v =
        v \wedge_{m-j-1} ( d_{m-j-1}v \wedge_{m-j-2} (\dots  \\
        \dots \wedge_{m-\ell-j+2}  (d_{m-\ell-j+2}^{\ell-2}v
        \wedge_{m-\ell-j+1} (d_{m-\ell-j+1}^{\ell-1}v \wedge_{m-\ell-j} u)\dots),
\end{multline*}
and is in particular an augmented directed chain map.
\end{lem}

The following lemma is the analog of~\cref{VvsW}.

\begin{lem}\label{VvsW2}
Consider $x$ and $j$ as in the situation of~\cref{closedformulas2}. 
Then the following hold.
\begin{enumerate}[leftmargin=*] 
    \item $v_jx= d_{m+1-j} w_jx$. 
    \item $d_{m-r+1}\alpha_jx=d_{m-r}^{r-j}v_jx$, so in particular $\alpha_jx \dwedgevar{r-j}{j}v_jx$ exists.
    \item $w_{j+1}x=\alpha_jx\dwedgevar{r-j}{j}v_jx$. 
\end{enumerate}
\end{lem}

The following lemma is the analog of~\cref{steiner4.8}.

\begin{lem}
\label{steiner4.8v2}
Let $x\colon O[m]\to O[\bot]\star D$ be an augmented directed chain map with rank $r$ and corank $s\geq 0$, and let $0\le j<r$. Then $\alpha_jx$ has rank $j+2$ and corank $s-1\geq 0$.
\end{lem}

The following lemma, which is the analog of~\cref{steiner5.1}, records the value of certain constructions on $d_ix$.

\begin{lem}
\label{steiner5.1v2}
Let $x\colon O[m]\to O[\bot]\star D$ be an augmented directed chain map with rank $r$. Let $s=m-r< i\leq m$.
The following hold.
\begin{itemize}[leftmargin=*]
    \item $d_ix$ has rank $r-1$.
    \item $\gamma d_ix=d_{i+r-m}\gamma x$.
    \item $\beta d_ix=d_i\beta x$.\label{betaandfacesv2}
    \item
    $ \alpha_jd_ix=\left\{\begin{array}{lll}
       \alpha_jx & \text{for $0\le j<m-i$} \\
       d_{i-r+j+2}\alpha_{j+1}x  & \text{for $m-i\le j<r-1$.} 
    \end{array}
    \right.$
\end{itemize}
\end{lem}

The following lemma, which is the analog of~\cref{steiner5.3} records certain features of $w_jx$ and $v_jx$.

\begin{lem}
\label{steiner5.3v2}
Let $x\colon O[m]\to O[\bot]\star D$ be an augmented directed chain map with rank $r$ and corank $s\geq 0$, and $0<j\le r$. Then $w_jx$ and $v_jx$ have corank~$s$.
\end{lem}

The following lemma is the analog of~\cref{vstabilize}.

\begin{lem}\label{vstabilizev2}
Consider $x$ as in the situation of~\cref{closedformulas2}. The following hold.

\begin{enumerate}[leftmargin=*]
\item $v_0x=\beta x$. 
\item $v_rx=x$.
\item if $v_jx=x$ then $v_{j+1}x=x$.
\end{enumerate}
\end{lem}

\section{Pairing of simplices in nerves of cones}

The goal of this section is to build a pairing on the set of simplices of $\NRS(\cD\star\cO[\top])$ that do not belong to $\NRS\cD\star\Delta[\top]$, where $\cD=\nu D$ is a strong Steiner $\omega$-category. This pairing will play a crucial role in the proof of the main theorem in~\cref{proofmaintheorem}.

Each ``suspect'' simplex will be paired up with a ``non-suspect'' simplex, which we now define.

\begin{defn}
\label{def:SuspectIndex}
The \emph{suspect index} of a chain map $x\colon O[m]\to D\star O[\top]$ is the minimal integer $p$ for which $v_px=x$. Alternatively, using~\cref{vstabilize}, it is any integer $p$ for which
\[v_px=x\quad\text{ and }\quad v_{p-1}x\neq x.\]
Note that $0\le p\le r$.
\end{defn}

The intuition is that, roughly speaking, the suspect index $p$ of a simplex $x$ counts how many factors are needed to express $x$ as a composite of generating cells.

\begin{rmk}
\label{SuspectIndexConical}
A non-degenerate augmented directed chain map $x\colon O[m]\to D\star O[\top]$ has suspect index $0$ if and only if it is conical,
as can be deduced combining~\cref{conecharacterization,vstabilize}.
\end{rmk}

\begin{defn}
\label{suspect}
Let $y\colon O[M]\to D\star O[\top]$
be a non-degenerate non-conical augmented directed chain map of suspect index $p$. The chain map $y$ is said to be a \emph{suspect} simplex of $\NRS(\cD\star[0])$ if $v_{p-1}y$ is degenerate at $p-1$, and it is said to be a \emph{non-suspect} simplex of $\NRS(\cD\star[0])$ otherwise.
\end{defn}

\begin{lem} \label{SuspectRankSuspectIndex}
Let $y\colon O[M]\to D\star O[\top]$ be a suspect simplex with rank $R$ and suspect index $p$. Then
\begin{enumerate}[leftmargin=*]
    \item $R> 1$,
    \item $R> p$.
\end{enumerate}
\end{lem}

\begin{proof}
For (1), we prove that $R\neq0,1$.
If $R=0$, then $y$ would be conical by~\cref{SuspectIndexConical}, and in particular non-suspect, contradicting the assumption.
If $R=1$, then $p=0$ or $p=1$. If $p=0$ then
\[y=v_0y=\beta y=s_1^{M-1}\gamma y\]
would be conical, contradicting the hypothesis.
If $p=1$, on the one hand by~\cref{suspect} $v_0y$ is degenerate at $0$, which implies
\[v_0y[0]=v_0y[1]\]
On the other hand (using~\cref{vstabilize}) we have
\[v_0y=\beta y=s_1^{M-1}\gamma y,\]
which implies
\[v_0y[1]=v_0y[2]=\dots=v_0y[M].\]
As a consequence, all vertices of $\beta y$ coincide, meaning the rank of $v_0y$ is $0$, meaning that $\beta y$ is a degeneracy of a $0$-simplex using~\cref{totallydegenerate}. In particular, also $\gamma y$ must be totally degenerate, and therefore has rank $0$.
This contradicts~\cref{gammaproperties}, as we know that the rank of $y$ is $1$.

For (2), we show that $p=R$ leads to a contradiction. The assumption that $p=R$ means that we have \[y=v_Ry\quad\text{ and }\quad v_{R-1}y=s_{R-1}z\]
with $z=d_{R-1}v_{R-1}y$.
It follows that
\[
\begin{array}{ccll}
     y&=&v_Ry\\
     &=&d_{R}w_Ry&\text{\cref{VvsW}}\\
     &=&d_{R}(\alpha_{R-1}y \wedge_{R-1} v_{R-1}y)& \text{\cref{VvsW}} \\
     &=&d_{R}(\alpha_{R-1}y \wedge_{R-1} s_{R-1}z)& \text{Assumption}\\
     &=&d_{R}(s_R\alpha_{R-1}y-s_{R-1}^2d_Rs_{R-1}z+s_{R-1}^2z)&\text{\cref{defwedge}}\\
     &=& \alpha_{R-1}y.&\text{\cref{SimpIdentities}}
\end{array}
\]
By \cite[Proposition~4.8]{SteinerUniversal}, the corank of $\alpha_{R-1}y$ must be strictly smaller than the corank of $y$. This contradicts the fact that $y=\alpha_{R-1}y$.
\end{proof}

\subsubsection{The statements for the dual case} 

One can develop a similar argument to the one from the previous section to build a pairing on the set of simplices of $\NRS(\cO[\bot]
\star\cD)$ that do not belong to $\Delta[\bot]\star\NRS\cD$. This can be done using the following variants of the definitions we just introduced.

\begin{defn}
\label{def:SuspectIndexv2}
The \emph{suspect index} of a chain map $x\colon O[m]\to O[\bot]\star D$ is the minimal integer $p$ for which $v_px=x$. Alternatively, using~\cref{vstabilizev2}, it is any integer $p$ for which
\[v_px=x\quad\text{ and }\quad v_{p-1}x\neq x.\]
Note that $0\le p\le r$.
\end{defn}

\begin{rmk}
A non-degenerate augmented directed chain map $x\colon O[m]\to O[\bot] \star D$ has suspect index $0$ if and only if it is conical, as can be deduced combining ~\cref{conecharacterizationv2,vstabilizev2}.
\end{rmk}

\begin{defn}
\label{suspectv2}
Let $y\colon O[M]\to O[\bot] \star D$
be a non-degenerate non-conical  augmented directed chain map of suspect index $p$. The chain map $y$ is said to be a \emph{suspect} simplex of $\NRS(\cO[\bot]\star\cD)$ if $v_{p-1}y$ is degenerate at $M-p$, and it is said to be a \emph{non-suspect} simplex of $\NRS(\cO[\bot]\star\cD)$ otherwise.
\end{defn}

\begin{lem} \label{SuspectRankSuspectIndexv2}
Let $y\colon O[M]\to O[\bot] \star D$ be a suspect simplex with rank $R$ and suspect index $p$. Then
\begin{enumerate}[leftmargin=*]
    \item $R> 1$,
    \item $R> p$.
\end{enumerate}
\end{lem}

\subsection{Features of \texorpdfstring{$w_px$}{w_p x}}
The goal of this subsection is to collect and prove all the relevant properties and features of $w_px$ we need, when $x$ is a non-suspect simplex of suspect index $p$ of the form

\begin{equation}
\label{bothsimplices}    
x\colon O[m]\to D\star O[\top]\quad\text{ or }\quad x\colon O[m]\to O[\bot]\star D
\end{equation}
We provide the proofs for the first case, and leave the details of the second case to the reader, indicating the relevant statements  in~\cref{dualw}.

\begin{prop}
\label{wpwelldef}
Let $x\colon O[m]\to D\star O[\top]$ be a non-degenerate non-suspect simplex with suspect index $p$. The following hold.
\begin{enumerate}[leftmargin=*, label={\normalfont (W\arabic*)}, ref={\normalfont (W\arabic*)}]
\item \label{SuspectPartnerDim} If $x$ has dimension $m$, the simplex $w_px$ has dimension $m+1$.
\item \label{SuspectPartnerDeg} The simplex $w_px$ is non-degenerate.
\item \label{SuspectPartnerTerm} If $x$ has non-negative corank, the simplex $w_px$ has non-negative corank.
\item \label{SuspectPartnerCone}
If $x$ is non-conical, then $w_px$ is also non-conical.
\item \label{SuspectPartnerRank} If $x$ has rank $r$, the simplex $w_px$ has rank $r+1$.
\item \label{SuspectPartnerIndex}
The simplex $w_px$ has suspect index $p$.
\item \label{SuspectPartnerSuspect} The simplex $w_px$ is a suspect simplex.
\end{enumerate}
\end{prop}

\begin{proof}[Proof of \ref{SuspectPartnerDim}]
This is by definition of the witness construction $w_p$.
\end{proof}

A careful but straightforward analysis shows the following lemma, of which we omit the proof, which will be used to prove \ref{SuspectPartnerDeg}.

\begin{lem}\label{DegWedge}
Let $x,x'$ be augmented directed chain maps of the form \eqref{bothsimplices}
with $d_ix=d_{i+1}x'$. The wedge simplex $x\wedge_i x'$ is degenerate at $i+1$ if and only if $x'=s_id_ix'$.
    
\end{lem}

\begin{proof}[Proof of \ref{SuspectPartnerDeg}]
Let $x\colon O[m]\to D\star O[\top]$ be an augmented directed chain map.
We can write
\begin{equation}
\label{auxiliaryformulaeta}
w_px=\alpha_{p-1}x \wedge^{r-p+1}_{p-1} v_{p-1}x=\eta_{p-1}x\wedge_{p-1}v_{p-1}x,
\end{equation}
where
\[
\eta_{p-1}x:=\alpha_{p-1}x \wedge_{p-1} d_p^{r-p}v_{p-1}x \wedge_{p-1} \ldots \wedge_{p-1}d_pv_{p-1}x.
\]
Assume by contradiction that $w_px$ were degenerate at $k$. If $k\notin \{p, p-1\}$, then $x=v_px=d_pw_px$ would also be degenerate, contradicting the assumption. If $k=p$, by~\cref{DegWedge} we would have that $v_{p-1}x$ is degenerate at $p-1$, which is excluded by assumption. If $k=p-1$, then $x=v_px=d_pw_px=d_{p-1}w_px=v_{p-1}x$, contradicting the minimality of $p$.
\end{proof}

\begin{proof}[Proof of \ref{SuspectPartnerTerm}]
This is a consequence of~\cref{steiner5.3}.
\end{proof}

\begin{proof}[Proof of \ref{SuspectPartnerCone}]
Let $x\colon O[m]\to D\star O[\top]$ be an augmented directed chain map.
If we assume by contradiction that $w_px=\rho z$ is conical, then using the fact that $p<m+1$ we would have that
\[x=v_px=d_pw_px=d_p\rho z=\rho d_pz,\]
which contradicts the fact that $x$ is non-conical.
\end{proof}

\begin{proof}[Proof of \ref{SuspectPartnerRank}]
Let $x\colon O[m]\to D\star O[\top]$ be an augmented directed chain map.
Using~\cref{gammaproperties,conecharacterization,steiner5.3}, we can conclude that the rank of $w_px$ is $r+1$ and the corank is $s$.
\end{proof}

To prove \ref{SuspectPartnerIndex} and \ref{SuspectPartnerSuspect} we use the following lemma.

\begin{lem}
 \label{vjidempotent}
 Let $x\colon O[m]\to D\star O[\top]$ be an augmented directed chain map of rank $r$ and corank $s=m-r$. Then the following holds:
 \begin{enumerate}
     \item For all $1\leq j \leq r$, we have $\beta w_jx=s_{j-1}\beta x$.
     \item For all $0 \leq j \leq r$, we have $\beta v_jx=\beta x$. 
     \item For all $0 \leq j \leq r$, we have $v_j v_jx=v_j x$.  
 \end{enumerate}
 \end{lem}

The proof of this lemma makes use of the following observation.

\begin{rmk} 
\label{cosimplicialrmk}
Let $s^j\colon[m]\to[m-1]$ and $d^j\colon[m]\to[m+1]$ with usual ranges denote coface and codegeneracy maps in $\Delta$.
The following identities hold.
\begin{itemize}[leftmargin=*]
\item
For any $0\leq a \leq m+k$, $0\leq j \leq m$, $0<k$ have
\[
(s^j)^{k}(a)=\begin{cases}
a, \mbox{ if } 0 \leq a \leq j,\\
j, \mbox{ if }j \leq a \leq j+k,\\
a-k, \mbox{ if }j+k \leq a \leq m+k.
\end{cases}
\]
\item For any $0\leq a \leq m$, $0\leq j \leq m$, $0<k$ have
\[
(d^j)^{k}(a)=\begin{cases}
a, \mbox{ if } 0 \leq a \leq j-1,\\
a+k, \mbox{ if }j \leq a \leq m.
\end{cases}
\]
\item 
Combining the two yields for $0\leq a \leq m$, $0<k$, $0\leq j \leq m-k-1$:
\[
(d^j)^{k}(s^j)^{k+1}(a)=\begin{cases}
a, \mbox{ if } 0 \leq a \leq j-1,\\
j+k, \mbox{ if }j \leq a \leq j+k+1,\\
a-1, \mbox{ if }j+k+1 \leq a \leq m.
\end{cases}
\]
\end{itemize}
\end{rmk}

\begin{proof}[Proof of~\cref{vjidempotent}]
As a first auxiliary step, we show that
\begin{equation}
    \label{auxiliaryformula}
\gamma w_jx=\gamma s_{j-1}\beta x.
\end{equation}
For this, by definition of $\gamma$ and using that $\gamma\beta=\gamma$ by~\cref{gammavsdegeneracies}(1), it is enough to argue that the chain $s_j^{r-j+1} d_{j}^{r-j}(x-\beta x)[\mathbf{a}, r+1]$, where $\mathbf{a}\in B[r]$,
does not contain terms of the form $b\star \top$.

If we set $\widetilde{\mathbf a}=(d^p)^{r-p}(s^p)^{r-p+1}[\mathbf{a}]$, we obtain that
\[\begin{array}{rll}
  s_j^{r-j+1} d_{j}^{r-j}(x-\beta x)[\mathbf{a}, r+1]=&(x-\beta x)[(d^j)^{r-j}(s^j)^{r-j+1}(\mathbf{a}), r]\\
  =&(x-\beta x)[\widetilde{\mathbf a},r]\\
  =&x[\widetilde{\mathbf a},r]-(\beta x)[\widetilde{\mathbf a},r]\\
  =&x[\widetilde{\mathbf a},r]-(\gamma x)[(s^r)^s(\widetilde{\mathbf a}),r]&\text{\cref{cosimplicialrmk}}\\
   =&x[\widetilde{\mathbf a},r]-(\gamma x)[\widetilde{\mathbf a},r]&\text{\cref{cosimplicialrmk}}
\end{array}\]
which by definition does not contain any summands of the form $b\star \top$, as desired.
Thus we have shown $\gamma w_jx =\gamma s_{j-1}\beta x$ for all $1\leq j \leq r$. 

Next, recalling that the rank of $\beta x$ is exactly $r$ (explained in~\cref{DefGammaBeta}) we show that 
$\gamma w_jx=s_{j-1}\gamma x$. For this, we have that
\[
\begin{array}{ccll}
\gamma w_j x&= &\gamma s_{j-1}\beta x &\text{Previous step} \\
&=& s_{j-1}\gamma \beta x & \text{\cref{gammavsdegeneracies}(2)}\\
&=& s_{j-1} \gamma s_r^s \gamma x & \text{\cref{DefGammaBeta}}\\
&=& s_{j-1}\gamma x, & \text{\cref{gammavsdegeneracies}(1)}
\end{array}
\]
as desired.

Finally, we prove (1) by showing that $\beta w_jx=s_{j-1}\beta x$. For this, we have that
\[
\begin{array}{ccll}
     s_{j-1}\beta x&=&s_{j-1}s_r^s\gamma x&  \text{\cref{DefGammaBeta}}\\
     &=&s_{r+1}^ss_{j-1}\gamma x& \text{\cref{SimpIdentities}}\\
     &=&s_{r+1}^s\gamma w_j x&\text{Previous step}\\
     &=&\beta w_jx,&\text{\ref{SuspectPartnerRank}}
\end{array}
\]
so done with (1). 

For (2), the statement for $j=0$ and $j=r$ is an immediate consequence of~\cref{vstabilize,gammavsdegeneracies}. For $1\leq j < r$, we have
\[
\begin{array}{cll}
     \beta v_jx&= s_r^s\gamma v_jx&  \text{\cref{steiner5.3}}\\
     &= s_r^s\gamma d_jw_jx& \text{\cref{VvsW}}\\
     &= s_r^sd_j\gamma w_jx&\text{\cref{steiner5.1}}\\
     &= d_js_{r+1}^s\gamma w_jx&\text{\cref{SimpIdentities}}\\
     &= d_j\beta w_jx&  \text{\cref{steiner5.3}}\\
     &= d_js_{j-1}\beta x&\text{Part (1)}\\
     &=\beta x &\text{\cref{SimpIdentities}}.
\end{array}
\]

For (3), use~\cref{steiner5.3} to obtain 
\[
\begin{array}{cll}
v_jv_jx&= s_j^{r-j}d_j^{r-j}(v_jx-\beta v_j x) + \beta v_jx& \text{\cref{closedformulas}}\\
&= s_j^{r-j}d_j^{r-j}(v_jx-\beta x) + \beta x & \text{Part (2)}\\
&= s_j^{r-j}d_j^{r-j}(s_j^{r-j}d_j^{r-j}(x-\beta  x) + \beta x-\beta x) + \beta x & \text{\cref{closedformulas}}\\
&= s_j^{r-j}d_j^{r-j}s_j^{r-j}d_j^{r-j}(x-\beta  x)  + \beta x &\\
&= s_j^{r-j}d_j^{r-j}(x-\beta  x)  + \beta x &\text{\cref{SimpIdentities}}\\
&= v_jx & \text{\cref{closedformulas}},
\end{array}
\]
as desired.
\end{proof}

\begin{proof}[Proof of \ref{SuspectPartnerIndex} and \ref{SuspectPartnerSuspect}]
Let $x\colon O[m]\to D\star O[\top]$ be an augmented directed chain map.
We show that the simplex $v_{p-1}w_px$ is degenerate (at $p-1$) while the simplex 
 $x$ is non-degenerate, so in particular $v_{p-1}w_px\neq x$.
We have that
\[
\begin{array}{ccll}
v_{p-1}w_{p}x&=&s_{p-1}^{(r+1)-(p-1)}d_{p-1}^{(r+1)-(p-1)}(w_px-\beta w_px)+\beta w_px&\text{\cref{closedformulas}}\\
&=&s_{p-1}^{r+2-p}d_{p-1}^{r+2-p}(w_px-s_{p-1}\beta x)+s_{p-1}\beta x&\text{\cref{vjidempotent}}\\
&=&s_{p-1}\left(s_{p-1}^{r+1-p}d_{p-1}^{r+2-p}(w_px-s_{p-1}\beta x)+\beta x\right) &\in \im s_{p-1},\\
\end{array}\]
as desired.
We already argued
that $v_{p-1}w_px$ is degenerate at $p-1$, so it only remains to show that $v_pw_px=w_px$. For this, we have that
\[\begin{array}{ccll}
v_pw_px&=&s_p^{r+1-p}d_p^{r+1-p}(w_px-\beta w_px)+\beta w_px&\text{\cref{closedformulas}}\\
&=&s_p^{r+1-p}d_p^{r+1-p}(w_px-s_{p-1}\beta x)+s_{p-1}\beta x&\text{\cref{vjidempotent}} \\
&=&s_p^{r+1-p}d_p^{r-p}(d_pw_px-\beta x)+s_{p-1}\beta x&\text{\cref{SimpIdentities}}\\
&=&s_p^{r+1-p}d_p^{r-p}(v_px-\beta x)+s_{p-1}\beta x&\text{\cref{VvsW}}\\
&=&s_p^{r+1-p}d_p^{r-p}(x-\beta x)+s_{p-1}\beta x&\\
&=&w_px.&\text{\cref{closedformulas}}
\end{array}
\]
In conclusion, the second chain of equalities shows that the suspect index of $w_px$ is at most $p$, while the first chain of equalities shows at the same time that the suspect index is exactly $p$ and furthermore that $w_px$ is a suspect simplex.
\end{proof}

This completes the proof of \cref{wpwelldef}, showing in particular that $w_p$ assigns a suspect simplex to any non-suspect simplex of $\NRS(\cD\star\cO[\top])$.

\subsubsection{The statements for the dual case}
\label{dualw}
We state what are the features of  $w_px$ for an augmented directed chain map
\[x\colon O[m]\to O[\bot]\star D.\]
The following is the analog of~\cref{wpwelldef}. It is almost identical to the original stamement, with the only exception of the item that describes the suspect index.

\begin{prop}
\label{wpwelldef2}
Let $x\colon O[m]\to O[\bot]\star D$ be a non-degenerate non-suspect simplex with suspect index $p$.
The following hold.
\begin{enumerate}[leftmargin=*, label={\normalfont (W\arabic*)}, ref={\normalfont (W\arabic*)}]
\item \label{SuspectPartnerDimv2} If $x$ has dimension $m$, the simplex $w_px$ has dimension $m+1$.
\item \label{SuspectPartnerDegv2} The simplex $w_px$ is non-degenerate.
\item \label{SuspectPartnerTermv2} If $x$ has non-negative corank, the simplex $w_px$ has non-negative corank.
\item \label{SuspectPartnerConev2}
If $x$ is non-conical, then $w_px$ is also non-conical.
\item \label{SuspectPartnerRankv2} If $x$ has rank $r$, the simplex $w_px$ has rank $r+1$.
\item \label{SuspectPartnerIndexv2} The simplex $w_px$ has suspect index $m+1-p$.
\item \label{SuspectPartnerSuspectv2} The simplex $w_px$ is a suspect simplex.
\end{enumerate}
\end{prop}

The details of the proof can be filled in with similar techniques to those employed in the original case.

For the proof of \ref{SuspectPartnerDegv2} it is useful to observe that applying~\cref{lem:dual_iterated_wedge} to~\cref{VvsW2}, we have 
\[
w_px=v_{p-1}x \wedge_{m-p} \widetilde{\eta}_{p-1}x
\]
for some $\widetilde{\eta}_{p-1}x$, a formula that plays a similar role to that of \eqref{auxiliaryformulaeta}.

The proofs of \ref{SuspectPartnerIndexv2} and \ref{SuspectPartnerSuspectv2} make use of the following lemma, which is the analog of~\cref{vjidempotent}.

\begin{lem}
 \label{vjidempotentv2}
 Let $x\colon O[m]\to O[\bot]\star D$ be an augmented directed chain map of rank $r$ and corank $s=m-r$. Then the following holds:
 \begin{enumerate}
     \item For all $1\leq j \leq r$, we have $\beta w_jx=s_{m-j+1}\beta x$.
     \item For all $0 \leq j \leq r$, we have $\beta v_jx=\beta x$. 
     \item For all $0 \leq j \leq r$, we have $v_j v_jx=v_j x$.  
 \end{enumerate}
 \end{lem}
To prove this lemma, the auxiliary formula replacing \eqref{auxiliaryformula} is \[\gamma w_jx=\gamma s_{m-j+1}\beta x= s_{r-j+1}\gamma x.\]

\subsection{Proof of \texorpdfstring{$w_pd_py=y$}{w_p d_p y = y}}

 \begin{lem} 
 \label{SuspectFaceRank}
 If a suspect simplex $y\colon O[M]\to D\star O[\top]$ has rank $R\ge1$ and suspect index $p$, then $d_py$ has rank $R-1$.
 \end{lem}
 
 \begin{proof}
Since $y$ is a suspect simplex, we know by~\cref{SuspectRankSuspectIndex} that its suspect index $p$ is strictly smaller than $R$. Thus we can apply~\cref{steiner5.1} 
to deduce that the rank of $d_py$ is $R-1$ (and the corank is $s=(M-1)-(R-1)=M-R$).
 \end{proof}

We obtain for free using~\cref{VvsW} that $ d_pw_px=x$, so $d_p$ and $w_p$ define partially inverse functions. The following lemma addresses the other composite.  Note that $p\leq R-1$, and the latter is the rank of $d_py$ by~\cref{SuspectFaceRank}, so the claim makes sense.

\begin{prop}\label{DifficultComposition} 
If $y\colon O[M]\to D\star O[\top]$ is a suspect simplex of rank $R$
and suspect index $p$, then
\[w_p d_py=y.\]
\end{prop}

Before proving the proposition, we give a lemma.

\begin{lem} \label{BetaSuspect}\label{BetaCriticalFace}
Let $y\colon O[M]\to D\star O[\top]$ be a
suspect simplex of rank $R$ and suspect index $p$. Then the following hold.
\begin{enumerate}
    \item $\beta d_py=d_p\beta y$.
    \item $\beta y=s_{p-1} d_p\beta y$.
\end{enumerate}
\end{lem}

\begin{proof}[Proof of~\cref{BetaCriticalFace}]
For (1), we know by~\cref{SuspectRankSuspectIndex} that
$p<R$. Then the statement follows from~\cref{betaandfaces}.

We now show (2).
On the one hand, by~\cref{closedformulas} we have that:
\[v_{p-1}y=s_{p-1}^{R-p+1}d_{p-1}^{R-p+1}(y-\beta y) + \beta y.
\]
On the other hand, by~\cref{closedformulas} and using the fact that $y$ is a suspect simplex of index $p$ and rank $R$, we have that
\[
\begin{array}{rcll}
v_{p-1}y&=&s_{p-1}d_{p}v_{p-1}y& y\text{ suspect}\\
&=&s_{p-1}d_{p}(s_{p-1}^{R-p+1}d_{p-1}^{R-p+1}(y-\beta y) + \beta y)& \text{\cref{closedformulas}}\\
&=&s_{p-1}d_{p}s_{p-1}^{R-p+1}d_{p-1}^{R-p+1}(y-\beta y) + s_{p-1}d_{p}\beta y&\\
&=&s_{p-1}^{R-p+1}d_{p-1}^{R-p+1}(y-\beta y) + s_{p-1} d_p\beta y.&\text{\cref{SimpIdentities}}
\end{array}\]
Since the first summands in the two expressions of $v_{p-1}y$ we obtained are equal, we conclude that the second summands are equal, namely:
\[\beta y=s_{p-1} d_p\beta y,\]
as desired.
\end{proof}

We can now prove the proposition.

\begin{proof}[Proof of~\cref{DifficultComposition}] 
Recall from~\cref{SuspectFaceRank} that the rank of $d_py$ is $R-1$, and that its corank is $s$. We have that
\[
\begin{array}{ccll}
    w_pd_py & = &s_p^{R-1-p+1}d_p^{R-1-p}(d_py-\beta d_py) + s_{p-1}\beta d_py& \text{\cref{closedformulas}}\\
     &=&s_p^{R-p}d_p^{R-1-p}(d_py-d_p\beta y) + s_{p-1}d_p\beta y &\text{\cref{BetaCriticalFace}(1)} \\
     &=&s_p^{R-p}d_p^{R-1-p}(d_py-d_p\beta y) + \beta y &\text{\cref{BetaCriticalFace}(2)} \\
    &=&s_{p}^{R-p}d_p^{R-p}(y-\beta y)+\beta y&\\
     &=&v_py=y,& \text{\cref{closedformulas}}
\end{array}\]
as desired.
\end{proof}

  \begin{lem} 
 \label{SuspectFaceIndex}
 If a suspect simplex $y\colon O[M]\to D\star O[\top]$ has rank $R\ge1$ and suspect index $p$, then $d_py$ has suspect index $p$.
 \end{lem}
 
 \begin{proof}
To compute the suspect index of $d_py$, we first have to show that $v_pd_py=d_py$. For this, we observe that
\[
\begin{array}{lllll}
v_pd_py&=&d_pw_pd_py&\text{\cref{VvsW}}\\
&=&d_py,&\text{\cref{DifficultComposition}}
\end{array}
\]
as desired. It remains to show that $v_{p-1}d_py\neq d_py$.

As an auxiliary fact, we show that $d_pv_{p-1}y=v_{p-1}d_py$.
For this, we observe that
\[
\begin{array}{ccll}
s_{p-1}v_{p-1}d_py&=& s_{p-1}(s_{p-1}^{R-1-p+1}d_{p-1}^{R-1-p+1}(d_py-\beta d_py)+\beta d_py)&\text{\cref{closedformulas}}\\ 
&=& s_{p-1}^{R-p+1}d_{p-1}^{R-1-p+1}(d_py-\beta d_py)+s_{p-1}\beta d_py&\\ 
&=& s_{p-1}^{R-p+1}d_{p-1}^{R-p}(d_py-d_p\beta y)+s_{p-1}d_p\beta y&\text{\cref{BetaCriticalFace}}\\ 
&=& s_{p-1}^{R-p+1}d_{p-1}^{R-p}(d_py-d_p\beta y)+\beta y& \text{\cref{BetaSuspect}}\\
&=& s_{p-1}^{R-p+1}d_{p-1}^{R-p+1}(y-\beta y)+\beta y& \text{\cref{IteratedSimpIdentities}}\\ 
&=& v_{p-1}y, & \text{\cref{closedformulas}}
\end{array}
\]
which implies
\[d_pv_{p-1}y=d_ps_{p-1}v_{p-1}d_py=v_{p-1}d_py,\] as desired.

Next, in order to show that $v_{p-1}d_py\neq d_py$, we assume that $v_{p-1}d_py= d_py$ and obtain a contradiction. 
By what we have just shown, we obtain that \[\begin{array}{lllll}
    v_{p-1}y &= & s_{p-1} v_{p-1}d_py &\text{Previous formula}\\
     &  = & s_{p-1}d_py.&\text{Assumption}
\end{array}\]

Next, we obtain that
\[
\begin{array}{ccll}
   y & = &v_py& \text{Assumption}\\
     &=&d_pw_py &\text{\cref{VvsW}(1)} \\
     &=&d_p(\alpha_{p-1}y\wedge^{R-p+1}_{p-1} v_{p-1}y) &\text{\cref{VvsW}(3)} \\
    &=&d_p(s_{p}^{R-p+1}\alpha_{p-1}y - s_{p-1}^{R-p+2}d_{p}^{R-p+1}v_{p-1}y+s_{p-1}v_{p-1}y)&\text{\cref{DefIteratedWedge}}\\
     &=&s_{p}^{R-p}\alpha_{p-1}y - s_{p-1}^{R-p+1}d_{p}^{R-p+1}v_{p-1}y+v_{p-1}y& \text{\cref{SimpIdentities}}\\
      &=&s_{p}^{R-p}\alpha_{p-1}y - s_{p-1}^{R-p+1}d_{p}^{R-p+1}s_{p-1}d_py+s_{p-1}d_py& \text{Formula above}\\
      &=&s_{p}^{R-p}\alpha_{p-1}y - s_{p-1}^{R-p+1}d_{p}^{R-p+1}y+s_{p-1}d_py,& \text{\cref{SimpIdentities}}\\
\end{array}
\]
so that we have the relation
\begin{equation}
\label{longformula}
y = s_{p}^{R-p}\alpha_{p-1}y - s_{p-1}^{R-p+1}d_{p}^{R-p+1}y+s_{p-1}d_py.
\end{equation}
By applying $d_p$ 
we obtain that
\[
d_py =
s_{p}^{R-p-1}\alpha_{p-1}y - s_{p-1}^{R-p}d_{p}^{R-p+1}y+d_py.
\]
By cancellation of $d_py$ in this last formula, we get that
\[
s_{p}^{R-p-1}\alpha_{p-1}y =s_{p-1}^{R-p}d_{p}^{R-p+1}y.
\]
Substituting this back into (\ref{longformula}) we obtain that
\[
\begin{array}{ccll}
y&=&s_{p}^{R-p}\alpha_{p-1}y - s_{p-1}^{R-p+1}d_{p}^{R-p+1}y+s_{p-1}d_py& \\
&=& s_ps_{p-1}^{R-p}d_{p}^{R-p+1}y - s_{p-1}^{R-p+1}d_{p}^{R-p+1}y+s_{p-1}d_py&\\
&=& s_{p-1}^{R-p+1}d_{p}^{R-p+1}y - s_{p-1}^{R-p+1}d_{p}^{R-p+1}y+s_{p-1}d_py&\text{\cref{IteratedSimpIdentities}}\\
&=& s_{p-1}(s_{p-1}^{R-p}d_{p}^{R-p+1}y - s_{p-1}^{R-p}d_{p}^{R-p+1}y+d_py)&\in \im s_{p-1},\\
\end{array}
\]
which contradicts the assumption that $y$ is non-degenerate. 
\end{proof}

\subsubsection{The statements for the dual case}

 The following lemma is the analog of \cref{SuspectFaceRank}
 
 \begin{lem} 
 \label{SuspectFaceRankv2}
 If a suspect simplex $y\colon O[M]\to O[\bot]\star D$ has rank $R\ge1$ and suspect index $p$, then $d_{M-p}y$ has rank $R-1$.
 \end{lem}

We obtain for free using~\cref{VvsW2} that $ d_{m+1-p}w_px=x$, so $d_{m+1-p}$ and $w_p$ define partially inverse functions. The following proposition addresses the other composite.  Note that $p\leq R-1$, and the latter is the rank of $d_{M-p}y$ by~\cref{SuspectFaceRankv2}, so the claim makes sense.

The following is the analog of \cref{DifficultComposition}.

\begin{prop}\label{DifficultCompositionv2} 
If $y\colon O[M]\to O[\bot]\star D$ is a suspect simplex of rank $R$
and suspect index $p$, then
\[w_p d_{M-p}y=y.\]
\end{prop}

The proposition relies on the following lemma, which is the analog of \cref{BetaSuspect}.

\begin{lem} \label{BetaSuspectv2}\label{BetaCriticalFacev2}
Let $y\colon O[M]\to O[\bot]\star D$ be a
suspect simplex of rank $R$ and suspect index $p$. Then the following hold.
\begin{enumerate}
    \item $\beta d_{M-p}y=d_{M-p}\beta y$.
    \item $\beta y=s_{M-p} d_{M-p}\beta y$.
\end{enumerate}
\end{lem}

The following is the analog of \cref{SuspectFaceIndex}.

\begin{lem} 
 \label{SuspectFaceIndexv2}
 If a suspect simplex $y\colon O[M]\to O[\bot]\star D$ has suspect index $p$, then $d_{M-p}y$ has suspect index $p$.
 \end{lem}

\subsection{Features of \texorpdfstring{$d_p y$}{d_p y}}

The goal of this subsection is to collect and prove all the relevant properties and features of $d_py$ we need, when $y$ is a suspect simplex. We include the relevant statements for the dual case in~\cref{dpwelldefv2}.

\begin{prop} 
\label{dpwelldef}
For any suspect simplex $y\colon O[M]\to D\star O[\top]$ with suspect index $p$, the following hold.
\begin{enumerate}[leftmargin=*, label={\normalfont (D\arabic*)}, ref={\normalfont (D\arabic*)}]
\item\label{SuspectFaceDim} If $y$ has dimension $M$, the simplex $d_py$ has dimension $M-1$.
\item\label{SuspectFaceDeg}
$ d_py$ is non-degenerate.
\item\label{SuspectFaceTerm}
If $y$ has non-negative corank, the simplex $d_py$ has non-negative corank.
\item \label{SuspectFaceCone} $d_py$ is non-conical.
\item \label{SuspectFaceRankRepeat} If $y$ has rank $R$, the simplex $d_py$ has rank $R-1$.
\item\label{SuspectFaceIndexRepeat} $d_py$ has suspect index $p$.
\item \label{SuspectFaceNonSuspect}$d_py$ is a non-suspect simplex.
\end{enumerate}
\end{prop}

\begin{proof}[Proof of \ref{SuspectFaceDim}]
This is by definition of the face map $d_p$.
\end{proof}

\begin{proof}[Proof of \ref{SuspectFaceRankRepeat}]
This was done in~\cref{SuspectFaceRank}.
\end{proof}

\begin{proof}[Proof of \ref{SuspectFaceDeg}]
Let $R=r+1$ be the rank of $y$, so that by~\cref{SuspectFaceRank} the rank of $d_py$ is $r$.
Notice that by~\cref{SuspectRankSuspectIndex} we have $p < R$.
We assume that $d_py=s_jz$ for some
$(M-2)$-dimensional simplex $z$ and some $j\leq M-2$ and deduce that $y$ is degenerate, contradicting the assumption.

We first treat the case $j<r$. For this, we observe that
\[\begin{array}{llllll}
    y & =&w_pd_py = w_ps_jz \\
     & =&s_p^{r-p+1} d_{p}^{r-p}(s_jz-\beta s_jz)+s_{p-1}\beta s_jz&\text{\cref{closedformulas}}\\
       & =&s_p^{r-p+1} d_{p}^{r-p}s_j(z-\beta z)+s_{p-1}s_j\beta z.&\text{\cref{gammavsdegeneracies}}
\end{array}\]
If furthermore $j< p-1$, we obtain that
\[\begin{array}{llllll}
y&=&s_p^{r-p+1} d_{p}^{r-p}s_j(z-\beta z)+s_{p-1}s_j\beta z\\
&=&s_p^{r-p+1} s_jd_{p-1}^{r-p}(z-\beta z)+s_{p-1}s_j\beta z&\text{\cref{SimpIdentities}}\\
&=&s_js_{p-1}^{r-p+1} d_{p-1}^{r-p}(z-\beta z)+s_js_{p-2}\beta z&\text{\cref{SimpIdentities}}\\
&=&s_j(s_{p-1}^{r-p+1} d_{p-1}^{r-p}(z-\beta z)+s_{p-2}\beta z)&\in\im s_j,
\end{array}\]
implying that $y$ is degenerate.
If instead $j=p-1$, then we obtain that 
\[\begin{array}{llllll}
y&=&s_p^{r-p+1} d_{p}^{r-p}s_{p-1}(z-\beta z)+s_{p-1}s_{p-1}\beta z\\
&=&s_p^{r-p+1} d_{p}^{r-p}s_{p-1}(z-\beta z)+s_ps_{p-1}\beta z&\text{\cref{SimpIdentities}}\\
&=&s_p(s_p^{r-p} d_{p}^{r-p}s_{p-1}(z-\beta z)+s_{p-1}\beta z)&\in\im s_{p},
\end{array}
\]
implying that $y$ is degenerate.
If instead $j>p-1$, we obtain that
\[\begin{array}{llllll}
y&=&s_p^{r-p+1} d_{p}^{r-p}s_j(z-\beta z)+s_{p-1}s_j\beta z\\
&=&s_p^{r-p+1} d_{p}^{r-p}(s_jz-s_j \beta z)+s_{j+1}s_{p-1}\beta z&\text{\cref{SimpIdentities}}\\
&=&s_{j+1}s_p^{r-p} d_{p}^{r-p}(s_jz-s_j\beta  z)+s_{j+1}s_{p-1}\beta z&\text{\cref{IteratedSimpIdentities}}\\
&=&s_{j+1}(s_p^{r-p} d_{p}^{r-p}(s_jz-s_j\beta z)+s_{p-1}\beta z)&\im s_{j+1},
\end{array}\]
implying that $y$ is degenerate.

Next, we treat the case $j\geq r$. Notice that we have $r \leq j \leq M-2 < M-1 = r +s$,
where $s$ is the corank of $d_py$, so that necessarily $s>0$.
Observe that
\[\begin{array}{llllll}
    y   & =& w_pd_py = w_ps_jz&\text{\cref{DifficultComposition} } \\
        & =& s_p^{r-p+1} d_{p}^{r-p}(s_jz-\beta s_jz)+s_{p-1}\beta s_jz&\text{\cref{closedformulas}}\\
        & =& s_p^{r-p+1} d_{p}^{r-p}(s_jz-s_r^{s}\gamma s_jz)+s_{p-1}s_r^{s}\gamma s_jz&\text{\cref{DefGammaBeta}}\\
        & =& s_p^{r-p+1} d_{p}^{r-p}(s_jz-s_js_r^{s-1}\gamma s_jz)+s_{p-1}s_js_r^{s-1}\gamma s_jz&\text{(*)}\\
        & =& s_p^{r-p+1} d_{p}^{r-p}s_j(z - s_r^{s-1}\gamma s_jz)+s_{p-1}s_js_r^{s-1}\gamma s_jz& \\
        & =& s_p^{r-p+1}s_{j-r+p} d_{p}^{r-p}(z-s_r^{s-1}\gamma s_jz)+s_{j+1}s_{p-1}s_r^{s-1}\gamma s_j z&\text{\cref{IteratedSimpIdentities}}\\
        & =& s_{j+1} s_p^{r-p+1}d_{p}^{r-p}(z-s_r^{s-1}\gamma s_jz)+s_{j+1}s_{p-1}s_r^{s-1}\gamma s_j z&\text{\cref{IteratedSimpIdentities}}\\
        & =& s_{j+1}( d_{p}^{r-p}(z-s_r^{s-1}\gamma s_jz)+s_{p-1}s_r^{s-1}\gamma s_j z)&\in\im s_{j+1},\\
\end{array}\]
where the equality (*) uses~\cref{IteratedSimpIdentities} if $s>1$ and
otherwise follows from the fact that $s=1$ implies $j=r$.
This entails that $y$ is degenerate, finishing the proof.
\end{proof}

\begin{proof}[Proof of \ref{SuspectFaceTerm}]
We know from~\cref{SuspectRankSuspectIndex} that the rank of $y$ must be strictly higher than $p$. This implies that, if $y$ is of dimension $M$, then $d_py[M-1]=y[M]$. In particular, $y$ is valued in $D$ if and only if 
$d_py$ does.
\end{proof}

\begin{proof}[Proof of \ref{SuspectFaceCone}]
Assuming that $d_py$ is a conical simplex, we obtain that
\[\begin{array}{ccll}
    y&=&w_pd_py&\text{\cref{DifficultComposition}}\\
    &=&w_p\beta d_py &\text{\cref{conecharacterization}}\\ 
    &=&  s_{p}^{R-1-p+1}d_{p}^{R-1-p}(\beta d_py-\beta^2d_py)+s_{p-1}\beta^2 d_py&\text{\cref{closedformulas}}\\
    &=&s_{p-1}\beta d_py,& \text{\cref{gammavsdegeneracies}}
\end{array}\]
contradicting the assumption that $y$ is a non-degenerate simplex.
\end{proof}

\begin{proof}[Proof of \ref{SuspectFaceIndexRepeat}]
This was proven as ~\cref{SuspectFaceIndex}.
\end{proof}

\begin{proof}[Proof of \ref{SuspectFaceNonSuspect}]
We have to show that $v_{p-1}d_py$ is not degenerate at $p-1$. For this, we assume that $v_{p-1}d_py=s_{p-1}z$ with $z=d_{p-1}v_{p-1}d_py$ and obtain a contradiction.

As a preliminary observation, we have that
\[
\begin{array}{ccll}
    s_{p-1}z&=&v_{p-1}d_py \\
    &=& s_{p-1}^{R-1-p+1}d_{p-1}^{R-1-p+1}(d_py-\beta d_py) + \beta d_py&\text{\cref{closedformulas}}\\
    &=& s_{p-1}^{R-p}d_{p-1}^{R-p}(d_py-\beta d_py) + d_p\beta y.&\text{\cref{BetaCriticalFace}}
\end{array}
\]
Solving for $d_p\beta y$, the last formula implies that
\begin{equation}
\label{previouscomputation}
d_p\beta y=s_{p-1}z',\end{equation}
where $z'=z-s_{p-1}^{R-p-1}d_{p-1}^{R-p}(d_py-\beta d_py)$.
We can then express $\beta y$ as follows:
\[
\begin{array}{llll}
\beta y&=&s_{p-1}d_p\beta y&\text{\cref{BetaSuspect}}\\
&=&s_{p-1}s_{p-1}z'&\text{(\ref{previouscomputation})}\\
&=&s_ps_{p-1}z'.
\end{array}
\]
Finally, we can express $y$ as follows:
\[\begin{array}{llll}
y&=&v_py\\
&=&s_p^{R-p}d_p^{R-p}(y-\beta y)+ \beta y\\
&=&s_p^{R-p}d_p^{R-p}(y-s_ps_{p-1}z')+ s_ps_{p-1}z'&\text{\cref{closedformulas}}\\
&=&s_p(s_p^{R-p-1}d_p^{R-p}(y-s_ps_{p-1}z')+s_{p-1}z')&\in\im s_p,
\end{array}\]
which would imply that $y$ is a degenerate simplex, contradicting the fact that $y$ is a suspect simplex and in particular non-degenerate.
\end{proof}

\begin{lem}
\label{BoundarySuspect}
Let $y=w_px\colon O[M]\to D\star O[\top]$ be a suspect $M$-simplex
of rank $R=r+1$
and suspect index $p$. The $i$-th face of $y$, unless it is degenerate or conical, has the property that
\[
d_iy\text{ is }\left\{
\begin{array}{llll}
\text{suspect of dimension $M-1$}     & \text{ if } i<p-1\\
\text{of suspect index at most $p-1$}       & \text{ if } i=p-1\\
\text{suspect of dimension $M-1$}       & \text{ if } i=p+1\\
\text{suspect of dimension $M-1$}        &\text{ if }  r\ge i>p+1\\
\text{of rank $r+1$}            &  \text{ if }i\ge r+1.\\
\end{array}
\right.\]
\end{lem}

\begin{proof}
We treat all cases by distinguishing in terms of the value of $i$.
\begin{itemize}[leftmargin=*]
\item We suppose $i<p-1$.
Then using that $\beta d_iy=d_i\beta y$ by~\cref{betaandfaces}
for $i\leq r$
we obtain the following auxiliary formula:
\[\begin{array}{llll}
  d_iy&=& d_iw_px &  \\
  &=&d_is_p^{r-p+1}d_p^{r-p}(x-\beta x)+d_is_{p-1}\beta x&\text{\cref{closedformulas}}\\
    &= &
    s_{p-1}^{r-p-1}d_{p-1}^{r-p}(d_ix-d_i\beta x)+s_{p-2}d_i\beta x & \text{\cref{SimpIdentities}}\\
    &= &
    s_{p-1}^{r-p-1}d_{p-1}^{r-p}(d_ix-\beta d_ix)+s_{p-2}\beta d_ix & \text{\cref{betaandfaces}}\\
    &=&w_{p-1}d_ix &\text{\cref{closedformulas}}.
    \end{array}\]
    The chain of equalities
\[\begin{array}{llll}
  v_{p-1}d_iy &=& s_{p-1}^{r-p+1}d_{p-1}^{r-p+1}(d_iy-\beta d_iy)+\beta d_iy&\text{\cref{closedformulas}}\\
    &= &
    d_is_{p}^{r-p+1}d_{p}^{r-p+1}(y-\beta y)+d_i\beta y & \text{\cref{SimpIdentities}}\\
    &=&d_i v_py &\text{\cref{closedformulas}}\\
    &=&d_iy\\
    \end{array}\]
    shows that the suspect index of $d_iy$ is at most $p-1$.
    
The chain of equalities
    \[\begin{array}{llll}
  v_{p-2}d_iy &=& s_{p-2}^{r-p+2}d_{p-2}^{r-p+2}(d_iy-\beta d_iy)+\beta d_iy&\text{\cref{closedformulas}}\\
    &= & s_{p-2}^{r-p+2}d_{p-2}^{r-p+2}(d_iy-\beta d_iy)+\beta w_{p-1}d_ix & \text{Formula above}\\
    &=& s_{p-2}^{r-p+2}d_{p-2}^{r-p+2}(d_iy-\beta d_iy)+ s_{p-2} \beta d_ix &\text{\cref{vjidempotent}}
    \end{array}\]
shows that $v_{p-2}d_iy$ is degenerate at $p-2$.

If $v_{p-2}d_iy$ equals $d_iy$, then this face is degenerate, and otherwise, we have shown that $d_iy$ is a suspect simplex of suspect index $p-1$ and of dimension $m$, as desired. 

\item We suppose $i=p-1$. Then we obtain that
\begin{equation}
\label{(p-1)Facey}
    \begin{array}{llll}
   d_{p-1}y  &= &d_{p-1}w_px& \\
&=     & 
d_{p-1}(\alpha_{p-1}x\wedge_{p-1}^{r-p+1}v_{p-1}x)&\text{\cref{VvsW}}\\
&=& v_{p-1}x.
\end{array}
\end{equation}
Next, we observe that
\[\begin{array}{lllll}
  v_{p-1}d_{p-1}y   & = & v_{p-1}v_{p-1}x&\text{Formula \eqref{(p-1)Facey}}\\
  & =  & v_{p-1}x&\text{\cref{vjidempotent}}, 
\end{array}\]
proving that the suspect index of $d_{p-1}y$ is at most $p-1$.

\item We suppose $i=p+1$. Then we obtain the following auxiliary formula using~\cref{steiner5.1}.
\[\begin{array}{llll}
 d_{p+1}y    & = &d_{p+1}w_px &\\
&=     & d_{p+1}s_p^{r-p+1}d_p^{r-p}(x-\beta x)+d_{p+1}s_{p-1}\beta x  & \text{\cref{closedformulas}}\\
&=     & s_p^{r-p}d_p^{r-p}(x-\beta x)+s_{p-1}d_p\beta x &\text{\cref{SimpIdentities}}\\
&=&s_p^{(r-1)-p+1}d_p^{(r-1)-p}(d_px-\beta d_px)+s_{p-1}\beta d_px &\text{\cref{betaandfaces}}\\
&=&w_pd_px & \text{\cref{closedformulas}},
\end{array}\]
Now we give an upper bound for the suspect index of $d_{p+1}y$. 
\[\begin{array}{llll}
  v_{p}d_{p+1}y &=& s_{p}^{r-p}d_{p}^{r-p}(d_{p+1}y-\beta d_{p+1}y)+\beta d_{p+1}y&\text{\cref{closedformulas}}\\
    &= &
    s_{p}^{r-p}d_{p}^{r-p+1}(y-\beta y)+\beta d_{p+1}y & \text{\cref{IteratedSimpIdentities}}\\
    &= &
    d_{p+1}s_{p}^{r-p+1}d_{p}^{r-p+1}(y-\beta y)+\beta d_{p+1}y & \text{\cref{SimpIdentities}}\\
   &= &  d_{p+1}s_{p}^{r-p+1}d_{p}^{r-p+1}(y-\beta y)+d_{p+1}\beta y & \text{\cref{betaandfaces}}\\
    &=&d_{p+1} v_py &\text{\cref{closedformulas}}\\
    &=&d_{p+1}y,\\
    \end{array}\]
        showing that the suspect index of $d_{p+1}y$ is at most $p$. Next, we compute
    \[\begin{array}{llll}
  v_{p-1}d_{p+1}y &=& s_{p-1}^{r-p+1}d_{p-1}^{r-p+1}(d_{p+1}y-\beta d_{p+1}y)+\beta d_{p+1}y&\text{\cref{closedformulas}}\\
    &= & s_{p-1}^{r-p+1}d_{p-1}^{r-p+1}(d_{p+1}y-\beta d_{p+1}y)+\beta w_{p}d_px & \text{Formula above}\\
    &=& s_{p-1}^{r-p+1}d_{p-1}^{r-p+1}(d_{p+1}y-\beta d_{p+1}y)+ s_{p-1} \beta d_px &\text{\cref{vjidempotent}}
    \end{array}\]
        This shows that $v_{p-1}d_{p+1}y$ is degenerate at $p-1$. If $v_{p-1}d_{p+1}y$ equals $d_{p+1}y$, then this face is degenerate, and otherwise, we have shown that $d_{p+1}y$ is a suspect simplex of suspect index $p$ and of dimension $m$, as desired.

\item We suppose $r\ge i>p+1$. Then we obtain the following auxiliary formula using~\cref{steiner5.1}.
\[\begin{array}{llll}
  d_iy&=& d_iw_px   \\
   &= &d_is_p^{r-p+1}d_p^{r-p}(x-\beta x)+d_is_{p-1}\beta x & \text{\cref{closedformulas}}\\
    &= &d_is_p^{r-p+1}d_p^{r-p}(x-\beta x)+s_{p-1}d_{i-1}\beta x & \text{\cref{SimpIdentities}}\\
 &=&
s_p^{r-p}d_p^{r-p}(x-\beta x)+s_{p-1}d_{i-1}\beta x&\text{\cref{IteratedSimpIdentities}}\\
&=&s_p^{r-p}d_p^{r-p-1}(d_{i-1}x-d_{i-1}\beta x)+s_{p-1}d_{i-1}\beta x &\text{\cref{IteratedSimpIdentities}}\\
&=&s_p^{r-p}d_p^{r-p-1}(d_{i-1}x-\beta d_{i-1}x)+s_{p-1}\beta d_{i-1}x&\text{\cref{betaandfaces}}\\
&=&w_pd_{i-1}x&\text{\cref{closedformulas}},
\end{array}\]
Now we give an upper bound for the suspect index of $d_{i}y$.
\[\begin{array}{llll}
  v_{p}d_{i}y &=& s_{p}^{r-p}d_{p}^{r-p}(d_{i}y-\beta d_{i}y)+\beta d_{i}y&\text{\cref{closedformulas}}\\
    &= &
    d_{i}s_{p}^{r-p+1}d_{p}^{r-p+1}(y-\beta y)+\beta d_{i}y & \text{\cref{SimpIdentities}}\\
   &= &  d_{i}s_{p}^{r-p+1}d_{p}^{r-p+1}(y-\beta y)+d_{i}\beta y & \text{\cref{betaandfaces}}\\
    &=&d_{i} v_py &\text{\cref{closedformulas}}\\
    &=&d_{i}y,\\
    \end{array}\]
        showing that the suspect index of $d_{p+1}y$ is at most $p$. Next, we compute
    \[\begin{array}{llll}
  v_{p-1}d_{i}y &=& s_{p-1}^{r-p+1}d_{p-1}^{r-p+1}(d_{i}y-\beta d_{i}y)+\beta d_{i}y&\text{\cref{closedformulas}}\\
    &= & s_{p-1}^{r-p+1}d_{p-1}^{r-p+1}(d_{i}y-\beta d_{i}y)+\beta w_{p}d_{i-1}x & \text{Formula above}\\
    &=& s_{p-1}^{r-p+1}d_{p-1}^{r-p+1}(d_{i}y-\beta d_{i}y)+ s_{p-1} \beta d_{i-1}x &\text{\cref{vjidempotent}}
    \end{array}\]
        This shows that $v_{p-1}d_{i}y$ is degenerate at $p-1$. If $v_{p-1}d_{i}y$ equals $d_{i}y$, then this face is degenerate, and otherwise, we have shown that $d_{i}y$ is a suspect simplex of suspect index $p$ and of dimension $m$, as desired. 

\item We suppose $i>r$. Then the rank of $d_iy$ equals the rank of $y$ (see~\cref{rmk:face_rank}), namely $r+1$,
as desired.
\qedhere
\end{itemize}
\end{proof}

\subsubsection{The statements for the dual case}

The goal of this subsection is to collect the statements describing the features of $d_{M-p}y$ we need, when $y$ is a suspect simplex of the form $y\colon O[M]\to O[\bot]\star D$.

The following is the analog of~\cref{dpwelldef}.

\begin{prop}
\label{dpwelldefv2}
For any suspect simplex $y\colon O[M]\to O[\bot]\star D$ with suspect index $p$, the following hold.
\begin{enumerate}[leftmargin=*, label={\normalfont (D\arabic*)}, ref={\normalfont (D\arabic*)}]
\item\label{SuspectFaceDimv2} If $y$ has dimension $M$, the simplex $d_{M-p}y$ has dimension $M-1$.
\item\label{SuspectFaceDegv2}
$d_{M-p}y$ is non-degenerate.
\item\label{SuspectFaceTermv2} 
If $y$ has non-negative corank, the simplex $d_{M-p}y$ has non-negative corank.
\item \label{SuspectFaceConev2} $d_{M-p}y$ is non-conical.
\item \label{SuspectFaceRankRepeatv2} If $y$ has rank $R$, the simplex $d_{M-p}y$ has rank $R-1$.
\item\label{SuspectFaceIndexv2} $d_{M-p}y$ has suspect index $p$.
\item \label{SuspectFaceNonSuspectv2}$d_{M-p}y$ is a non-suspect simplex.
\end{enumerate}
\end{prop}

The following is the analog of~\cref{BoundarySuspect}.

\begin{lem}
\label{BoundarySuspectv2}
Let $y=w_px\colon O[M]\to O[\bot]\star D$ be a suspect $M$-simplex
of rank $R=r+1$
and suspect index $p$. The $i$-th face of $y$, unless it is degenerate or conical, has the property that
\[
d_iy\text{ is }\left\{
\begin{array}{llll}
\text{of suspect index at most $p-1$}     & \text{ if } m-p+2<i\leq m+1\\
\text{of suspect index at most $p-1$}       & \text{ if } i=m-p+2\\
\text{suspect of dimension $M-1$}        &\text{ if } m-r<i\leq m-p\\
\text{of rank $r+1$}            &  \text{ if }i\le m-r.\\
\end{array}
\right.\]
\end{lem}

\section{Proof of the theorem}
\label{proofmaintheorem}

We are finally ready to prove the main result(s):~\cref{maintheorem,maintheoremdual}.
We treat the former in detail, and briefly indicate in~\cref{maindual} how to treat the latter.

After having understood $N_m(\cD\star\cO[\top])$ (resp.~$N_m(\cO[\bot\star\cD])$) as the set of augmented directed chain maps $x\colon O[m]\to D\star O[\top]$ (resp.~$x\colon O[m]\to O[\bot]\star \cD$), and having defined what it means for one such simplex to be a suspect or not, we can show that there is a pairing between suspect simplices and non-suspect simplices, as follows.

\begin{lem}
\label{lemmabijection}
Let $\cD$ be a strong Steiner category, and recall the inclusion $\NRS\cD\star\Delta[\top]\hookrightarrow\NRS(\cD\star\cO[\top])$ from~\cref{comparisonmap}. Let $d\ge1$.
\begin{enumerate}[leftmargin=*]
    \item The non-degenerate $d$-simplices in $\NRS\cD\star\Delta[\top]$, regarded as $d$-simplices of $\NRS(\cD\star\cO[\top])$, are precisely the conical simplices and the simplices contained in $\NRS\cD$.
    \item The $p$-th face defines a bijective correspondence between the (necessarily non-degenerate and non-conical) suspect $(d+1)$-simplices of suspect index $p$ and the non-degenerate non-conical non-suspect $d$-simplices of suspect index $p$. 
\end{enumerate}
\end{lem}

\begin{proof}
We leave the verification of (1) to the reader 
and we now show (2). For this, the collection of results from the previous section allows us to conclude that the $p$-th face and the construction $w_p$ define inverse bijections
\[w_p\colon\left\{
x\ \begin{array}{cc}
     \text{non-degenerate}&  \\
    \text{non-conical} & \\
    \text{non-suspect}
\end{array}
\begin{array}{c}
\text{with}\\
\text{suspect}\\
\text{index $p$}
\end{array} \right\}\cong\left\{
y\ \text{suspect}
\begin{array}{c}
\text{with}\\
\text{suspect}\ \\
\text{index $p$}
\end{array}\right\}\colon d_p.\]
Indeed, we showed that $ d_p$ is well-defined in~\cref{dpwelldef}, that $w_p$ is well-defined in~\cref{wpwelldef}, that the composite $d_p  w_p$ is an identity in~\cref{VvsW}, and finally that the composite $w_p  d_p$ is an identity in~\cref{DifficultComposition}. 
\end{proof}

We are now able to prove the first case of main theorem of the paper.

\begin{thm}
\label{maintheorem}
Let $\cD$ be a
strong Steiner $\omega$-category. The inclusion \[\NRS\cD\star\Delta[\top]\hookrightarrow\NRS(\cD\star\cO[\top])\]
from~\cref{comparisonmap} is a complicial inner anodyne extension, and in particular an acyclic cofibration in the model structure(s) from~\cref{modelstructurewithsaturation}.
\end{thm}

The strategy of the proof consists in adding all simplices of $\NRS(\cD\star\cO[\top])$
missing from $\NRS\cD\star\Delta[\top]$ inductively on their (ascending) dimension $d$, their
(descending) rank $r$, and their (ascending) suspect index $p$.
This provides us with encapsulated chains of inner anodyne extensions.
For each choice of dimension $d$, rank $r$ and suspect index $p$, the inner anodyne
extension in question is constructed by filling $p$-horns where the principal simplices
are $(d+1)$-dimensional suspects of index $p$ and rank $r$ and so the corresponding missing $p$-faces
are given by the bijection established in the previous lemma.

\begin{proof}
 In order to show that the inclusion $\NRS\cD\star\Delta[\top]\hookrightarrow \NRS(\cD\star\cO[\top])$ is a complicial inner anodyne extension, we will realize it as a transfinite composite of intermediate complicial inner anodyne extensions
\[\NRS\cD\star\Delta[\top]=:X_0\hookrightarrow X_1\hookrightarrow X_2\hookrightarrow\dots\hookrightarrow X_{d-1}\hookrightarrow X_{d}\hookrightarrow\dots \hookrightarrow  \NRS(\cD\star\cO[\top]).\]
For $d\geq 1$, we let $X_d$ be the smallest regular subsimplicial set of $\NRS(\cD\star\cO[\top])$ 
containing $X_{d-1}$, all $d$-simplices of $\NRS(\cD\star\cO[\top])$, 
as well as
the suspect $(d+1)$-simplices of $\NRS(\cD\star\cO[\top])$.
Note that $X_0$ already contains all $0$-simplices of $\NRS(\cD\star\cO[\top])$ and that there are no non-degenerate suspect $1$-simplices by~\cref{SuspectRankSuspectIndex}.
We see that the difference between $X_{d-1}$ and $X_d$ are the non-degenerate non-suspect $d$-simplices of rank at most $d$ and the non-degenerate suspect $(d+1)$-simplices of rank at most $d+1$.

In order to show that the inclusion $X_{d-1}\hookrightarrow X_{d}$ is a complicial inner anodyne extension for all $d\geq 1$, we realize it as a composite of intermediate complicial inner anodyne extensions
\[X_{d-1}=:Y_{d+1} \hookrightarrow Y_{d}\hookrightarrow Y_{d-1} \hookrightarrow \ldots \hookrightarrow Y_{r+1}\hookrightarrow Y_{r}\hookrightarrow\ldots \hookrightarrow Y_{1}=X_{d}.\]
For $1\le r\leq d$,
let $Y_r$ be the
smallest regular subset of $X_{d-1}$ containing $Y_{r+1}$ as well as 
all (necessarily non-degenerate)
suspect $(d+1)$-simplices $y$ of $\NRS(\cD\star\cO[\top])$ of rank $r+1$ and all non-degenerate non-suspect $d$-simplices of rank $r$. 

We see using~\cref{BoundarySuspect,lemmabijection}
that the difference between $Y_{r}$ and $Y_{r+1}$ are
the non-degenerate suspect $(d+1)$-simplices of rank $r$ and possibly some of their faces (precisely those that are neither suspect nor of rank $r$ or higher).

In order to show that the inclusion $Y_{r+1}\hookrightarrow Y_r$ is a complicial inner anodyne extension for $1\le r\le d$, we realize it as a filtration made by intermediate complicial inner anodyne extensions
\[Y_{r+1}=:W_0\hookrightarrow W_1 \hookrightarrow \ldots \hookrightarrow W_{p-1}\hookrightarrow W_p\hookrightarrow\ldots\hookrightarrow W_{r}=Y_{r}.\]
For $0<p\leq r$, we let $W_p$ be the smallest regular simplicial subset of $Y_{r+1}$
containing $W_{p-1}$ as well as 
all suspect $(d+1)$-simplices of $ \NRS(\cD\star\cO[\top])$ of rank $r$ and suspect index $p$.
Note that any $d$-simplex of suspect index $0$ is either degenerate or conical and thus are already in $X_{d-1}\subseteq W_0$.
We see using~\cref{BoundarySuspect,lemmabijection} and~\cref{SuspectFaceIndex} that the difference between $W_{p-1}$ and $W_p$ are the non-degenerate suspect $(d+1)$-simplices $y$ of rank $r$ and suspect index $p$ and the non-degenerate non-suspect $d$-simplices $x$
of rank $r-1$ and suspect index $p$.
There is a bijective correspondence between the $(d+1)$- and $d$-simplices mentioned above, as shown in~\cref{lemmabijection}.
 
We now record some relevant properties of the $(d+1)$-simplices $y$ as above. 
\begin{itemize}[leftmargin=*]
    \item We use~\cref{BoundarySuspect} to argue that the $p$-horn of $y$ belongs to $W_{p-1}$; in particular, the $p$-horn defines a map of (underlying) simplicial sets 
  \[\Lambda^{p}[d+1]\to W_{p-1}.\]
Indeed, using~\cref{BoundarySuspect} we see that the $i$-th face of $y$ is already in $W_{p-1}$ for $i\neq p$ since:
\begin{itemize}[leftmargin=*]
\renewcommand\labelitemii{$\diamondsuit$}
    \item if $0\leq i \leq p-2$, the face $d_iy$ is a suspect $d$-simplex, and in particular it belongs to $X_{d-1}\subseteq W_{p-1}$.
     \item if $i=p-1$, the face $d_iy$ has suspect index at most $p-1$, and in particular it belongs to $W_{p-1}$ (even in $X_{d-1}$ if $p=1$).
        \item if $p+1\leq i \leq r+1$, the face $d_iy$ is a suspect $d$-simplex, and in particular it belongs to $X_{d-1}\subseteq W_{p-1}$.
     \item if $r+1\leq i\leq d+1$, the face $d_iy$ is of rank $r+1$,
     and in particular it belongs to $Y_{r+1}\subseteq W_{p-1}$.
        \end{itemize}
    \item We argue that the $r$-th horn of $y$ defines a map of simplicial sets
    \[\Lambda^p[d+1]\to W_{p-1}\]
    with  marking.
To this end, we need to show that any face of $\Delta[d+1]$ containing
$\{p-1, p, p + 1\}$ is mapped to a marked simplex of $W_{p-1}$. Evoking~\cref{descriptionnerve}, it is enough to show that $y$ maps any basis element $[\mathbf{a}] = [\mathbf{a}_{<p-1}, p-1,p,p+1, \mathbf{a}_{>p+1}]$ of $O[d+1]$ to $0$. Using the closed formula for $w_p$ from~\cref{closedformulas} and the alternative description of degeneracies in~\cref{DefFaceDeg}, we obtain that
\[y[\mathbf{a}_{<p-1}, p-1,p,p+1, \mathbf{a}_{>p-1}]=(w_pd_py) [\mathbf{a}_{<p-1}, p-1,p,p+1, \mathbf{a}_{>p-1}]=0
  \]
as desired.
    \item If furthermore $d_py$ is marked, we argue that the $p$-th horn of $y$ defines a map of simplicial sets with  marking
        \[\Lambda^p[d+1]'\to W_{p-1},\]
        with the simplicial set with marking $\Lambda^p[d+1]'$ defined in~\cref{CompMarkAtOnce}.
        
         To this end, we need to show that the $(p-1)$-st and $(p+1)$-faces
  of the top $(d+1)$-dimensional
  simplex are mapped to a marked simplex of $W_{p-1}$.

As an auxiliary fact, we observe that
\[y[0,1,\ldots, d+1]=(w_pd_py) [0,1,\ldots, d+1]=0
  \]
  so the top $(d+1)$-dimensional simplex is mapped to a marked simplex of $W_{p-1}$.
     Moreover we have already shown that for all $j \notin \{p-1,p,p+1\}$ the face $d_jy$ is marked, namely $(d_jy)[0,1,2,\ldots, d]=0$, and we additionally assumed that $d_py$ is marked, namely $(d_py)[0,1,2,\ldots, d]=0$.
It then follows that
\begin{equation}
\label{equation}
\begin{array}{lll}
0&=&\partial y[0,1,2,\ldots, d+1]\\
&=&y\partial [0,1,2,\ldots, d+1]\\
&=&(-1)^{p+1}(d_{p+1}y)[0,1,2,\ldots, d] +\\
&&+ (-1)^{p-1}(d_{p-1}y)[0,1,2,\ldots, d].
\end{array}
\end{equation}
Since both faces $d_{p+1}y$ and $d_{p-1}y$ are directed chain maps, it must be that
\[
(d_{p+1}y)[0,1,2,\ldots, d]=(d_{p-1}y)[0,1,2,\ldots, d]=0.
\]
This means precisely that both faces $d_{p-1}y$ and $d_{p+1}y$ are marked, as desired.
\end{itemize}
 By filling all $p$-horns of suspect $(d+1)$-simplices $y$ of $W_p$, we then obtain their $p$-th face $x$, which was missing in $W_{p-1}$, as well as the suspect $(d+1)$-simplex $y$ itself.
This can be rephrased by saying that there is a pushout square 
\[
\begin{tikzcd}[column sep=1.0cm]
\underset{\mathclap{\begin{subarray}{c}
  x \\
  \mbox{\tiny{non-marked}}
  \end{subarray}} }{\coprod} \Lambda^{p}[d+1]\amalg\underset{\mathclap{\begin{subarray}{c}
 x \\
  \mbox{\tiny{marked}}
  \end{subarray} }}{\coprod} \Lambda^{p}[d+1]'\arrow[d]\arrow[r] &
 \underset{\mathclap{\begin{subarray}{c}
  x \\
  \mbox{\tiny{non-marked}}
  \end{subarray} }}{\coprod} \Delta^{p}[d+1]\amalg \underset{\mathclap{\begin{subarray}{c}
x\\
  \mbox{\tiny{marked}}
  \end{subarray}} }{\coprod} \Delta^{p}[d+1]''\arrow[d]\\
W_{p-1} \arrow[r]& W_p.
\end{tikzcd}
\]
Since the involved horn inclusions are in fact inner horn inclusions,
the inclusions of simplicial sets with marking $\Lambda^p[d+1]\hookrightarrow\Delta^p[d+1]$ and $\Lambda^p[d+1]'\hookrightarrow\Delta^p[d+1]''$ are complicial inner anodyne extensions by~\cref{CompMarkAtOnce}.

It follows that
the inclusion  $W_{{p-1}}\hookrightarrow W_p$ for any $1\le p\le r$, the inclusion $Y_{r+1}\hookrightarrow Y_r$ for any $1\le r\le d$,
the inclusion $X_{d-1}\hookrightarrow X_d$ for any $d\ge1$, and the inclusion $\NRS\cD\star\Delta[\top]\hookrightarrow \NRS(\cD\star\cO[\top])$ are complicial inner anodyne extensions, as desired.
\end{proof}

We conclude with two interesting consequences of the main theorem. The first one recovers a case treated by Maehara in \cite{MaeharaOrientals}.

\begin{cor}
\label{orientalsimplex}
For every $m\ge0$,
the canonical inclusion
\[\Delta[n]\hookrightarrow\NRS\cO[n]\]
from~\cref{comparisonmap} is a complicial inner anodyne extension, and in particular an acyclic cofibration in the model structure(s) from~\cref{modelstructurewithsaturation}.
\end{cor}

\begin{proof}
We prove the statement by induction on $n\ge0$. For $n=0$ the map in question is the isomorphism of simplicial sets with marking
\[\Delta[\top]\cong\NRS\cO[\top].\]
We now prove the statement for $m>0$. By induction hypothesis, the canonical map
\[\Delta[n-1]\hookrightarrow\NRS\cO[n-1]\]
is a complicial inner anodyne extension. Since \cite[Observation 40]{VerityComplicialI} implies that the join construction $\star\colon\msset\times\msset\to\msset$ preserves complicial inner anodyne extensions in each variable,
we then obtain that the canonical map
\[\Delta[n-1]\star\Delta[\top]\hookrightarrow\NRS\cO[n-1]\star\Delta[\top]\]
is a complicial inner anodyne extension.
By~\cref{maintheorem}, we then obtain that the canonical inclusion
\[\Delta[n]\cong\Delta[n-1]\star\Delta[\top]\hookrightarrow\NRS(\cO[n-1]\star\cO[\top])\cong\NRS\cO[n].\]
is a complicial inner anodyne extension,
as desired.
\end{proof}

\subsubsection{The statements for the dual case}

Finally, with techniques similar to those employed in \cref{maintheorem} and using all the adjusted preliminary results that we described along the way, one can prove the following theorem.

\label{maindual}
\begin{thm}
\label{maintheoremdual}
Let $\cD$ be a
strong Steiner $\omega$-category. The inclusion \[\Delta[\bot]\star\NRS\cD\hookrightarrow\NRS(\cO[\bot]\star\cD)\]
from~\cref{comparisonmap} is a complicial inner anodyne extension, and in particular an acyclic cofibration in the model structure(s) from~\cref{modelstructurewithsaturation}.
\end{thm}



 \bibliographystyle{amsalpha}
\bibliography{ref}

\end{document}